\documentclass[11pt]{article}
\usepackage[colorlinks,linkcolor=blue,citecolor=orange]{hyperref}

\usepackage[letterpaper]{geometry}
\usepackage{amsmath,amsthm,amsfonts,amssymb,booktabs}
\usepackage{enumerate,color,xcolor}
\usepackage{graphicx}
\usepackage{subfigure}
\usepackage{url}

\numberwithin{equation}{section}
\theoremstyle{plain}

\newtheorem{theorem}{Theorem}

\newtheorem{corollary}[theorem]{Corollary}

\newtheorem{example}[theorem]{Example}
\newtheorem{lemma}[theorem]{Lemma}

\newtheorem{proposition}[theorem]{Proposition}

\newtheorem{assumption}[theorem]{Assumption}

\usepackage{comment}

\newcommand{\R}{\mathbb{R}}
\newcommand{\bz}{\mathbf{z}}
\newcommand{\beq}{\begin{eqnarray}}
\newcommand{\eeq}{\end{eqnarray}}

\newcommand{\heps}{\hat{\varepsilon}}

\newcommand{\bigO}{\mathcal{O}}
\newcommand{\tildeO}{\tilde{\mathcal{O}}}
\newcommand{\tildeOm}{\tilde{\Omega}}
\newcommand{\E}{\mathbb{E}}
\newcommand{\bZ}{\mathbf{Z}}

\makeatletter
\def\munderbar#1{\underline{\sbox\tw@{$#1$}\dp\tw@\z@\box\tw@}}
\makeatother

\title{Global Convergence of Stochastic Gradient Hamiltonian Monte Carlo for Non-Convex Stochastic Optimization: \\Non-Asymptotic Performance Bounds and\\ Momentum-Based Acceleration}
\author{{Xuefeng
  Gao}\,\footnote{The authors are in alphabetical order.}~\footnote{Department of Systems
    Engineering and Engineering Management, The Chinese University of Hong Kong, Shatin, N.T. Hong Kong;
    xfgao@se.cuhk.edu.hk.}, 
    {Mert G\"{u}rb\"{u}zbalaban}\,$^*$\footnote{Department of Management Science
and Information Systems and the DIMACS Institute, Rutgers University, Piscataway, NJ-08854, United States of America;
    mg1366@rutgers.edu.},
  Lingjiong Zhu\,$^*$\footnote{Department of Mathematics, Florida State University, 1017 Academic Way, Tallahassee, FL-32306, United States of America; zhu@math.fsu.edu.
  }}
\date{}

\begin{document}
 
\maketitle
\vspace{-1pc}
\begin{center}
 	\today
 \end{center}
\vspace{1pc}
\begin{abstract} 
Stochastic gradient Hamiltonian Monte Carlo (SGHMC) is a variant of stochastic gradient with momentum  where a controlled and properly scaled Gaussian noise is added to the stochastic gradients to steer the iterates towards a global minimum. Many works reported its empirical success in practice for solving stochastic non-convex optimization problems, in particular it has been observed to outperform overdamped Langevin Monte Carlo-based methods such as stochastic gradient Langevin dynamics (SGLD) in many applications. Although asymptotic global convergence properties of SGHMC are well known, its finite-time performance is not well-understood. 
In this work, we study two variants of SGHMC based on two alternative discretizations of the underdamped Langevin diffusion. We provide finite-time performance bounds for the global convergence of both SGHMC variants for solving stochastic non-convex optimization problems with explicit constants. Our results lead to non-asymptotic guarantees for both population and empirical risk minimization problems. For a fixed target accuracy level, on a class of non-convex problems, we obtain complexity bounds for SGHMC that can be tighter than those available for SGLD.
\end{abstract}

\author{}

\section{Introduction}
We consider the stochastic non-convex optimization problem
\begin{equation}\label{opt-pbm}
\min_{x\in\mathbb{R}^{d}}F(x):=\mathbb{E}_{Z \sim \mathcal{D}}[f(x,Z)]\,,
\end{equation}
where $Z$ is a random variable whose probability distribution $\mathcal{D}$ is unknown, supported on some unknown set $\mathcal{Z}$, the objective $F$ is the expectation of a random function $f:\R^d \times \mathcal{Z}\to\R$ where the functions $x\mapsto f(x,z)$ are continuous and non-convex. Having access to independent and identically distributed (i.i.d.) samples $\bZ = (Z_1, Z_2, \dots, Z_n)$ where each $Z_i$ is a random variable distributed with the population distribution $\mathcal{D}$, the goal is to compute an approximate minimizer $\hat{x}$ (possibly with a randomized algorithm) of the \emph{population risk}, i.e. we want to compute $\hat{x}$ such that $\mathbb{E} F(\hat{x}) - F^* \leq \hat\varepsilon$ for a given target accuracy $\hat\varepsilon>0$, 
where $F^* = \min_{x\in\mathbb{R}^{d}} F(x)$ is the minimum value and the expectation is taken with respect to both $\bZ$ and the randomness encountered (if any) during the iterations of the algorithm to compute $\hat{x}$. This formulation arises frequently in several contexts including machine learning. A prominent example is deep learning
where $x$ denotes the set of trainable weights for a deep learning model and $f(x,z_i)$ is the penalty (loss) of prediction using weight $x$ with the individual sample value $Z_i = z_i \in \mathcal{Z}$.

Because the population distribution $\mathcal{D}$ is unknown, a common popular approach is to consider the \emph{empirical risk minimization} problem	\beq\label{pbm-finite-sum} \min_{x\in\mathbb{R}^d} F_\bz(x) := \frac{1}{n} \sum_{i=1}^n f(x,z_i) \,,
    \eeq
based on the dataset $\bz := (z_1, z_2, \dots, z_n) \in \mathcal{Z}^n$ as a proxy to the problem \eqref{opt-pbm} and minimize the \emph{empirical risk}
	\beq\label{def-emp-risk} \mathbb{E} F_\bz(x) - \min_{x\in\R^d} F_\bz(x)
    \eeq
instead, where the expectation is taken  with respect to any randomness encountered during the algorithm to generate $x$.\footnote{We note that in our notation $\bZ$ is a random vector, whereas $\bz$ is deterministic vector associated to a dataset that corresponds to a realization of the random vector $\bZ$.} Many algorithms have been proposed to solve the problem \eqref{opt-pbm} and its finite-sum version \eqref{pbm-finite-sum}. Among these, gradient descent, stochastic gradient and their variance-reduced or momentum-based variants come with guarantees for finding a local minimizer or a stationary point for non-convex problems. In some applications, convergence to a local minimum can be satisfactory (\cite{ge2017learning,du2017gradient}). However, in general, methods with global convergence guarantees are also desirable and preferable in many settings (\cite{hazan2016graduated,simsekli-async-MCMC-18}).

It has been well known that sampling from a distribution which concentrates around a global minimizer of $F$ is a similar goal to computing an approximate global minimizer of $F$. For example such connections arise in the study of simulated annealing algorithms in optimization which admit several asymptotic convergence guarantees (see e.g. \cite{gidas1985nonstationary, hajek1985tutorial,gelfand1991recursive,kirkpatrick1983optimization,bertsimas1993simulated,Belloni,Borkar}). Recent studies made such connections between the fields of statistics and optimization stronger, justifying and popularizing the use of Langevin Monte Carlo-based methods in stochastic non-convex optimization and large-scale data analysis further (see e.g. \cite{Chaudhari,Dalalyan, Raginsky,chen2016bridging,simsekli2016stochastic,simsekli-async-MCMC-18,welling2011bayesian,wibisono2018sampling}).

Stochastic gradient algorithms based on Langevin Monte Carlo are popular variants of stochastic gradient which admit asymptotic global convergence guarantees where a properly scaled Gaussian noise is added to the gradient estimate. Two popular Langevin-based algorithms that have demonstrated empirical success are stochastic gradient Langevin dynamics (SGLD) (\cite{welling2011bayesian, carin-2015-langevin-integrators}) and stochastic gradient Hamiltonian Monte Carlo (SGHMC) (\cite{emilyfox-sghmc,carin-2015-langevin-integrators,neal2010mcmc,duane1987hybrid}) and their variants to improve their efficiency and accuracy (\cite{ahn2012bayesian,ma2015complete,patterson-teh,ding2014bayesian,wibisono2018sampling}). In particular, SGLD can be viewed as the analogue of stochastic gradient in the Markov Chain Monte Carlo (MCMC) literature whereas SGHMC is the analogue of stochastic gradient with momentum (see e.g. \cite{emilyfox-sghmc}). SGLD iterations consist of
\begin{equation*}\label{eq-sgld}
X_{k+1}=X_{k}-\eta g_{k}+\sqrt{2\eta\beta^{-1}}\xi_{k}\,,
\end{equation*}
where $\eta>0$ is the stepsize parameter, $\beta>0$ is the inverse temperature, $g_k$ is a conditionally unbiased estimate of the gradient of $F_\bz$ and $\xi_k \in \R^d$ is a sequence of i.i.d. centered Gaussian random vector with unit covariance matrix. When the gradient variance is zero, SGLD dynamics corresponds to (explicit) Euler discretization of the first-order (a.k.a. overdamped) Langevin stochastic differential equation (SDE)
\begin{equation}\label{eq:overdamped}
dX(t)=-\nabla F_{\mathbf{z}}(X(t))dt+\sqrt{2\beta^{-1}}dB(t)\,, \quad t\geq 0\,,
\end{equation}
where $\{B(t): t \ge 0 \}$ is the standard Brownian motion in $\mathbb{R}^d$. The process $X$ admits a unique stationary distribution $\pi_\bz (dx) \propto \exp(-\beta F_\bz(x))dx$, also known as the \emph{Gibbs measure}, under some assumptions on $F_{\mathbf{z}}$ (see e.g. \cite{chiang1987diffusion,stroock-langevin-spectrum}). For $\beta$ chosen properly (large enough), it is easy to see that this distribution will concentrate around approximate global minimizers of $F_{\mathbf{z}}$. Recently, \cite{Dalalyan} established novel theoretical guarantees for the convergence of the overdamped Langevin MCMC and the SGLD algorithm for sampling from a smooth and log-concave density and these results have direct implications to stochastic convex optimization; see also \cite{DK2017}. In a seminal work, \cite{Raginsky} showed that SGLD iterates track the overdamped Langevin SDE closely and obtained finite-time performance bounds for SGLD. Their results show that SGLD converges to $\varepsilon$-approximate global minimizers after $\mathcal{O}(\text{poly}(\frac{1}{\lambda_*},\beta,d, \frac{1}{\varepsilon}))$ iterations where $\lambda_*$ is the uniform spectral gap that controls the convergence rate of the overdamped Langevin diffusion which is in general exponentially small in both $\beta$ and the dimension $d$ (\cite{Raginsky,pmlr-v75-tzen18a}).  A related result of
\cite{zhang-sgld} shows that a modified version of the SGLD algorithm will find an $\varepsilon$-approximate local minimum after polynomial time (with respect to all parameters). Recently, \cite{xu2018global} improved the $\varepsilon$ dependency of the upper bounds of \cite{Raginsky} further in the mini-batch setting, and obtained several guarantees for the gradient Langevin dynamics and variance-reduced SGLD algorithms. 

On the other hand, the SGHMC algorithm is based on the underdamped (a.k.a. second-order or kinetic) Langevin diffusion
\begin{align}
&dV(t)=-\gamma V(t)dt- \nabla F_{\mathbf{z}}(X(t))dt+\sqrt{2\gamma \beta^{-1}}dB(t), \label{eq:VL}
\\
&dX(t)=V(t)dt, \label{eq:XL}
\end{align}
where $\gamma>0$ is the friction coefficient, $X(t),V(t) \in \mathbb{R}^d$ models the position and the momentum of a particle moving in a field of force (described by the gradient of $F_{\mathbf{z}}$) plus a random (thermal) force described by Brownian noise, first derived by \cite{kramers1940brownian}. It is known that under some assumptions on $F_{\mathbf{z}}$, the Markov process $(X(t), V(t))_{t\geq 0}$ is ergodic and admits a unique stationary distribution
\begin{equation} \label{eq:inv-distr}
\pi_{\mathbf{z}}(dx, dv)= \frac{1}{\Gamma_{\mathbf{z}}} \exp\left(-\beta \left(\frac{1}{2} \Vert v\Vert^2 + F_{\mathbf{z}}(x)\right) \right) dx dv,
\end{equation}
(see e.g. \cite{herau-nier-underdamped,pavliotis2014stochastic}) where $\Gamma_{\mathbf{z}}$ is the normalizing constant:
\begin{equation*}
\Gamma_{\mathbf{z}} = \int_{\mathbb{R}^d \times \mathbb{R}^d} \exp\left(-\beta \left(\frac{1}{2} \Vert v\Vert^2 + F_{\mathbf{z}}(x)\right) \right) dx dv = \left(\frac{2 \pi}{\beta}\right)^{d/2} \int_{\mathbb{R}^d} e^{- \beta F_{\mathbf{z}}(x)} dx.
\end{equation*}
Hence, the $x$-marginal distribution of stationary distribution $\pi_{\mathbf{z}}(dx, dv)$ is exactly the invariant distribution of the overdamped Langevin diffusion.\footnote{With slight abuse of notation, we use $\pi_\bz(dx)$ to denote the $x$-marginal of the equilibrium distribution $\pi_\bz(dx,dv)$.} SGHMC dynamics correspond to the discretization of the underdamped Langevin SDE where the gradients are replaced with their unbiased estimates. Although various discretizations of the underdamped Langevin SDE has also been considered and studied (\cite{carin-2015-langevin-integrators,leimkuhler2015computation}), the following first-order Euler scheme is the simplest approach that is easy to implement, and a common scheme among the practitioners (\cite{teh2016consistency,chen2016bridging,carin-2015-langevin-integrators}):
\begin{align}
&V_{k+1}=V_{k}-\eta[\gamma V_{k}+g(X_{k},U_{{\mathbf{z}},k})]+\sqrt{2\gamma \beta^{-1} \eta}\xi_{k},\label{eq:V-iterate}
\\
&X_{k+1}=X_{k}+\eta V_{k}, \label{eq:X-iterate}
\end{align}
where $(\xi_k)_{k=0}^{\infty}$
is a sequence of i.i.d standard Gaussian random vectors in $\mathbb{R}^d$,
$\{U_{{\mathbf{z}},k}: k =0, 1, \ldots\}$ is a sequence of i.i.d random elements such that
\begin{equation*}
\mathbb{E}g(x,U_{{\mathbf{z}}, k}) = \nabla F_{\mathbf{z}}(x) \quad \text{for any $x\in\mathbb{R}^d$}.
\end{equation*}
In this paper, we focus on the unadjusted dynamics (without Metropolis-Hastings type of correction) that works well in many applications (\cite{emilyfox-sghmc,carin-2015-langevin-integrators}), as Metropolis-Hastings correction is typically computationally expensive for applications in machine learning and large-scale optimization when the size of the dataset $n$ is large and low to medium accuracy is enough in practice (see e.g. \cite{welling2011bayesian,emilyfox-sghmc}). 

There is also an alternative discretization to \eqref{eq:V-iterate}-\eqref{eq:X-iterate}, recently proposed by \cite{Cheng} which leads to state-of-the-art estimates in the special case that improves upon the Euler discretization when the objective is strongly convex (\cite{Cheng}). To introduce this alternative discretization by  \cite{Cheng}, 
we first define a sequence of functions $\psi_k$ by $\psi_{0}(t)=e^{-\gamma t}$ and $\psi_{k+1}(t)=\int_{0}^{t}\psi_{k}(s)ds$, $k\geq 0$. The iterates $(\hat{X}_{k}, \hat{V}_{k})$ are then defined by the following recursion:
\begin{align}
&\hat{V}_{k+1}=\psi_{0}(\eta)\hat{V}_{k}-\psi_{1}(\eta)g(\hat{X}_{k},U_{\mathbf{z},k})+\sqrt{2\gamma\beta^{-1}}\xi_{k+1}, \label{SGHMC2:V}
\\
&\hat{X}_{k+1}=\hat{X}_{k}+\psi_{1}(\eta)\hat{V}_{k}-\psi_{2}(\eta)g(\hat{X}_{k},U_{\mathbf{z},k})+\sqrt{2\gamma\beta^{-1}}\xi'_{k+1}, \label{SGHMC2:X}
\end{align}
where $(\xi_{k+1},\xi'_{k+1})$ is a $2d$-dimensional centered Gaussian vector so that
$(\xi_{j},\xi'_{j})$'s are independent and identically distributed (i.i.d.) and independent of the initial condition,
and for any fixed $j$, the random vectors $((\xi_{j})_{1},(\xi'_{j})_{1})$, $((\xi_{j})_{2},(\xi'_{j})_{2})$,
$\ldots$ $((\xi_{j})_{d},(\xi'_{j})_{d})$ are i.i.d. with the covariance matrix:
\begin{equation}
C(\eta)=\int_{0}^{\eta}\left[\psi_{0}(t),\psi_{1}(t)\right]^{T}
\left[\psi_{0}(t),\psi_{1}(t)\right]dt. 
\end{equation}
In the rest of the paper, we refer to Euler discretization \eqref{eq:V-iterate}-\eqref{eq:X-iterate} as SGHMC1 whereas the alternative discretization \eqref{SGHMC2:V}-\eqref{SGHMC2:X} as SGHMC2.

Recently, \cite{Eberle} show that the underdamped SDE converges to its stationary distribution faster than {that of the best known convergence rate} of overdamped SDE in the 2-Wasserstein metric under some assumptions, where $F_\bz$ can be non-convex. Their result is for the continuous-time underdamped dynamics. 
This raises the natural question whether the discretized underdamped dynamics (SGHMC), can lead to better guarantees than the SGLD method for solving stochastic non-convex optimization problems. Indeed, experimental results show that SGHMC can outperform SGLD dynamics in many applications (see e.g. \cite{Eberle,carin-2015-langevin-integrators,emilyfox-sghmc}). Although asymptotic convergence guarantees for SGHMC exist (see e.g. \cite{emilyfox-sghmc} \cite[Section 3]{mattingly2002ergodicity}, \cite{leimkuhler2015computation}), there is a lack of finite-time explicit performance bounds for solving  non-convex stochastic optimization problems with SGHMC in the literature including risk minimization problems.

\subsection{Contributions}
Our main contributions can be summarized as follows:
\begin{itemize}
    \item We provide for the first time the non-asymptotic provable guarantees for SGHMC to find approximate minimizers of both empirical and population risks with explicit constants. We establish the results under some regularity and growth assumptions for the component functions $f(x,z)$ and the noise in the gradients, but we do not assume $f$ is strongly convex in any region.
    
    \item We show that for a class of non-convex problems,
SGHMC2 can improve upon 
the (vanilla) SGLD algorithm in terms of the \emph{gradient complexity}, i.e. the total number of stochastic gradients required to achieve a global minimum.
Here, ``improvement'' means the best available bounds for
SGHMC2, which we prove in our paper, are better than the best available bounds for SGLD for some class of problems; see Section~\ref{sec-compare} for details. As a consequence, our analysis gives further theoretical justification to the success of momentum-based methods for solving non-convex machine learning problems, empirically observed in practice (see e.g. \cite{sutskever2013importance}).

\item We illustrate the applications of our theoretical results using two examples including binary linear classification and robust ridge regression. 

\item On the technical side,
we adapt the proof techniques of \cite{Raginsky} developed for the overdamped dynamics to the underdamped dynamics and combine it with the analysis of \cite{Eberle} which quantifies the convergence rate of the underdamped Langevin SDE to its equilibrium. The main new technical results we derive in this paper, relative to these studies, include controlling the discretization errors between SGHMC and
the continuous-time underdamped Langevin SDE, and bounding the moments of underdamped dynamics.

\end{itemize}

\subsection{Related Work and Comparison to Existing Literature}
In a recent work, \cite{simsekli-async-MCMC-18} obtained a finite-time performance bound for the ergodic average of the SGHMC iterates in the presence of delays in gradient computations. Their analysis highlights the dependency of the optimization error on the delay in the gradient computations and the stepsize explicitly, however it hides some implicit constants which can be exponential both in $\beta$ and $d$ in the worst case. A comparison with the SGLD algorithm is also not given. On the contrary, in our paper, we make all the constants explicit. This allows us to make gradient complexity comparisons with respect to overdamped MCMC approaches such as SGLD.

 \cite{cheng-nonconvex} considered the problem of sampling from a target distribution $p(x)\propto \exp(-F(x))$ where $F:\R^d \to \R$ is $L$-smooth everywhere and $m$-strongly convex outside a ball of finite radius $\mathcal{R}$. They proved upper bounds for the time required to sample from a distribution that is within $\varepsilon$ of the target distribution with respect to the $1$-Wasserstein distance for both underdamped and overdamped methods that scales polynomially in $\varepsilon$ and $d$. They also show that underdamped MCMC has a better dependency with respect to $\varepsilon$ and $d$ by a square root factor. Compared to this paper, in our analysis, we consider a larger class of non-convex functions $F(x)$ that satisfy the dissipativity condition, a weaker condition that does not require strong convexity outside a region. Under our assumptions, the best known bounds are such that the distance to the invariant distribution scales exponentially with dimension $d$ in the worst-case but not polynomially in $d$ (see e.g. \cite{Raginsky,xu2018global}). When $F$ is globally strongly convex (or equivalently when the target distribution $p(x) \propto \exp(-F(x))$ is strongly log-concave), there is also a growing interesting literature that establish performance bounds for both overdamped MCMC (see e.g. \cite{Dalalyan}) and underdamped MCMC methods (see e.g. \cite{Cheng,Mangoubi-Smith17}). In this particular setting, the fact that underdamped Langevin MCMC (also known as Hamiltonian MCMC) can improve upon the best available bounds for overdamped Langevin MCMC algorithms has also been proven (\cite{Cheng,Mangoubi-Smith17,dalalyan2018kinetic,chatterji2018theory}. Similar results have also been established when $F(x)$ is convex but not strongly convex (\cite{dalalyan2019}). Compared to these papers where $F(x)$ is convex, our assumptions are weaker as we allow $F(x)$ to be non-convex as long as it is dissipative.

A related paper \cite{xu2018global} applies variance reduction techniques to overdamped MCMC to improve performance when the empirical risk can be non-convex satisfying the same dissipativity assumption considered in our paper. However, these results do not give guarantees for the risk minimization problem \eqref{opt-pbm}. Furthermore, such variance reduction techniques require objectives in the form of a finite sum and do not apply to the \emph{streaming data setting} when each data point is used only once. In this work, we obtain guarantees for both the risk minimization problem and the empirical risk minimization and our results apply to the streaming data setting. Also, the convergence guarantees provided in \cite{xu2018global} depends on a spectral gap-related parameter that is not provided explicitly; whereas all our results are explicit and this allows us to have explicit performance comparisons between the upper bounds of SGLD and SGHMC algorithms.

We also note that underdamped Langevin MCMC (also known as Hamiltonian MCMC) and its practical applications have also been analyzed further in a number of recent works (see e.g. \cite{lee2018convergence,betancourt2017geometric,betancourt2017conceptual,betancourt2014optimizing,MPS2018}). In particular, \cite{MPS2018} provide a characterization of the conductance of Hamiltonian Monte Carlo (HMC) in continuous time using Liouville's theorem and invoking the Cheeger's inequality, they obtain upper and lower bounds on the spectral gap of HMC in continuous-time. Although the formula provided in \cite{MPS2018} for the conductance of HMC is elegant, it is not an explicit formula. In our analysis, our focus is to obtain performance bounds with explicit constants and therefore we build on the coupling techniques of \cite{Eberle} which leads to explicit constants for the class of problems we consider.

We also note that \cite{MPS2018} consider sampling from the target distribution $\frac{1}{2}\mathcal{N}(-1,\sigma^2) + \frac{1}{2}\mathcal{N}(1,\sigma^2) $ in dimension one and estimate the spectral gap of HMC in the regime as $\sigma \to 0$ . This is a mixture of two Gaussians with the same variance $\sigma^2$ centered at $-1$ and $1$ respectively where they argue that for this specific example HMC does not lead to much improvement over the Random Walk approach for sampling. In our paper, our results apply to more general targets that are not necessarily mixture of Gaussians. However, if we consider sampling from the distribution $\frac{1}{2}\mathcal{N}(-a,\sigma^2) + \frac{1}{2}\mathcal{N}(a,\sigma^2) $ as $a \to \infty$ for $\sigma^2$ fixed, Proposition \ref{prop:comparison} is applicable and it implies that HMC will be more efficient than overdamped Langevin dynamics in terms of dependency to $a$ (which measures the distance between the modes) in the sense that the mixing time will be $\mathcal{O}(a)$ in HMC whereas it will be $\mathcal{O}(a^2)$ in Random Walk. This does not contradict results of \cite{MPS2018} because we consider different scaling regimes: We fix $\sigma>0$ and let $a\to \infty$ whereas \cite{MPS2018} fix $a=1$ and let $\sigma\to 0$. 

There are also some connections of our work to existing momentum-based optimization algorithms. More specifically, if the term with $dB(t)$ involving the Brownian noise is removed in the underdamped SDE \eqref{eq:VL}--\eqref{eq:XL}, this results in a second-order ODE in $X(t)$. Momentum-based algorithms for strongly convex objectives such as Polyak's heavy ball method as well as Nesterov's accelerated gradient method can be both viewed as (alternative) discretizations of this ODE (see e.g. \cite{polyak1987introduction, su2014differential,shi2018understanding,wilson2016lyapunov}). It is known (\cite{su2014differential,shi2018understanding,wilson2016lyapunov}) that Nesterov's accelerated gradient method tracks this second-order ODE (also referred to as the Nesterov's ODE in the literature), whereas the first-order non-accelerated methods such as the classical gradient descent are known to track a first-order ODE in $X(t)$ called the \emph{gradient flow} dynamics. Furthermore, existing analysis shows that Nesterov's ODE converges to its equilibrium faster (in time) than the first-order gradient flow ODE in terms of upper bounds and this speed-up is also inherited by the discretized dynamics. Roughly speaking, our results can be interpreted as the analogue of these results in the non-convex optimization setting where we deal with SDEs instead of ODEs building on the theory of Markov processes and show that SGHMC tracks the second-order (underdamped) Langevin SDE closely and inherits its favorable convergence guarantees (in terms of upper bounds on the expected suboptimality) compared to that of overdamped Langevin SDE.

Acceleration of first-order gradient or stochastic gradient methods and their variance-reduced versions for finding a local stationary point (a point with a gradient less than $\varepsilon$ in norm) are also studied in the literature (see e.g. \cite{carmon18,nesterov1983method,ghadimi2016,JT2019,allen2016variance}). It has also been shown that under some assumptions momentum-based accelerated methods can escape saddle points faster (see e.g. \cite{wright2017escape,liu2018toward}). In contrast, in this work, our focus is obtaining performance guarantees for convergence to global minimizers instead. 

\section{Preliminaries and Assumptions}
In our analysis, we will use the following 2-Wasserstein distance: For any two probability measures $\nu_1,\nu_2$ on $\mathbb{R}^{2d}$, we define
\begin{equation*}
\mathcal{W}_2(\nu_1,\nu_2)=
\left(\inf_{Y_1 \sim\nu_{1}, Y_2 \sim\nu_{2}}\mathbb{E}\left[ \Vert Y_1 -Y_2 \Vert^2\right]\right)^{1/2},
\end{equation*}
where $\Vert \cdot \Vert$ is the usual Euclidean norm, $\nu_1, \nu_2$ are two Borel probability measures on $\mathbb{R}^{2d}$ with finite second moments, and the infimum is taken over all random couples $(Y_1, Y_2)$ taking values in $\mathbb{R}^{2d} \times \mathbb{R}^{2d}$ with marginals $Y_1 \sim\nu_1, Y_2 \sim\nu_2$ (see e.g. \cite{villani2008optimal}). 
We let $C^1(\mathbb{R}^d)$ denote the set of continuously differentiable functions on $\mathbb{R}^d$ and $L^2(\pi_{\bz})$ denote the space of square-integrable functions on $\mathbb{R}^d$ with respect to the measure $\pi_\bz$.

We first state the assumptions used in this paper below in Assumption~\ref{assumptions}. Note that we do not assume the component functions $f(x,z)$ to be convex; they can be non-convex. The first assumption of non-negativity of $f$ can be assumed without loss of generality by subtracting a constant and shifting the coordinate system as long as $f$ is bounded below. The second assumption of Lipschitz gradients is in general unavoidable for discretized Langevin algorithms to be convergent (see e.g. \cite{mattingly2002ergodicity}), and the third assumption is known as the \emph{dissipativity condition} (see e.g. \cite{hale1988asymptotic}) and is standard in the literature to ensure the convergence of Langevin diffusions to the stationary distribution (see e.g. \cite{Raginsky,Eberle,mattingly2002ergodicity}). The fourth assumption is regarding the amount of noise present in the gradient estimates and allows not only constant variance noise but allows the noise variance to grow with the norm of the iterates (which is the typical situation in mini-batch methods in stochastic gradient methods, see e.g. \cite{Raginsky}). Finally, the fifth assumption is a mild assumption saying that the initial distribution $\mu_0$ for the SGHMC dynamics should have a reasonable decay rate of the tails to ensure convergence to the stationary distribution. For instance, if the algorithm is started from any arbitrary point $(x_0,v_0) \in \R^{2d}$, then the Dirac measure $\mu_0(dx,dv) = \delta_{(x_0,v_0)}(dx,dv)$ would work. If the initial distribution $\mu_{0}(dx,dv)$ is supported
on a Euclidean ball with radius being some universal constant, it would also work. Similar assumptions on the initial distribution $\mu_0$ is also necessary to achieve convergence to a stationary measure in continuous-time underdamped dynamics as well (see e.g. \cite{herau-nier-underdamped}).



\begin{assumption}\label{assumptions}
We impose the following assumptions. 
\begin{itemize}
\item [$(i)$]
The function $f$ is continuously differentiable, takes non-negative real values, and there exist constants $A_{0},B\geq 0$ so that
\begin{equation*}
|f(0, z)|\leq A_{0},
\qquad
\Vert\nabla f(0,z)\Vert\leq B,
\end{equation*}
for any $z\in\mathcal{Z}$.
\item [$(ii)$]
For each $z\in\mathcal{Z}$, the function $f(\cdot,z)$ is $M$-smooth:
\begin{equation*}
\Vert\nabla f(w,z)-\nabla f(v,z)\Vert
\leq
M\Vert w-v\Vert.
\end{equation*}
\item [$(iii)$]
For each $z\in\mathcal{Z}$, the function $f(\cdot,z)$ is $(m,b)$-dissipative:
\begin{equation*}
\langle x,\nabla f(x,z)\rangle
\geq m\Vert x\Vert^{2}-b\,.
\end{equation*}
\item [$(iv)$]
There exists a constant $\delta\in[0,1)$ such that for every ${\mathbf{z}}$:
\begin{equation*}
\mathbb{E}[\Vert g(x,U_{\mathbf{z}})-\nabla F_{\mathbf{z}}(x)\Vert^{2}]
\leq
2\delta(M^{2}\Vert x\Vert^{2}+B^{2})\,.
\end{equation*}
\item [$(v)$]
The probability law $\mu_{0}$ of the initial state $(X_{0},V_{0})$ satisfies:
\begin{equation*}
\int_{\mathbb{R}^{2d}}e^{\alpha\mathcal{V}(x,v)}\mu_{0}(dx,dv)<\infty\,,
\end{equation*}
where $\mathcal{V}$ is a Lyapunov function to be used repeatedly for the rest of the paper:
\begin{equation} \label{eq:lyapunov}
\mathcal{V}(x,v):=\beta F_{\mathbf{z}}(x)
+\frac{\beta}{4}\gamma^{2}(\Vert x+\gamma^{-1}v\Vert^{2}+\Vert\gamma^{-1}v\Vert^{2}-\lambda\Vert x\Vert^{2})\,,
\end{equation}
and $\gamma$ is the friction coefficient as in \eqref{eq:VL},
$\lambda$ is a positive constant less than $\min(1/4,m/(M+\gamma^{2}/2))$, and $\alpha= \lambda(1- 2\lambda)/12$.
\end{itemize}
\end{assumption}

We note that the Lyapunov function $\mathcal{V}$ is used in \cite{Eberle} to study the rate of convergence to equilibrium for underdamped Langevin diffusion, 
which itself is motivated by e.g. \cite{mattingly2002ergodicity}. 
It follows from the above assumptions (applying Lemma~\ref{lem:gradient-bound}) that
there exists a constant $A\in(0,\infty)$ so that
\begin{equation}
x\cdot\nabla F_{\mathbf{z}}(x)\geq m\Vert x\Vert^{2}-b
\geq 2\lambda( F_{\mathbf{z}}(x)+ \gamma^{2}\Vert x\Vert^{2}/4)-2A/\beta\,.
\label{eq:drift}
\end{equation}
This drift condition, which will be used later, guarantees the stability and the existence of Lyapunov function $\mathcal{V}$ for the underdamped Langevin diffusion in \eqref{eq:VL}--\eqref{eq:XL}, see \cite{Eberle}. 

                            

\section{Main Results for SGHMC1 Algorithm} \label{sec:mainresult}

Our first result shows SGHMC1 iterates $(X_k,V_k)$ in \eqref{eq:V-iterate}--\eqref{eq:X-iterate} track the underdamped Langevin SDE in the sense that the expectation of the empirical risk $F_{\bz}$ with respect to the probability law of $(X_k,V_k)$ conditional on the sample $\bz$, denoted by $\mu_{k, \bz}$, and the stationary distribution $\pi_\bz$ of the underdamped SDE is small when $k$ is large enough. The difference in expectations decomposes as a sum of two terms $\mathcal{J}_0(\bz,\varepsilon)$ and $\mathcal{J}_1(\varepsilon)$ while the former term quantifies the dependency on the initialization and the dataset $\bz$ whereas the latter term is controlled by the discretization error and the amount of noise in the gradients which depends on the parameter $\delta$. We also note that the parameter $\mu_*$ (see Table~\ref{table_constants}) in our bounds governs the speed of convergence to the equilibrium of the underdamped Langevin diffusion.

\begin{theorem}\label{thm-first-term}Consider the SGHMC1 iterates $(X_k,V_k)$ defined by the recursion \eqref{eq:V-iterate}--\eqref{eq:X-iterate} from the initial state $(X_0, V_0)$ which has the law $\mu_0$. If Assumption~\ref{assumptions} is satisfied, then for $\beta,\varepsilon> 0$, we have
	\begin{align*}
\left|\mathbb{E}F_{\bz}(X_{k}) - \E_{(X,V) \sim \pi_\bz}(F_{\bz}(X)) \right| &= \left|\int_{\mathbb{R}^{d}\times\mathbb{R}^{d}} F_{\bz}(x) \mu_{k,\bz}(dx,dv) - \int_{\mathbb{R}^d\times\mathbb{R}^d} F_{\bz}(x) \pi_{\bz}(dx,dv) \right|
\\
&\leq \mathcal{J}_0(\bz,\varepsilon) + \mathcal{J}_1(\varepsilon)  \,,
\end{align*}
where
\begin{align}
&\mathcal{J}_0(\bz,\varepsilon) := (M \sigma +B) \cdot C\sqrt{\mathcal{H}_{\rho}(\mu_{0},\pi_{\mathbf{z}})}\cdot\varepsilon, \label{eq:J0z}
\\
&\mathcal{J}_1(\varepsilon) := (M \sigma +B) \cdot \Bigg( 
 \left(\frac{C_{0}}{\mu_{\ast}^{3/2}}(\log(1/\varepsilon))^{3/2}\delta^{1/4}
 +\frac{C_{1}}{\mu_{\ast}^{3/2}}\varepsilon
 \right)\sqrt{\log(\mu_{\ast}^{-1}\log(\varepsilon^{-1}))}+\frac{C_{2}}{\mu_{\ast}}\frac{\varepsilon^{2}}{(\log(1/\varepsilon))^{2}}\Bigg),\label{eq:J1}
\end{align}
with $\sigma$ defined by \eqref{defn:sigma} provided that
\begin{equation}\label{bound:eta}
\eta\leq\min\left\{\left(\frac{\varepsilon}{(\log(1/\varepsilon))^{3/2}}\right)^{4}
,1,\frac{\gamma}{K_{2}}(d/\beta+A/\beta),\frac{\gamma\lambda}{2K_{1}},\frac{2}{\gamma\lambda}\right\},
\end{equation}
\noindent and
\begin{equation}\label{bound:keta}
k\eta=\frac{1}{\mu_*}\log\left(\frac{1}{\varepsilon}\right)\geq e.
\end{equation}
Here $\mathcal{H}_{\rho}$ is a semi-metric for probability distributions defined by \eqref{def-H-rho}. All the constants are made explicit and are summarized in Table~\ref{table_constants}.
\end{theorem}  

The proof of Theorem~\ref{thm-first-term} will be presented in details in Section~\ref{sec:proof} in the Appendix. In the following subsections, we discuss that this theorem combined with some basic properties of the equilibrium distribution $\pi_\bz$ leads to a number of results which provide performance guarantees for both the empirical risk and population risk minimization.
\subsection{Performance bound for the empirical risk minimization}\label{sec:per-ERM}
In order to obtain guarantees for the empirical risk given in \eqref{def-emp-risk}, in light of Theorem~\ref{thm-first-term}, one has to control the quantity
    \begin{equation*}\label{ineq-almost-erm} \int_{\mathbb{R}^d\times\mathbb{R}^d} F_{\bz}(x) \pi_{\bz}(dx,dv) - \min_{x\in\R^d} F_{\bz}(x) \,,
    \end{equation*}
which is a measure of how much the $x-$marginal of the equilibrium distribution $\pi_\bz$ concentrates around a global minimizer of the empirical risk. As $\beta$ goes to infinity, it can be verified that this quantity goes to zero. For finite $\beta$, \cite{Raginsky} (see Proposition 11) derives an explicit bound of the form
	 \beq \int_{\mathbb{R}^d\times\mathbb{R}^d} F_{\bz}(x) \pi_{\bz}(dx,dv) - \min_{x\in\R^d} F_{\bz}(x)  \leq \mathcal{J}_2 := \frac{d}{2\beta}\log\left(\frac{eM}{m}\left(\frac{b\beta}{d}+1\right)\right),\label{def-J2}
    \eeq
(which is also provided in the Appendix for the sake of completeness, see Lemma~\ref{lem:thirdbound}). This combined with Theorem~\ref{thm-first-term} immediately leads to the following performance bound for the empirical risk minimization. The proof is omitted. 

\begin{corollary}[Empirical risk minimization with SGHMC1]\label{coro:main-exp} Under the setting of Theorem~\ref{thm-first-term}, the empirical risk minimization problem admits the performance bounds:
\begin{align*}
\mathbb{E}F_{\bz}(X_{k})-\min_{x\in\mathbb{R}^{d}}F_{\bz}(x)
\leq \mathcal{J}_0(\varepsilon,\bz) + \mathcal{J}_1(\varepsilon) + \mathcal{J}_2 \,,
\end{align*}
provided that conditions \eqref{bound:eta} and \eqref{bound:keta} hold where the terms $\mathcal{J}_0(\bz,\varepsilon)$, $\mathcal{J}_1(\varepsilon)$ and $\mathcal{J}_2$ are defined by $\eqref{eq:J0z}$, $\eqref{eq:J1}$ and $\eqref{def-J2}$ respectively.
\end{corollary}

\begin{table}{\emph{Constants} \hfill \emph{Source}}

\centering

\rule{\textwidth}{\heavyrulewidth}

\vspace{-1.6\baselineskip}
\begin{flalign*}
&
C_{x}^{c}=\frac{\int_{\mathbb{R}^{2d}}\mathcal{V}(x,v)\mu_{0}(dx,dv)
+\frac{(d+A)}{\lambda}}{\frac{1}{8} (1-2 \lambda) \beta \gamma^2 },
\quad
C_{v}^{c}=\frac{\int_{\mathbb{R}^{2d}}\mathcal{V}(x,v)\mu_{0}(dx,dv)
+\frac{(d+A)}{\lambda}}{\frac{\beta}{4}(1-2\lambda)}
&&&\eqref{def-Cx-c},\eqref{def-Cxv-c}
\end{flalign*}

\vspace{-1.6\baselineskip}

\begin{flalign*}
&
K_{1}
=\max\left\{\frac{32M^{2}\left(\frac{1}{2}+\gamma+\delta\right)}{(1-2 \lambda) \beta \gamma^2 },
\frac{8\left(\frac{1}{2}M+\frac{1}{4}\gamma^{2}-\frac{1}{4}\gamma^{2}\lambda+\gamma\right)}{\beta(1-2\lambda)}\right\}
&&&\eqref{def-K1}
\end{flalign*}

\vspace{-1.6\baselineskip}
\begin{flalign*}
&
K_{2}=B^{2}\left(1+2\gamma+2\delta\right)
&&&\eqref{def-K2}
\end{flalign*}

\vspace{-1.6\baselineskip}
\begin{flalign*}
&
C_{x}^{d}=\frac{\int_{\mathbb{R}^{2d}}\mathcal{V}(x,v)\mu_{0}(dx,dv)
+\frac{4(d+A)}{\lambda}}{\frac{1}{8} (1-2 \lambda) \beta \gamma^2 },
\quad
C_{v}^{d}=\frac{\int_{\mathbb{R}^{2d}}\mathcal{V}(x,v)\mu_{0}(dx,dv)
+\frac{4(d+A)}{\lambda}}{\frac{\beta}{4}(1-2\lambda)}
&&\eqref{def-Cx-d},\eqref{def-Cxv-d}
\end{flalign*}

\vspace{-1.6\baselineskip}
\begin{flalign*}
&
\sigma=\max\left\{\sqrt{C_{x}^{c}},\sqrt{C_{x}^{d}}\right\}= \sqrt{C_{x}^{d}}
&&&\eqref{defn:sigma}
\end{flalign*}

\vspace{-1.6\baselineskip}
\begin{flalign*}
& C_{0}=\hat{\gamma}
\cdot\left(\left(M^{2}C_{x}^{d}+B^{2}\right)\beta/\gamma+\sqrt{\left(M^{2}C_{x}^{d}+B^{2}\right)\beta/\gamma}\right)^{1/2}
&&&\eqref{def-hat-Czero}
\end{flalign*}

\vspace{-1.6\baselineskip}
\begin{flalign*}
&C_{1}=\hat{\gamma}
\cdot\left(\beta M^{2}(C_{2})^{2}/(2\gamma)+\sqrt{\beta M^{2}(C_{2})^{2}/(2\gamma)}\right)^{1/2}
&&&\eqref{def-hat-Cone}
\end{flalign*}

\vspace{-1.6\baselineskip}
\begin{flalign*}
&C_{2}=\left(2\gamma^{2}C_{v}^{d}
+(4+2\delta)\left(M^{2}C_{x}^{d}+B^{2}\right)
+2\gamma\beta^{-1}\right)^{1/2}
&&&\eqref{def-hat-Ctwo}
\end{flalign*}

\vspace{-1.6\baselineskip}
\begin{flalign*}
&\hat{\gamma}=\frac{2\sqrt{2}}{\sqrt{\alpha_{0}}}\left(\frac{5}{2}+\log\left(\int_{\mathbb{R}^{2d}}e^{\frac{1}{4}\alpha\mathcal{V}(x,v)}\mu_{0}(dx,dv)
+\frac{1}{4}e^{\frac{\alpha(d+A)}{3\lambda}}\alpha\gamma(d+A)\right)
\right)^{1/2}
&&&\eqref{def-hat-gamma}
\end{flalign*}

\vspace{-1.6\baselineskip}
\begin{flalign*}
&\alpha_{0}=\frac{\alpha(1-2\lambda)\beta\gamma^{2}}{64+32\gamma^{2}},
\qquad
\alpha=\frac{\lambda(1-2\lambda)}{12}
&&&\eqref{eq:alpha0alpha}
\end{flalign*}

\vspace{-1.6\baselineskip}
\begin{flalign*}
&
\mu_* =\frac{\gamma}{768}\min\{\lambda M \gamma^{-2},\Lambda^{1/2}e^{-\Lambda}M\gamma^{-2},\Lambda^{1/2}e^{-\Lambda}\}
&&&\eqref{def-mu-star}
\end{flalign*}

\vspace{-1.6\baselineskip}
\begin{flalign*}
&
C=\frac{(1+\gamma)\sqrt{2}e^{1+\frac{\Lambda}{2}}}{\min\{1,\alpha_{1}\}}
\sqrt{\max\{1,4(1+2\alpha_{1}+2\alpha_{1}^{2})(d+A)\beta^{-1}\gamma^{-1}\mu_*^{-1}/\min\{1,R_{1}\}\}}
&&&\eqref{def-capital-C}
\end{flalign*}

\vspace{-1.6\baselineskip}
\begin{flalign*}
&
\Lambda=\frac{12}{5}(1+2\alpha_{1}+2\alpha_{1}^{2})(d+A)M\gamma^{-2}\lambda^{-1}(1-2\lambda)^{-1},
\qquad
\alpha_{1}=(1+\Lambda^{-1})M\gamma^{-2}
&&&\eqref{def-capital-Lambda-alpha_1}
\end{flalign*}

\vspace{-1.6\baselineskip}
\begin{flalign*}
&
\varepsilon_{1}=4\gamma^{-1}\mu_{\ast}/(d+A)
&&&\eqref{def-eps_1}
\end{flalign*}

\vspace{-1.6\baselineskip}
\begin{flalign*}
&
R_{1}=4\cdot(6/5)^{1/2}(1+2\alpha_{1}+2\alpha_{1}^{2})^{1/2}(d+A)^{1/2} \beta^{-1/2}\gamma^{-1}(\lambda-2\lambda^{2})^{-1/2}
&&&\eqref{def-R1}
\end{flalign*}

\vspace{-1.6\baselineskip}
\begin{flalign*}
&
\overline{\mathcal{H}}_{\rho}(\mu_{0})=
R_{1}+R_{1}\varepsilon_{1}\max\left\{M+\frac{1}{2}\beta\gamma^{2},\frac{3}{4} \beta \right\}\Vert(x,v)\Vert_{L^{2}(\mu_{0})}^{2}
\\
&\qquad
+R_{1}\varepsilon_{1}\left(M+\frac{1}{2} \beta \gamma^{2}\right) \frac{b+d/\beta}{m}
+R_{1}\varepsilon_{1}\frac{3}{4} d
+2R_{1}\varepsilon_{1}\left(\beta A_{0}+\frac{\beta B^{2}}{2M}\right)
&&&\eqref{ineq-H-rho-upper}
\end{flalign*}

\rule{\textwidth}{\heavyrulewidth}

\caption{Summary of the constants and where they are defined in the text. }\label{table_constants}
\end{table}
\subsection{Performance bound for the population risk minimization}\label{sec:per-popu-risk}
By exploiting the fact that the $x-$marginal of the invariant distribution for the underdamped dynamics is the same as it is in the overdamped case, it can be shown that the generalization error $F(X_k)-F_{\bZ}(X_k)$ is no worse than that of the available bounds for SGLD given in \cite{Raginsky}, and therefore, we have the following corollary. A more detailed proof will be given in Section~\ref{sec:proof} in the Appendix.

\begin{corollary}[Population risk minimization with SGHMC1]\label{coro:main} Under the setting of Theorem~\ref{thm-first-term}, the expected population risk of $X_k$ (the iterates in \eqref{eq:X-iterate}) is bounded by
\begin{align*}
&\mathbb{E}F(X_k)-F^{\ast} \leq \overline{\mathcal{J}}_0(\varepsilon) + \mathcal{J}_1(\varepsilon) + \mathcal{J}_2 + \mathcal{J}_3(n) \,,
\end{align*}
with
\begin{align}
&\overline{\mathcal{J}}_0(\varepsilon) :=
(M \sigma +B)\cdot C\cdot\sqrt{\overline{\mathcal{H}}_{\rho}(\mu_{0})}\cdot\varepsilon, \label{def-J0-bar}\\
&\mathcal{J}_3(n) := \frac{4\beta c_{LS}}{n}
\left(\frac{M^{2}}{m}(b+d/\beta)+B^{2}\right),\label{def-J3}
\end{align}
where $\sigma$ is defined by \eqref{defn:sigma}, $\overline{\mathcal{H}}_{\rho}(\mu_{0})$ is defined by \eqref{ineq-H-rho-upper}, $\mathcal{J}_1(\varepsilon)$ and $\mathcal{J}_2$ are defined by \eqref{eq:J1} and \eqref{def-J2} respectively and $c_{LS}$ is a constant satisfying
\begin{equation*}
c_{LS}
\leq\frac{2m^{2}+8M^{2}}{m^{2}M\beta}
+\frac{1}{\lambda_{\ast}}\left(\frac{6M(d+\beta)}{m}+2\right),
\end{equation*}
and $\lambda_{\ast}$ is the uniform spectral gap for overdamped Langevin dynamics
\footnote{In \cite{Raginsky}, their formula for $\lambda_{\ast}$ missed $\beta^{-1}$ factor.}:
\begin{equation} \label{eq:lambda_star}
\lambda_{\ast}:=
\inf_{\bz\in\mathcal{Z}^n}
\inf
\left\{
\frac{\beta^{-1}\int_{\mathbb{R}^{d}}\Vert\nabla g\Vert^{2}d\pi_{\mathbf{z}}}{\int_{\mathbb{R}^{d}}g^{2}d\pi_{\mathbf{z}}}:
g\in C^{1}(\mathbb{R}^{d})\cap L^{2}(\pi_{\mathbf{z}}),
g\neq 0,
\int_{\mathbb{R}^{d}}gd\pi_{\mathbf{z}}=0\right\}.
\end{equation}
\end{corollary}


\subsection{Generalization error of SGHMC1 in the one pass regime} \label{sec:gene-error-1}
Since the predictor $X_k$ is random, it is natural to consider the expected generalization error $\E F(X_k) - \E F_{\bZ}(X_k)$ (see e.g. \cite{hardt2016train}) which admits the decomposition
\begin{align} \E F_{\bZ}(X_k) - \E F(X_k) =& \left(\E F_{\bZ}(X_k) - \E F_\bZ(X^{\pi})\right) + \left(\E F_\bZ(X^{\pi}) - \E F(X^{\pi})\right) \label{eq-general-part1} \\
&\quad + \left(\E F(X^{\pi}) - \E F(X_k)\right)\,, \nonumber  
\end{align}
where $X^\pi$ is the Gibbs output, i.e. its distribution conditional on $\bZ = \bz$ is given by $\pi_\bz$. If every sample is used once, i.e. if only one pass is made over the dataset, then the second term in \eqref{eq-general-part1} disappears. As a consequence, the generalization error is controlled by the bound
\beq | \E F_{\bZ}(X_k) - \E F(X_k)|  \leq \left |\E F_{\bZ}(X_k) - \E F_\bZ(X^{\pi})\right | + \left|\E F(X^{\pi}) - \E F(X_k)\right | \label{ineq-gen-error}.
\eeq
The following result provides a bound on this quantity. The proof is similar to the proof of Theorem~\ref{thm-first-term} and its corollaries, and hence omitted.

\begin{theorem}[Generalization error of SGHMC1]\label{thm-gen-error-sghmc} 
Under the setting of Theorem~\ref{thm-first-term}, we have
	\begin{align*}
 \left| \E F(X_k)-\E F(X^{\pi})\right|
\leq \overline{\mathcal{J}}_0(\varepsilon) + \mathcal{J}_1(\varepsilon)\,,\\
 \left| \E F_\bZ(X_k)-\E F_\bZ(X^{\pi})\right|
 \leq \overline{\mathcal{J}}_0(\varepsilon) + \mathcal{J}_1(\varepsilon),
\end{align*}
provided that \eqref{bound:eta} and \eqref{bound:keta} hold where $X^\pi$ is the output of the underdamped Langevin dynamics, i.e. its distribution conditional on $\bZ = \bz$ is given by $\pi_\bz$ and $ \overline{\mathcal{J}}_0(\varepsilon)$ is defined by \eqref{def-J0-bar}. Then, it follows from \eqref{ineq-gen-error} that if each data point is used once, the expected generalization error satisfies
 \begin{equation*} | \E F_{\bZ}(X_k) - \E F(X_k)|  \leq 2 \overline{\mathcal{J}}_0(\varepsilon) + 2 \mathcal{J}_1(\varepsilon). \end{equation*}
\end{theorem}

\section{Main Results for SGHMC2 Algorithm}\label{sec-sghmc2}

Recall the SGHMC2 algorithm $(\hat{X}_{k},\hat{V}_{k})$ defined
in \eqref{SGHMC2:V}-\eqref{SGHMC2:X}, and denote the probability law of $(\hat{X}_{k},\hat{V}_{k})$ conditional on the sample $\bz$ by $\hat{\mu}_{k,\bz}(dx,dv)$.
Similar to our analysis for SGHMC1, we can derive similar performance guarantees for SGHMC2 in terms of empirical risk, population risk and the generalization error. The main difference is that the term $\mathcal{J}_1(\varepsilon)$ is controlled by the accuracy of the discretization and has to be replaced by another term ${\hat{\mathcal{J}}}_1(\varepsilon)$, as SGHMC2 algorithm is based on an alternative discretization. In particular, the performance bounds we get for SGHMC2 are tighter than SGHMC1, as will be elaborated further in the Section~\ref{sec-compare}.  

\begin{theorem}\label{thm-first-term-2}
Consider the SGHMC2 iterates $(\hat{X}_{k}, \hat{V}_{k})$ defined by the recursion \eqref{SGHMC2:V}--\eqref{SGHMC2:X} from the initial state $(X_0, V_0)$ which has the law $\mu_0$. If Assumption~\ref{assumptions} is satisfied, then for $\beta,\varepsilon> 0$, we have
\begin{align*}
\left|\mathbb{E}F_{\bz}(\hat{X}_{k}) - \E_{(X,V) \sim \pi_\bz}(F_{\bz}(X)) \right| &= \left|\int_{\mathbb{R}^{d}\times\mathbb{R}^{d}} F_{\bz}(x) \hat{\mu}_{k,\bz}(dx,dv) - \int_{\mathbb{R}^d\times\mathbb{R}^d} F_{\bz}(x) \pi_{\bz}(dx,dv) \right|
\\
&\leq \mathcal{J}_0(\bz,\varepsilon) + \hat{\mathcal{J}}_1(\varepsilon)  \,,
\end{align*}
where $\mathcal{J}_{0}(\bz,\varepsilon)$ is defined in \eqref{eq:J0z} and
\begin{equation}\label{eq:J1:Jordan}
\hat{\mathcal{J}}_1(\varepsilon) := (M \sigma +B) \cdot \left(\frac{C_{0}}{\sqrt{\mu_{\ast}}}\sqrt{\log(1/\varepsilon)}\delta^{1/4}
 +\frac{\hat{C}_{1}}{\sqrt{\mu_{\ast}}}\varepsilon\right)\sqrt{\log(\mu_{\ast}^{-1}\log(\varepsilon^{-1}))}, 
\end{equation}
with $\sigma$ defined by \eqref{defn:sigma} provided that
\begin{equation}\label{bound:eta:Jordan}
\eta\leq\min\left\{\left(\frac{\varepsilon}{\sqrt{\log(1/\varepsilon)}}\right)^{2}
,1,\frac{\gamma}{\hat{K}_{2}}(d/\beta+A/\beta),\frac{\gamma\lambda}{2\hat{K}_{1}},\frac{2}{\gamma\lambda}\right\},
\end{equation}
\noindent and
\begin{equation}\label{bound:keta:Jordan}
k\eta=\frac{1}{\mu_*}\log\left(\frac{1}{\varepsilon}\right)\geq e.
\end{equation}
Here $\mathcal{H}_{\rho}$ is a semi-metric for probability distributions defined by \eqref{def-H-rho}. All the constants are made explicit and are summarized in Table~\ref{table_constants} and Table~\ref{table_constants-2}.
\end{theorem}


\begin{table}{\emph{Constants} \hfill \emph{Source}}

\centering

\rule{\textwidth}{\heavyrulewidth}

\vspace{-1.6\baselineskip}

\begin{flalign*}
&
\hat{K}_{1}
=K_{1}+Q_{1}\frac{4}{1-2\lambda}+Q_{2}\frac{8}{(1-2\lambda)\gamma^{2}}
&&&\eqref{def-K1-hat}
\end{flalign*}

\vspace{-1.6\baselineskip}
\begin{flalign*}
&
\hat{K}_{2}=K_{2}+Q_{3}
&&&\eqref{def-K2-hat}
\end{flalign*}

\vspace{-1.6\baselineskip}
\begin{flalign*}
&
Q_{1}=\frac{1}{2}c_{0}
\Bigg((5M+4-2\gamma+(c_{0}+\gamma^{2}))
+(1+\gamma)\left(\frac{5}{2}+c_{0}(1+\gamma)\right)
+2\gamma^{2}\lambda\Bigg)
&&&\eqref{def-Q1}
\end{flalign*}

\vspace{-1.6\baselineskip}
\begin{flalign*}
&
Q_{2}=\frac{1}{2}c_{0}
\Bigg[\Bigg((1+\gamma)\left(c_{0}(1+\gamma)+\frac{5}{2}\right)
+c_{0}+2
+\lambda\gamma^{2}
+2(Mc_{0}+M+1)\Bigg)\left(2(1+\delta)M^{2}\right)\nonumber
\\
&\qquad\qquad\qquad
+\left(2M^{2}+\gamma^{2}\lambda+\frac{3}{2}\gamma^{2}(1+\gamma)\right)\Bigg]
&&&\eqref{def-Q2}
\end{flalign*}

\vspace{-1.6\baselineskip}
\begin{flalign*}
&
Q_{3}=c_{0}
\Bigg((1+\gamma)\left(c_{0}(1+\gamma)+\frac{5}{2}\right)
+c_{0}+2
+\lambda\gamma^{2}
+2(Mc_{0}+M+1)\Bigg)(1+\delta)B^{2}+c_{0}B^{2}\nonumber
\\
&\qquad
+\frac{1}{2}\gamma^{3}\beta^{-1}c_{22}
+\gamma^{2}\beta^{-1}c_{12}+M\gamma\beta^{-1}c_{22}
&&&\eqref{def-Q3}
\end{flalign*}

\vspace{-1.6\baselineskip}
\begin{flalign*}
&
c_{0}=1+\gamma^{2},\qquad c_{12}=\frac{d}{2},\qquad c_{22}=\frac{d}{3}
&&&\eqref{def-c0-c12-c22}
\end{flalign*}

\vspace{-1.6\baselineskip}
\begin{flalign*}
&\hat{C}_{1}=\hat{\gamma}
\cdot\Bigg(\frac{3\beta M^{2}}{2\gamma}
\bigg(C_{v}^{d}+\left(2(1+\delta)M^{2}C_{x}^{d}+2(1+\delta)B^{2}\right)+\frac{2d\gamma\beta^{-1}}{3}\bigg)
\nonumber
\\
&\qquad\qquad\qquad
+\sqrt{\frac{3\beta M^{2}}{2\gamma}
\bigg(C_{v}^{d}+\left(2(1+\delta)M^{2}C_{x}^{d}+2(1+\delta)B^{2}\right)+\frac{2d\gamma\beta^{-1}}{3}\bigg)}\Bigg)^{1/2}
&&&\eqref{def-hat-Jone-2}
\end{flalign*}

\rule{\textwidth}{\heavyrulewidth}

\caption{Summary of the constants and where they are defined in the text. }\label{table_constants-2}
\end{table}

The proof of Theorem~\ref{thm-first-term-2} is given in Section~\ref{sec:proof:2} in the Appendix. 
Relying on Theorem~\ref{thm-first-term-2}, one can readily derive the following result on the
performance bound for the empirical risk minimization with the SGHMC2 algorithm. The proof follows a similar argument as discussed in Section~\ref{sec:per-ERM}, and is omitted.

\begin{corollary}[Empirical risk minimization with SGHMC2]\label{coro:main-exp-2} Under the setting of Theorem~\ref{thm-first-term-2}, the empirical risk minimization problem admits the performance bounds:
\begin{align*}
\mathbb{E}F_{\bz}(\hat X_{k})-\min_{x\in\mathbb{R}^{d}}F_{\bz}(x)
\leq \mathcal{J}_0(\bz,\varepsilon) + \hat{\mathcal{J}}_1(\varepsilon) + \mathcal{J}_2 \,,
\end{align*}
provided that conditions \eqref{bound:eta:Jordan} and \eqref{bound:keta:Jordan} hold where the terms $\mathcal{J}_0(\bz,\varepsilon)$, $\hat{\mathcal{J}}_1(\varepsilon)$ and $\mathcal{J}_2$ are defined by $\eqref{eq:J0z}$, $\eqref{eq:J1:Jordan}$ and $\eqref{def-J2}$ respectively.
\end{corollary}

Next, we present the performance bound for the population risk minimization with the SGHMC2 algorithm. Similar as in Section~\ref{sec:per-popu-risk}, to control the population risk during SGHMC2 iterations, one needs to control the difference between the finite sample size problem \eqref{pbm-finite-sum} and the original problem \eqref{opt-pbm} in addition to the empirical risk. This leads to the following result. The details of the proof are given in Section~\ref{sec:proof:2} in the Appendix.

\begin{corollary}[Population risk minimization with SGHMC2]\label{coro:main-2} 
Under the setting of Theorem~\ref{thm-first-term-2}, the expected population risk of $\hat{X}_k$ (the iterates in \eqref{SGHMC2:X}) is bounded by
\begin{align*}
&\mathbb{E}F(\hat X_k)-F^{\ast} \leq \overline{\mathcal{J}}_0(\varepsilon) + \hat{\mathcal{J}}_1(\varepsilon) + \mathcal{J}_2 + \mathcal{J}_3(n) \,,
\end{align*}
where $\overline{\mathcal{J}}_0(\varepsilon)$, $\hat{\mathcal{J}}_1(\varepsilon)$, $\mathcal{J}_2$, $\mathcal{J}_3(n)$ are
defined in \eqref{def-J0-bar}, \eqref{eq:J1:Jordan}, \eqref{def-J2} and \eqref{def-J3}.
\end{corollary}

Finally, we present a result on the generalization error of the SGHMC2 algorithm in the one pass regime. The proof follows from Theorem~\ref{thm-first-term-2} and the discussion for SGHMC1 algorithm in Section~\ref{sec:gene-error-1}, and hence is omitted.

\begin{theorem}[Generalization error of SGHMC2]\label{thm-gen-error-sghmc-2} 
Under the setting of Theorem~\ref{thm-first-term-2}, we have
	\begin{align*}
 \left| \E F(\hat{X}_k)-\E F(X^{\pi})\right|
\leq \overline{\mathcal{J}}_{0}(\varepsilon)+\hat{\mathcal{J}}_{1}(\varepsilon)\,,\\
 \left| \E F_\bZ(\hat{X}_k)-\E F_\bZ(X^{\pi})\right|
 \leq \overline{\mathcal{J}}_{0}(\varepsilon)+\hat{\mathcal{J}}_{1}(\varepsilon),
\end{align*}
provided that \eqref{bound:eta:Jordan} and \eqref{bound:keta:Jordan} hold where $X^\pi$ is the output of the underdamped Langevin dynamics, i.e. its distribution conditional on $\bZ = \bz$ is given by $\pi_\bz$ and $ \overline{\mathcal{J}}_0(\varepsilon)$ is defined by \eqref{def-J0-bar}. Then, it follows from \eqref{ineq-gen-error} that if each data point is used once, the expected generalization error satisfies
 \begin{equation*} | \E F_{\bZ}(\hat X_k) - \E F(\hat X_k)|  \leq 2 \overline{\mathcal{J}}_0(\varepsilon) + 2 \hat{\mathcal{J}}_1(\varepsilon). 
 \end{equation*} 
\end{theorem}

\section{Performance comparison with respect to SGLD algorithm}\label{sec-compare}
In this section, we compare our performance bounds for SGHMC1 and SGHMC2 to SGLD. 
We use the notations $\tilde{\bigO}$ and $\tildeOm$ to give explicit dependence on the parameters $d, \beta, \mu_*$ but it hides factors that depend (at worst polynomially) on other parameters $m,M,B, \lambda, \gamma, b$ and $A_{0}$. Without loss of generality, we assume here the initial distribution $\mu_{0}(dx,dv)$ is supported
on a Euclidean ball with radius being some universal constant for the simplicity of performance comparison. 

\paragraph{Generalization error in the one-pass setting.} A consequence of
Theorem~\ref{thm-gen-error-sghmc} is that the generalization error of the SGHMC1 iterates $| \E F_{\bZ}(X_k) - \E F(X_k)|$ in the one-pass setting satisfy
\beq\label{ineq-gen-bound-sghmc}
\tilde{\bigO}\left(\frac{(d+\beta)^{3/2}}{\mu_{\ast} \beta^{5/4}} \varepsilon +
\frac{(d+\beta)^{3/2}}{\beta(\mu_{\ast})^{3/2}}
\left((\log(1/\varepsilon))^{3/2}\delta^{1/4}
+\varepsilon\right)
\sqrt{\log(\mu_{\ast}^{-1}\log(\varepsilon^{-1}))}
+\frac{d+\beta}{\beta}\frac{\varepsilon^{2}}{\mu_{\ast}(\log(1/\varepsilon))^{2}}\right)\,,
\eeq
for $k = K_{SGHMC1}:=\tildeOm\left(\frac{1}{ \mu_* \varepsilon^4}\log^7(1/\varepsilon) \right)$ iterations, 
and similarly, Theorem~\ref{thm-gen-error-sghmc-2} implies the generalization error of the SGHMC2 iterates $| \E F_{\bZ}(\hat{X}_k) - \E F(\hat{X}_k)|$ in the one-pass setting satisfy
\beq\label{ineq-gen-bound-sghmc-2}
\tilde{\bigO}\left(\frac{(d+\beta)^{3/2}}{\mu_{\ast} \beta^{5/4}} \varepsilon +
\frac{(d+\beta)^{3/2}}{\beta\sqrt{\mu_{\ast}}}
\left(\sqrt{\log(1/\varepsilon)}\delta^{1/4}
+\varepsilon\right)
\sqrt{\log(\mu_{\ast}^{-1}\log(\varepsilon^{-1}))}\right)\,,
\eeq
for $k = K_{SGHMC2}:=\tildeOm\left(\frac{1}{ \mu_* \varepsilon^2}\log^2(1/\varepsilon) \right)$ iterations
(see the discussion in Section~\ref{app-parameter-dependence} in the Appendix for details).
On the other hand, the results of Theorem 1 in \cite{Raginsky} imply that the generalization error for the SGLD algorithm after $K_{SGLD}$ iterations in the one-pass setting scales as
\beq\label{ineq-gen-bound-sgld}
\tildeO\left(\frac{\beta(\beta+d)^{2}}{\lambda_*} \left( \delta^{1/4} \log(1/\varepsilon) + \varepsilon \right)\right) \quad \mbox{for} \quad K_{SGLD}=\tildeOm\left(\frac{\beta(d+\beta)}{\lambda_* \varepsilon^4}\log^5(1/\varepsilon) \right).
\eeq
The constants $\lambda_*$ (see \eqref{eq:lambda_star}) and $\mu_*$ (see Table~\ref{table_constants}) are exponentially small in both $\beta$ and $d$ in the worst case, but under some extra assumptions the dependency on $d$ can be polynomial (see e.g. \cite{Cheng}) although the exponential dependence to $\beta$ is unavoidable in the presence of multiple minima in general (see \cite{bovier2005metastability}). One can readily see that $K_{SGHMC2}$ has better dependency on $\epsilon$ than $K_{SGHMC1}$, and 
infer from \eqref{ineq-gen-bound-sghmc}--\eqref{ineq-gen-bound-sghmc-2}  that the performance of SGHMC2 is better than SGHMC1. Hence, in the rest of the section,
we will only focus on the comparison between SGHMC2 and SGLD.

{We see that the generalization error for 
SGHMC2 \eqref{ineq-gen-bound-sghmc-2} is bounded by 
\begin{equation}\label{ineq-sghmc-general-nodelta}
\tildeO\left(\frac{(d+\beta)^{3/2}}{\beta\mu_*}\left(\sqrt{\log(1/\varepsilon)}\delta^{1/4}+\varepsilon\right) \sqrt{\log\log(1/\varepsilon)}  \right)\,,
\end{equation}
as $\mu_*$ is small, and if we ignore the $\sqrt{\log\log(1/\varepsilon)}$ factor
\footnote{We emphasize that the effect of the last term $\sqrt{\log\log(1/\varepsilon)}$ appearing in \eqref{ineq-sghmc-general-nodelta} is typically negligible compared to other parameters. For instance even if $\varepsilon=2^{-2^{16}}$ is double-exponentially small, we have $\sqrt{\log\log(1/\varepsilon)}\leq 4$.},
then, we get
\begin{equation}\label{ineq-sghmc-general-nodelta-K}
\tildeO\left(\frac{(d+\beta)^{3/2}}{\beta\mu_*}\left(\sqrt{\log(1/\varepsilon)}\delta^{1/4}+\varepsilon\right)\right)
\quad
\text{for}\quad
K_{SGHMC2} = \tildeOm\left(\frac{1}{ \mu_* \varepsilon^2}\log^2(1/\varepsilon) \right),
\end{equation}
iterations of the SGHMC2 algorithm whereas the corresponding bound for SGLD from \cite[Theorem 1]{Raginsky} is
\beq\label{ineq-sgld-general-nodelta}\tildeO\left(\frac{\beta(\beta+d)^2}{\lambda_*}   
\left(\log(1/\varepsilon)\delta^{1/4}+\varepsilon\right) \right)\quad \mbox{for} \quad K_{SGLD}=\tildeOm\left(\frac{\beta(d+\beta)}{\lambda_* \varepsilon^4}\log^5(1/\varepsilon) \right)
\eeq
iterations of the SGLD algorithm. Note that $K_{SGHMC2}$ and $K_{SGLD}$ do not have the same dependency to $\varepsilon$ up to $\log$ factors (the former scales with
$\varepsilon$ as $\log^2(1/\varepsilon)\varepsilon^{-2}$ and the latter $\log^5(1/\varepsilon)\varepsilon^{-4}$), 
and this improvement on $\varepsilon$ dependency is due to better diffusion 
approximation of SGHMC2 (see Lemma~\ref{diff-app-2}) compared to SGLD and the exponential integrability estimate
we have in Lemma~\ref{lem:exp-integrable} which improves the estimate in \cite{Raginsky}
and using the same argument, one can improve the $\log^5(1/\varepsilon)/\varepsilon^{4}$ term in \eqref{ineq-sgld-general-nodelta} to $\log^3(1/\varepsilon)/\varepsilon^{4}$.}

{To make the comparison to SGLD simpler, {we notice that 
in both expressions \eqref{ineq-sghmc-general-nodelta-K}
and \eqref{ineq-sgld-general-nodelta}, we see a term scaling with $\delta^{1/4}$ due to the gradient noise level $\delta$} ($\delta$ is fixed in the one-pass setting),
and we fix the error in \eqref{ineq-sghmc-general-nodelta-K}
and \eqref{ineq-sgld-general-nodelta} {without the $\delta$ term} to be the same order, and then compare
the number of iterations $K_{SGHMC2}$ and $K_{SGLD}$. 
More precisely, given $\hat{\varepsilon}>0$ and we choose 
$\varepsilon>0$ such that $\frac{(d+\beta)^{3/2}}{\beta\mu_*}\varepsilon=\hat{\varepsilon}$
in \eqref{ineq-sghmc-general-nodelta-K} so that the generalization error
for SGHMC2 is 
\begin{equation} \label{eq:k-sghmc2}
\tilde{\mathcal{O}}\left(\hat{\varepsilon}+\frac{(d+\beta)^{3/2}}{\beta\mu_{\ast}}
\sqrt{\log\left(\frac{(d+\beta)^{3/2}}{\beta\mu_{\ast}\hat{\varepsilon}}\right)}\delta^{1/4}\right)
\quad
\text{for}
\quad
K_{SGHMC2}=\tilde{\Omega}\left(\frac{(d+\beta)^{3}}{\beta^{2}\mu_{\ast}^{3}\hat{\varepsilon}^{2}}
\log^{2}\left(\frac{(d+\beta)^{3/2}}{\beta\mu_{\ast}\hat{\varepsilon}}\right)\right).
\end{equation}
Similarly, the generalization error for SGLD is
\begin{equation}\label{eq:k-sgld}
\tilde{\mathcal{O}}\left(\hat{\varepsilon}+\frac{\beta(\beta+d)^{2}}{\lambda_{\ast}}
\log\left(\frac{\beta(\beta+d)^{2}}{\lambda_{\ast}\hat{\varepsilon}}\right)\delta^{1/4}\right)
\quad
\text{for}
\quad
K_{SGLD}=\tilde{\Omega}\left(\frac{\beta^{5}(d+\beta)^{9}}{\lambda_{\ast}^{5}\hat{\varepsilon}^{4}}
\log^{5}\left(\frac{\beta(\beta+d)^{2}}{\lambda_{\ast}\hat{\varepsilon}}\right)\right).
\end{equation}
When $\lambda_*$ and $\mu_*$ are on the same order or $\mu_*$ is larger, since typically $\beta \geq 1$, the term involving $\delta$ in the generalization error
for SGHMC2 above is (smaller) better than the counterpart for SGLD, and this is guaranteed to be achieved in a less number of iterations ignoring the log factors and universal constants for $K_{SGHMC2}$ in \eqref{eq:k-sghmc2} and $K_{SGLD}$ in \eqref{eq:k-sgld}. 

Comparing $\lambda_*$ and $\mu_*$ on arbitrary non-convex functions seems not trivial, however we give a class of non-convex functions (see Proposition~\ref{prop:comparison} and Example~\ref{example-sym-well}) where $\frac{1}{\mu_*} = \tildeO \left(\sqrt{\frac{1}{\lambda_*}}\right)$. For this class, we can infer from \eqref{eq:k-sghmc2} that $K_{SGHMC2}$ has a dependency $1/\mu_{\ast}^{3}= \tildeO (1/\lambda_{\ast}^{3/2})$ which is much smaller in contrast to $1/\lambda_{\ast}^{5}$ for $K_{SGLD}$ in \eqref{eq:k-sgld}. 


\paragraph{Empirical risk minimization.} The empirical risk minimization bound given in Corollary~\ref{coro:main-exp-2} has an additional term $\mathcal{J}_2$ compared to the $\overline{\mathcal{J}}_0(\varepsilon)$ and $\hat{\mathcal{J}}_1(\varepsilon)$ terms appearing in the one-pass generalization bounds. Note also that $\mathcal{J}_0(\bz,\varepsilon) \le \overline{\mathcal{J}}_0(\varepsilon)$. 
As a consequence, SGHMC2 algorithm has expected empirical risk
\begin{align}\label{ineq-exp-emp-risk-compare-sghmc}
&\tilde{\bigO}\left(\frac{(d+\beta)^{3/2}}{\mu_{\ast} \beta^{5/4}} \varepsilon +
\frac{(d+\beta)^{3/2}}{\beta\sqrt{\mu_{\ast}}}
\left(\sqrt{\log(1/\varepsilon)}\delta^{1/4}
+\varepsilon\right)
\sqrt{\log(\mu_{\ast}^{-1}\log(\varepsilon^{-1}))}
+ d\cdot \frac{\log(1+\beta)}{\beta}\right)\ \,,
\end{align}
after $K_{SGHMC2} = \tildeOm\left(\frac{1}{ \mu_* \varepsilon^2}\log^2(1/\varepsilon) \right)$
iterations as opposed to
\beq\label{ineq-exp-emp-risk-compare-sgld} \tildeO\left(\frac{\beta(\beta+d)^2}{\lambda_*} \left( \delta^{1/4} \log(1/\varepsilon) + \varepsilon \right) + d\cdot \frac{\log(1+\beta)}{\beta}\right)\,, \eeq
after
$K_{SGLD}=\tildeOm\left(\frac{\beta(d+\beta)}{\lambda_* \varepsilon^4}\log^5(1/\varepsilon) \right)
$ iterations required in \cite{Raginsky}. Comparing \eqref{ineq-exp-emp-risk-compare-sghmc} and \eqref{ineq-exp-emp-risk-compare-sgld}, we see that the last terms are the same. If this term is the dominant term, then the empirical risk upper bounds for SGLD and SGHMC2 will be similar except that $K_{SGHMC2}$ can be smaller than $K_{SGLD}$ for instance when $\frac{1}{\mu_*} = \tilde{\bigO} \left(\sqrt{\frac{1}{\lambda_*}}\right)$. Otherwise, if the last term is not the dominant one and can be ignored with respect to other terms, then, the performance comparison will be similar to the discussion about the generalization bounds \eqref{ineq-sghmc-general-nodelta} and \eqref{ineq-sgld-general-nodelta} discussed above.

We next briefly discuss the comparisons of SGHMC2 and SGLD based on the total number of stochastic gradient evaluations (gradient complexity), and we compare with a recent work \cite{xu2018global} which established a faster convergence result and improved the gradient complexity for SGLD in the mini-batch setting compared with \cite{Raginsky}. Here, the total number of stochastic gradient evaluations of an algorithm is defined as
the number of stochastic gradients calculated per iteration (which is equal to the batch size in the mini-batch setting) times the total number of iterations. 
 \cite{xu2018global} showed that it suffices to take
\begin{equation}
\hat{K}_{SGLD}=\tilde{\Omega}\left(\frac{d^7}{\hat{\lambda}^5 \heps^5} \right)\label{ineq-iter-sgld-erm}
\end{equation} 
stochastic gradient evaluations to converge to an $\heps$ neighborhood of an almost ERM where $\tilde{\Omega}(\cdot)$ hides some factors in $\beta$ and $\hat{\lambda}$ is the spectral gap of the discrete overdamped Langevin dynamics, i.e. SGLD with zero gradient noise. This improves upon the result in \cite{Raginsky} which showed that the same task requires $
\tilde{\Omega}\left(\frac{d^{17}}{\lambda_*^9 \heps^8} \right)$ stochastic gradient evaluations.
Our results show that (see e.g. \eqref{ineq-exp-emp-risk-compare-sghmc}) for SGHMC2, it suffices to have
\begin{equation} 
\hat{K}_{SGHMC2} = \tilde{\Omega}\left(\frac{d^9}{ \mu_*^4 \hat \varepsilon^6 }\right)
\label{ineq-iter-sghmc-erm}
\end{equation}
stochastic gradient evaluations, ignoring the $\log$ factors in the parameters $\hat \varepsilon, \mu_*, d$ and hiding factors in $\beta$ that can be made explicit. To see \eqref{ineq-iter-sghmc-erm}, we infer from \eqref{ineq-exp-emp-risk-compare-sghmc} that for fixed precision 
$\hat{\varepsilon}>0$ and dimension $d$, by ignoring the log factors and $\beta$,
we can choose $\varepsilon$ so that $d^{3/2} \varepsilon/\mu_* =\hat{\varepsilon}$ and choose the gradient noise level $\delta$ so that $d^{3/2}\delta^{1/4}/\sqrt{\mu_*} =\hat{\varepsilon}$. So the number of SGHMC2 iterations is
\begin{equation*} 
K_{SGHMC2} = \tildeOm\left(\frac{1}{ \mu_* \varepsilon^2}\right) = \tildeOm\left(\frac{d^3}{ \mu_*^2 \hat{\varepsilon}^2}\right).
\end{equation*}
On the other hand, the mini-batch size to achieve gradient noise level $\delta$ is given by $1/\delta$ (see \cite{Raginsky}), which is equal to $d^6/(\mu_*^2 \hat{\varepsilon}^4)$. Hence, we obtain \eqref{ineq-iter-sghmc-erm} which is the product of the mini-batch size and number of iterations.

It is hard to compare $\hat{\lambda}$ in \eqref{ineq-iter-sgld-erm} and $\mu_*$ in \eqref{ineq-iter-sghmc-erm} in general since $\hat{\lambda}$ is the spectral gap
of the discrete overdamped Langevin dynamics (i.e. SGLD with zero gradient noise) without a simple closed-form formula.
However, when the stepsize is small enough, we expect $\hat{\lambda}$ will be similar to $\lambda_*$, which is the spectral gap of the continuous-time overdamped Langevin diffusion. 
As a consequence, when the stepsize $\eta$ is small enough (which is the case for instance, when target accuracy $\heps$ is small enough), we will have $\hat{\lambda} \approx \lambda_*$ and 
$\frac{1}{\mu_*} = \mathcal{O}\left(\sqrt{\frac{1}{\lambda_*}}\right) = \mathcal{O}\left(\sqrt{\frac{1}{\hat\lambda}}\right) $ 
for the class of non-convex functions we discuss in Proposition~\ref{prop:comparison} and Example~\ref{example-sym-well}. 
For this class of problems, comparing \eqref{ineq-iter-sgld-erm} and \eqref{ineq-iter-sghmc-erm}, we see that we obtain an improvement in the spectral gap parameter ($\mu_*^4$ vs. $\hat {\lambda}^5$), however $\heps$ and $d$ dependency of the bound \eqref{ineq-iter-sgld-erm} is better than \eqref{ineq-iter-sghmc-erm}.

\paragraph{Population risk minimization.} If samples are recycled and multiple passes over the dataset is made, then one can see from Corollary~\ref{coro:main} that there is an extra term $\mathcal{J}_3$ that needs to be added to the bounds given in \eqref{ineq-exp-emp-risk-compare-sghmc} and \eqref{ineq-exp-emp-risk-compare-sgld}. This term satisfies
$$\mathcal{J}_3 = \tildeO\left(\frac{(\beta+d)^2}{\lambda_* n}\right).$$
If this term is dominant compared to other terms $\overline{\mathcal{J}}_0, \mathcal{J}_1$ and $\mathcal{J}_2$, for instance this may happen if the number of samples $n$ is not large enough, then the performance guarantees for population risk minimization via SGLD and SGHMC2 will be similar. Otherwise, if $n$ is large and $\beta$ is chosen in a way to keep the $\mathcal{J}_2$ term on the order $\overline{\mathcal{J}_0}$, 
then similar improvement can be achieved.
\paragraph{Comparison of $\lambda_*$ and $\mu_*$.} 
The parameters $\lambda_*$ (see \eqref{eq:lambda_star}) and $\mu_*$ (see Table~\ref{table_constants}) govern the convergence rate to the equilibrium of the overdamped and underdamped Langevin SDE, they can be both exponentially small in dimension $d$ and in $\beta$. They appear naturally in the complexity estimates of SGHMC2 and SGLD method as these algorithms can be viewed as discretizations of Langevin SDEs (when the discretization step is small and the gradient noise $\delta=0$, the discrete dynamics will behave similarly as the continuous dynamics). 
Next, to get further intuition, first we discuss some toy examples of non-convex functions below 
where $\frac{1}{\mu_*} = \mathcal{O}\left(\sqrt{\frac{1}{\lambda_*}}\right)$. For these examples if the other parameters $(\beta,d,\delta)$ are fixed, then SGHMC2 can lead to an improvement upon the SGLD performance. We will then show in Proposition~\ref{prop:comparison} that these examples generalize to a more general class of non-convex functions.

\begin{example}\label{example-sym-well} Consider the following symmetric double-well potential in $\R^d$ studied previously in the context of Langevin diffusions (\cite{Eberle}): 
	 \begin{equation*} f_a(x) = U(x/a) \quad \mbox{with} \quad U(x) := \begin{cases}
	\frac{1}{2}(\|x\| - 1)^2 & \mbox{for} \quad \|x\| \geq \frac{1}{2}, \\
    \frac{1}{4} - \frac{\|x\|^2}{2}  & \mbox{for} \quad \|x\| \leq \frac{1}{2},
\end{cases}
\end{equation*}
where $a>0$ is a scaling parameter which is illustrated in the left panel of Figure~\ref{fig-double-well}. For this example, there are two minima that are apart at a distance $\mathcal{R} = \bigO(a)$. For simplicity, we assume there is only one sample, i.e. $\bz = (z_1)$ and $F_\bz(x) = f(x,z_1) = f_a(x)$. 
We consider the non-convex optimization problem \eqref{pbm-finite-sum} with both the SGHMC2 algorithm and the SGLD algorithm. 
\cite{Eberle} showed that $\mu_* \geq \Theta(\frac{1}{a})$ for this example whereas $\lambda_* \leq \Theta(\frac{1}{a^2})$ making the constants hidden by the $\Theta$ explicit. 
This shows that the contraction rate of the underdamped diffusion $\mu_*$ is (faster) larger than that of the overdamped diffusion $\lambda_*$ by a square root factor when $a$ is large where all the constants can be made explicit. 
Such results extend to a more general class of non-convex functions with multiple-wells and higher dimensions as long as the gradient of the objective satisfies a growth condition (see Example 1.1, Example 1.13 in \cite{Eberle} for a further discussion). 

For computing an $\varepsilon$-approximate global minimizer of $f_a=f(x,z_1)$ (or more generally for a non-convex problem satisfying Assumption~\ref{assumptions}), $\beta$ is chosen large enough so that the stationary measure concentrates around the global minimizer. 
Using the tight characterization of $\lambda_*$ from Theorem 1.2 in \cite{bovier2005metastability} for $\beta$ large, further comparisons with similar conclusions between the rate of convergence to the equilibrium distribution between the underdamped and overdamped dynamics can also be made. For example, consider the non-convex objective $F_\bz(x) = \tilde{f}_a(x)=\tilde{U}(x/a)$ instead, illustrated in the right panel of Figure~\ref{fig-double-well} for $a=4$ where  
\beq \tilde{U}(x) = \begin{cases}
	\frac{1}{2}(x - 1)^2 & \mbox{for} \quad x \geq \frac{1}{2}, \\
    \frac{1}{4} - \frac{x^2}{2}  & \mbox{for} \quad -\frac{1}{8} \leq x \leq \frac{1}{2}, \\
    \frac{1}{2}(x+\frac{1}{4})^2 + \frac{15}{64} & \mbox{for} \quad x\leq -\frac{1}{8}, 
\end{cases} \nonumber
\eeq
is the asymmetric double well potential in dimension one. It follows from Theorem~\ref{exp-converg} (see also \cite{Eberle}) that the contraction rate satisfies 
$\mu_* = \Theta\left(a^{-1}\right)\,,$ whereas it follows from Theorem 1.2 in \cite{bovier2005metastability} that $\lambda_* = \Theta(1/a^2)$. This shows that when the separation between minima, or alternatively the scaling factor $a$ is large enough, $\mu_*$ is larger than $\lambda_*$ by a square root factor up to constants.
\begin{figure}
\begin{center}
\includegraphics[scale=0.55]{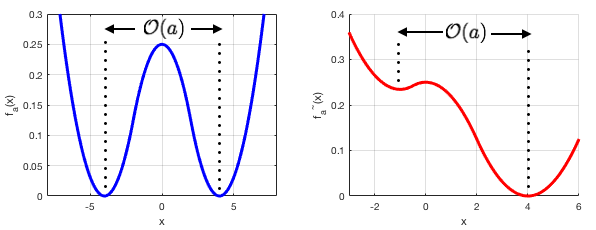}
\caption{\label{fig-double-well}The illustration of the functions $f_a(x)$ (left) and $\tilde{f}_a(x)$ (right) for $a=4$.}
\end{center}
\end{figure}
\end{example} 

The behavior in these toy examples can be generalized to more general non-convex objectives with a finite-sum structure satisfying Assumption~\ref{assumptions}. Proposition~\ref{prop:comparison} below gives a class of functions where $\mu_*$ is on the order of the square root of $\lambda_*$. The proof will be presented in details in Section~\ref{sec:comparison}.  

\begin{proposition}\label{prop:comparison}
Suppose that the functions $f_a(x,z)$ indexed by $a$ satisfies Assumption~\ref{assumptions} (i)-(iii) 
with $m=m_{1}a^{-2}$, $M=M_{1}a^{-2}$ and $B=B_{1}a^{-1}$ for some fixed constants $m_{1}$, $M_{1}$, and $B_{1}$. 
Then, we have as $a \rightarrow \infty,$
\begin{equation}
\lambda_{\ast}=\mathcal{O}(a^{-2}),
\qquad
\mu_{\ast}=\Theta(a^{-1}).
\end{equation}
\end{proposition}

This result is more general than the previous
example. In particular, if $f(x,z)$ satisfies Assumption~\ref{assumptions} (i)-(iii) with $m,M,B$ replaced
by $m_{1},M_{1},B_{1}$, then $f_a(x,z):=f(x/a, z)$ satisfies Assumption~\ref{assumptions} (i)-(iii) with $m=m_{1}a^{-2}$, $M=M_{1}a^{-2}$ and $B=B_{1}a^{-1}$.
 Proposition~\ref{prop:comparison} essentially
says that if we consider the normalized empirical risk objective $F_\bz(x/a) = \frac{1}{n} \sum_{i=1}^n f(x/a, z_i)$ where $a$ is a (normalization) scaling parameter and $f(x, z)$ satisfies Assumption~\ref{assumptions}, then for large enough values of $a$, the empirical risk surface will be relatively flat and the convergence rate of momentum variant SGHMC2 to an $\varepsilon$-neighborhood of the global minimum will be governed by the parameter $\mu_*$ which will be larger than that of the parameter $\lambda_*$ of SGLD when $a$ is sufficiently large. This will lead to improved performance bounds for SGHMC2 compared to known performance bounds for SGLD. 

\section{Applications}

We note that several non-convex stochastic optimization problems of interest can satisfy Assumption~\ref{assumptions} under appropriate noise assumptions for the underlying dataset. For example, Lasso problems with non-convex regularizers (see e.g. \cite{hu2017group}), non-convex formulations of the phase retrieval problem around global minimum (see e.g. \cite{zhang2017nonconvex}) or non-convex stochastic optimization problems defined on a compact set including but not limited to dictionary learning over the sphere (see e.g. \cite{sun2016complete}), training deep learning models subject to norm constraints in the model parameters (see e.g. \cite{anil2019sorting}).
In this section, we discuss some applications of our results where we provide two specific examples.

\subsection{Binary linear classification}\label{example-binary}
In linear binary classification, the aim is to learn a predictive model of the form $\mathbb{P}(Y = 1 | A_{in} =a) = \sigma_c(\langle \tilde{x},a \rangle )$, where $\tilde{x}\in \mathbb{R}^d$ is a parameter vector to be learned, $A_{in}$ is the input variable (feature vector), $Y$ is the binary output and $\sigma_c:\mathbb{R}\to[0,1]$ is a threshold function. Binary classification 
arises in many data-driven applications in operations research from diagnosing patients in healthcare \cite{wu2007robust} to predicting directions in the stock market \cite{james2013introduction}.

A number of empirical studies have demonstrated that non-convex choices of the $\sigma_c$ function can often lead to superior classification accuracy and robustness properties compared to convex choices of $\sigma_c$ such as the hinge loss \cite{chapelle2009tighter,collobert2006trading,wu2007robust,nguyen2013algorithms}. Given access to a dataset of input-output pairs $z_i = (a_i,y_i)$, a standard way of estimating $\tilde{x}$ is based on minimizing the \emph{regularized squared loss} over the dataset, i.e. 
\beq \min_{x\in\mathbb{R}^d}  \frac{1}{n} \sum_{i=1}^n \left(y_i - \sigma_c(\langle x,a_i \rangle )\right)^2 + \frac{\lambda_{r}}{2}\|x\|^2,
\label{eq-erm-binary}
\eeq
where $\lambda_{r}>0$ is a regularization parameter that may depend on the number of samples $n$. By Lagrangian duality, this problem is equivalent to the constrained optimization problem
$$ \min_{x\in\mathbb{R}^d}  \frac{1}{n} \sum_{i=1}^n \left(y_i - \sigma_c(\langle x,a_i \rangle )\right)^2 \quad \mbox{subject to}  \quad \|x\| \leq R,
$$
for some $R$, which has also been considered in the literature (see e.g. \cite{mei2018landscape,foster2018uniform,wang2019differentially}).
For non-convex $\sigma_c(\cdot)$, this problem is also non-convex in general. We consider minimizing the objective \eqref{eq-erm-binary} in the mini-batch setting where the gradients in SGHMC iterations are estimated from $n_b$ data points sampled with replacement, i.e. the gradient is estimated as \looseness=-1
\begin{equation} g(x,U_{{\mathbf{z}}}) = \frac{1}{n_b}\sum_{j=1}^{n_b} \nabla f(x,z_{j})\,, 
\label{eq-grad-stoc}
\end{equation}
where $z_{j}$ are i.i.d. with a uniform distribution over the set of indices $\{1,2,\dots,n\}$. We also consider the following assumption for the threshold function $\sigma_c$ which are satisfied by many choices of $\sigma_c$ in practice. A prominent example is the logistic (or sigmoid) function in which case $\sigma_c(z) = 1/(1+e^{-z})$ which is also used in deep learning. Another possible choice is the \emph{probit function} which corresponds to $\sigma_c(t) = \Phi(t)$ where $\Phi$ is the cumulative distribution function of the standard normal distribution. 
\begin{assumption}\label{assump-binary} The threshold function $\sigma_c$ is twice continuously differentiable on $\mathbb{R}$. It is bounded and has bounded first and second derivatives, i.e. there exists a constant $L_{\sigma_c}>0$ such that $\max \big\{ \|\sigma_c\|_\infty, \|\sigma_c'\|_\infty, \|\sigma_c''\|_\infty \big\} \leq L_{\sigma_c}.$ The distribution of the input data $A_{in}$ has compact support, i.e. $\|A_{in}\|\leq D$ for some $D>0$.
\end{assumption}
We show in the next lemma that if Assumption~\ref{assump-binary} holds, then Assumption~\ref{assumptions} holds with explicit constants $A_0, B, M, m, b$ and $\sigma_c$ that we can precise. The proof can be found in the Appendix.

\begin{lemma}\label{lem-example-classification} In the setting of binary linear classification, consider the SGHMC method applied to the objective \eqref{eq-erm-binary} where gradients are estimated according to \eqref{eq-grad-stoc} where the probability law $\mu_0$ of the initial state has compact support. If Assumption~\ref{assump-binary} holds; then Assumption~\ref{assumptions} hold for any $\delta \in [\frac{1}{4n_b},1)$ with the following constants:
\beq 
A_0 &=&  (1 + \|\sigma_c(0)\|)^2, \quad B =  2D \left(1 + \|\sigma_c\|_\infty \right) \|\sigma_c'\|_\infty, \\
 M &=& 2D^2 \|\sigma_c'\|_\infty^2  + 2D^2(1+\|\sigma_c\|_\infty) \|\sigma_c''\|_\infty  + 5\lambda_r, \\
 m &=& \lambda/2, \quad b = 8 (1+\|\sigma_c\|_\infty)^2 \|\sigma_c'\|_\infty^2 D^2 / \lambda_r.
\eeq
\end{lemma}

We conclude from Lemma \ref{lem-example-classification} that the objective is dissipative and our main results for SGHMC1 and SGHMC2 algorithms described in Sections \ref{sec:mainresult}--\ref{sec-compare} apply to binary linear classification under Assumption \ref{assump-binary} with the constants given in Lemma \ref{lem-example-classification} and where $\mu_*$ is given by the formula in Table \ref{table_constants}. For example, if $D = \mathcal{O}(1)$, then we have $\frac{1}{\mu_*} = \tilde{\Theta}(\sqrt{d+\beta} e^{\tilde{\Theta}(d+\beta)})$ (see \eqref{eq:mu-ast}) and we conclude from \eqref{ineq-iter-sghmc-erm} that it suffices to have\looseness=-1 \begin{equation*} 
\hat{K}_{SGHMC2} = \tilde{\Omega}\left(\frac{d^9}{ \mu_*^4 \hat \varepsilon^6 }\right) = \tilde{\Omega}\left(\frac{d^9 e^{\tilde{\Theta}(d+\beta)}}{ (d^2 + \beta^2) \hat \varepsilon^6 }\right)
\end{equation*}
stochastic gradient evaluations to converge to an $\heps$ neighborhood of an almost ERM ignoring the $\log$ factors in the parameters $\hat \varepsilon, \mu_*, d$ and hiding other constants that can be made explicit based on Lemma \ref{lem-example-classification}.\footnote{We also note that under further assumptions on the statistical nature of the input and if the number of data points is large enough, it can be shown that the objective \eqref{eq-erm-binary} admits a unique minimizer and the objective is strongly convex in some regions \cite{mei2018landscape}. However, our assumptions here are weaker, therefore such arguments are not directly applicable.}

\subsection{Robust Ridge Regression}\label{ex-robust-reg} 
Given an input (feature) vector $A_{in}\in\mathbb{R}^d$, the aim is to predict the output $Y\in \mathbb{R}$. Given access to a dataset of input-output pairs $z_i = (a_i,y_i)$, we assume a linear model 
$y_i = a_i^T \tilde{x} + \varepsilon_i$ where the errors $\varepsilon_i$ are i.i.d. with mean zero. The standard \emph{ridge regression} estimate of $\tilde{x}$ minimizes a penalized residual sum of squares \cite{hoerl1970ridge}, i.e. minimizes $\sum_{i=1}^n \|y_i - \langle x,a_i \rangle \|^2 +\lambda_r \|x\|^2$ where $\lambda_r>0$ is a regularization parameter.\footnote{See \cite{kilmer2001choosing} for details regarding the choice of the parameter $\lambda_r$.} However, this formulation can be sensitive to outliers. Robust formulations of the ridge regression \cite{razavi2012robust} can be obtained if one solves instead the following problem 
  \beq \min_{x\in\mathbb{R}^d}\frac{1}{n} \sum_{i=1}^n f(x,z_i), \quad f(x,z_i) = \rho\left(y_i - \langle x,a_i \rangle \right) + \frac{\lambda_r}{2} \|x\|^2,
 \label{eq-robust-reg-erm}
 \eeq
where $\lambda_r>0$ is a regularization parameter 
and $\rho:\mathbb{R}\to\mathbb{R}$ is a suitably chosen loss function. In particular, for achieving robustness to outliers, the non-convex choices of the function 
$\rho$ that are either bounded or slowly growing near infinity has been considered in the literature (as opposed to the standard ridge regression setting which corresponds to $\rho(t) = \|t\|^2$). For example, popular choices of the function $t\mapsto \rho(t)$ include \emph{Tukey's bisquare loss} defined as 
$$ \rho_{Tukey}(t) =  \begin{cases}
    1 -( 1 - (t/t_0)^2)^3  & \mbox{for} \quad \|t\|\leq t_0, \\ 
    1 &  \mbox{for} \quad | t| \geq t_0, 
\end{cases}
$$ 
(see e.g. \cite{mei2018landscape}) and exponential squared loss \cite{exp-loss}: $\rho_{exp}(t) = 1 - e^{-\|t\|^2/t_0}$, where $t_0>0$ is a tuning parameter. 
In the following, similar to \cite{wang2019differentially}, we assume that the data $A_{in}$ is bounded and the threshold function and its derivatives up to order two are bounded, similar to  \cite{mei2018landscape}. This assumption for $\rho$ is satisfied in several cases, including Tukey's bisquare loss and exponential squares loss mentioned above.
\begin{assumption}\label{assump-robust} The function $\rho$ is twice continuously differentiable on $\mathbb{R}$. The function $\rho$ is bounded and has bounded first and second derivatives; i.e. there exists a constant $L_\rho$ such that $\max(\|\rho\|_\infty,\|\rho'\|_\infty, \|\rho''\|_\infty) \leq L_\rho$. Furthermore, the distribution of the input data $A_{in}$ has compact support, i.e. there exists $D$ such that $\|A_{in}\|\leq D$.
\end{assumption}

The following lemma shows that under Assumption \ref{assump-robust}, our assumptions (Assumption \ref{assumptions}) for analyzing SGHMC methods hold with proper initialization.

\begin{lemma}\label{lem:2} 
In the setting of robust regression, consider the objective \eqref{eq-erm-binary} where gradients are estimated according to \eqref{eq-grad-stoc} where the probability law $\mu_0$ of the initial state has compact support. If Assumption~\ref{assump-binary} holds; then Assumption~\ref{assumptions} hold for both SGHMC1 and SGHMC2 methods for any choice of $\delta \in [\frac{1}{4n_b},1)$ with the following constants:
 \beq 
 A_0 &=& \|\rho\|_\infty, \qquad B =  4\|\rho'\|_\infty D,\\
 M &=& \|\rho''\|_\infty D^2 + \lambda_r, 
 \qquad m=\lambda_r/2, \qquad b = \frac{2 \|\rho'\|_\infty^2 D^2}{\lambda_r}.
 \eeq
\end{lemma}  

Similarly, we conclude from Lemma \ref{lem:2} that our main results for SGHMC1 and SGHMC2 algorithms described in Sections \ref{sec:mainresult}--\ref{sec-compare} apply to the problem of robust regression under Assumption~\ref{assump-robust}.

\section{Outline of the Proof}\label{sec:proof:overview}

To obtain the main results in this paper, we adapt the proof techniques of \cite{Raginsky} developed for the overdamped dynamics to the underdamped dynamics and combine it with the analysis of \cite{Eberle} which quantifies the convergence rate of the underdamped Langevin SDE to its equilibrium. In an analogy to the fact that momentum-based first-order optimization methods require a different Lyapunov function and a quite different set of analysis tools (compared to their non-accelerated variants) to achieve fast rates (see e.g. \cite{lu2018accelerating,su2014differential,nesterov1983method}), our analysis of the momentum-based SGHMC1 and SGHMC2 algorithms requires studying a different Lyapunov function $\mathcal{V}$ defined in \eqref{eq:lyapunov} that also depends on the objective $f$ as opposed to the classic Lyapunov function $\mathcal{H}(x) = \| x \|^2$ arising in the study of the SGLD algorithm (see e.g. \cite{mattingly2002ergodicity,Raginsky}). This fact introduces some challenges for the adaptation of the existing analysis techniques for SGLD to SGHMC. For this purpose, we take the following steps:

First, we show that SGHMC1 and SGHMC2 iterates track the underdamped Langevin diffusion closely in the 2-Wasserstein metric. As this metric requires finiteness of second moments, we first establish uniform (in time) $L^2$ bounds for both the underdamped Langevin SDE and SGHMC1 and SGHMC2 iterates (see Lemma~\ref{lem:L2bound} and Lemma ~\ref{lem:L2bound:Jordan} in Appendix), exploiting the structure of the Lyapunov function $\mathcal{V}$. Second, we obtain a bound for the Kullback-Leibler divergence between the discrete and continuous underdamped dynamics making use of the Girsanov theorem, which is then converted to bounds in the 2-Wasserstein metric by an application of an optimal transportation inequality of \cite{BV}. This step requires proving a certain exponential integrability property of the underdamped Langevin diffusion (Lemma~\ref{lem:exp-integrable} in Appendix). 
We show in Lemma~\ref{lem:exp-integrable} that the exponential moments grow at most linearly
in time, which strictly improves the exponential growth in time in Lemma~4 in \cite{Raginsky}.
\footnote{The method that is used in the proof of Lemma~\ref{lem:exp-integrable} in Appendix can indeed be adapted
to improve the exponential integrability and hence the overall estimates in \cite{Raginsky} for
SGLD as well.} As a result, the method improves upon the $\varepsilon$ dependence of
the number of iterates (see equations~\eqref{ineq-sghmc-general-nodelta-K} and \eqref{ineq-sgld-general-nodelta}).

Second, we apply the seminal result of \cite{Eberle} which showed that the continuous-time underdamped Langevin SDE is geometrically ergodic with an explicit rate $\mu_*$ in the 2-Wasserstein metric. In order to get explicit performance guarantees, we derive new bounds that make the dependence of the constants to the initialization in \cite{Eberle} explicit (see Lemma~\ref{lem:nu0-piz} in Appendix).

As the $x$-marginal of the equilibrium distribution $\pi_\bz(dx,dv)$ of the underdamped Langevin SDE concentrates around the global minimizers of $F_{\bz}$ for $\beta$ appropriately chosen, and we can control the error between the discrete-time SGHMC1 and SGHMC2 dynamics and the underdamped SDE by choosing the step size accordingly; this leads to performance bounds for the empirical risk minimizations for SGHMC1 and SGHMC2 algorithms in Corollary~\ref{coro:main-exp} and Corollary~\ref{coro:main-exp-2}.
For controlling the population risk during SGHMC iterations, in addition to the empirical risk, one has to control the \emph{generalization error} $F(X_k)-F_{\bZ}(X_k)$ that accounts for the differences between the finite sample size problem \eqref{pbm-finite-sum} and the original problem \eqref{opt-pbm}. By exploiting the fact that the $x-$marginal of the invariant distribution for the underdamped dynamics is the same as it is in the overdamped case, we control the generalization error 
in Corollary~\ref{coro:main} and Corollary~\ref{coro:main-2} which is no worse than that of the available bounds for SGLD given in \cite{Raginsky}.

\section{Conclusion}
SGHMC is a momentum-based popular variant of stochastic gradient where a controlled amount of isotropic Gaussian noise is
added to the gradient estimates for optimizing a non-convex function. We obtained first-time finite-time guarantees for the
convergence of SGHMC1 and SGHMC2 algorithms to the $\varepsilon$-global minimizers under
some regularity assumption on the non-convex objective $f$.
We also show that on a class of non-convex problems, SGHMC2 can be faster than overdamped Langevin MCMC approaches such as
SGLD in the sense that the best available bounds for
SGHMC2, which we prove in our paper, are better than the best available bounds for SGLD. This
effect is due to the momentum term in the underdamped SDE. Furthermore, our results show that momentum-based acceleration is possible on a class of non-convex problems under some conditions if we compare known upper bounds between SGLD and SGHMC.
Finally, we mention a few limitations in our work that may lead to
some future research directions.
In our paper, the performance dependence on dimension is exponential in general.
In the future, we will investigate for what class of (non-convex) target functions
$f$ we can obtain performance bound independent of dimension $d$ or has polynomial dependence on $d$. In addition, our results {suggest} 
{that momentum-based SGHMC methods will work particularly well} when the (non-convex) target functions have relatively flat landscapes.
In the future, we will investigate whether we can obtain {theoretical results for SGHMC} on a wider class of non-convex problems.

\section*{Acknowledgements}
We thank Agostino Capponi, Xiuli Chao, Wenbin Chen, Jim Dai, Murat A. Erdogdu, Fuqing Gao, 
Jianqiang Hu, Jin Ma, Sanjoy Mitter, Asuman Ozdaglar, Pablo Parrilo, Umut \c{S}im\c{s}ekli, and S. R. S. Varadhan for helpful discussions. 
Xuefeng Gao acknowledges support from Hong Kong RGC Grants 24207015 and 14201117.
Mert G\"{u}rb\"{u}zbalaban's research is supported in part by the grants NSF DMS-1723085 and NSF CCF-1814888.
Lingjiong Zhu is grateful to the support from the grant NSF DMS-1613164.


\bibliographystyle{alpha}
\bibliography{langevin}

\newcommand{\etalchar}[1]{$^{#1}$}
\begin{thebibliography}{{\c{S}im\c{s}ekli}YN{\etalchar{+}}18}

\bibitem[AKW12]{ahn2012bayesian}
Sungjin Ahn, Anoop Korattikara, and Max Welling.
\newblock Bayesian posterior sampling via stochastic gradient {F}isher scoring.
\newblock In {\em International Conference on Machine Learning}, pages
  1771--1778, 2012.

\bibitem[ALG19]{anil2019sorting}
Cem Anil, James Lucas, and Roger Grosse.
\newblock Sorting out {L}ipschitz function approximation.
\newblock In {\em International Conference on Machine Learning}, pages
  291--301, 2019.

\bibitem[AZH16]{allen2016variance}
Zeyuan Allen-Zhu and Elad Hazan.
\newblock Variance reduction for faster non-convex optimization.
\newblock In {\em International Conference on Machine Learning}, pages
  699--707, 2016.

\bibitem[BBG14]{betancourt2014optimizing}
MJ~Betancourt, Simon Byrne, and Mark Girolami.
\newblock Optimizing the integrator step size for {H}amiltonian {M}onte
  {C}arlo.
\newblock {\em arXiv preprint arXiv:1411.6669}, 2014.

\bibitem[BBLG17]{betancourt2017geometric}
Michael Betancourt, Simon Byrne, Sam Livingstone, and Mark Girolami.
\newblock The geometric foundations of {H}amiltonian {M}onte {C}arlo.
\newblock {\em Bernoulli}, 23(4A):2257--2298, 2017.

\bibitem[Bet17]{betancourt2017conceptual}
Michael Betancourt.
\newblock A conceptual introduction to {H}amiltonian {M}onte {C}arlo.
\newblock {\em arXiv preprint arXiv:1701.02434}, 2017.

\bibitem[BGK05]{bovier2005metastability}
Anton Bovier, V{\'e}ronique Gayrard, and Markus Klein.
\newblock Metastability in reversible diffusion processes {II}: Precise
  asymptotics for small eigenvalues.
\newblock {\em Journal of the European Mathematical Society}, 7(1):69--99,
  2005.

\bibitem[BLNR15]{Belloni}
Alexandre Belloni, Tengyuan Liang, Hariharan Narayanan, and Alexander Rakhlin.
\newblock Escaping the local minima via simulated annealing: Optimization of
  approximately convex functions.
\newblock In {\em Conference on Learning Theory}, pages 240--265, 2015.

\bibitem[BM99]{Borkar}
Vivek~S Borkar and Sanjoy~K Mitter.
\newblock A strong approximation theorem for stochastic recursive algorithms.
\newblock {\em Journal of Optimization Theory and Applications},
  100(3):499--513, 1999.

\bibitem[BT93]{bertsimas1993simulated}
Dimitris Bertsimas and John Tsitsiklis.
\newblock Simulated annealing.
\newblock {\em Statistical Science}, 8(1):10--15, 1993.

\bibitem[BV05]{BV}
Fran{\c{c}}ois Bolley and C{\'e}dric Villani.
\newblock Weighted {C}sisz{\'a}r-{K}ullback-pinsker inequalities and
  applications to transportation inequalities.
\newblock {\em Annales de la Facult\'e des sciences de Toulouse:
  Math\'{e}matiques}, 14(3):331--352, 2005.

\bibitem[CCA{\etalchar{+}}18]{cheng-nonconvex}
X.~{Cheng}, N.~S. {Chatterji}, Y.~{Abbasi-Yadkori}, P.~L. {Bartlett}, and M.~I.
  {Jordan}.
\newblock {Sharp Convergence Rates for {L}angevin Dynamics in the Nonconvex
  Setting}.
\newblock {\em arXiv preprint arXiv: 1805.01648}, 2018.

\bibitem[CCBJ18]{Cheng}
Xiang Cheng, Niladri~S Chatterji, Peter~L Bartlett, and Michael~I Jordan.
\newblock Underdamped {L}angevin {MCMC}: A non-asymptotic analysis.
\newblock In {\em Proceedings of the 31st Annual Conference on Learning
  Theory}, 2018.

\bibitem[CCG{\etalchar{+}}16]{chen2016bridging}
Changyou Chen, David Carlson, Zhe Gan, Chunyuan Li, and Lawrence Carin.
\newblock Bridging the gap between stochastic gradient {MCMC} and stochastic
  optimization.
\newblock In {\em Proceedings of the 19th International Conference on
  Artificial Intelligence and Statistics (AISTATS)}, pages 1051--1060, 2016.

\bibitem[CCS{\etalchar{+}}17]{Chaudhari}
P~Chaudhari, Anna Choromanska, S~Soatto, Yann LeCun, C~Baldassi, C~Borgs,
  J~Chayes, Levent Sagun, and R~Zecchina.
\newblock Entropy-{SGD}: Biasing gradient descent into wide valleys.
\newblock In {\em International Conference on Learning Representations (ICLR)},
  2017.

\bibitem[CDC15]{carin-2015-langevin-integrators}
Changyou Chen, Nan Ding, and Lawrence Carin.
\newblock On the convergence of stochastic gradient {MCMC} algorithms with
  high-order integrators.
\newblock In {\em Advances in Neural Information Processing Systems (NIPS)},
  pages 2278--2286, 2015.

\bibitem[CDHS18]{carmon18}
Y.~Carmon, J.~Duchi, O.~Hinder, and A.~Sidford.
\newblock Accelerated methods for nonconvex optimization.
\newblock {\em SIAM Journal on Optimization}, 28(2):1751--1772, 2018.

\bibitem[CDT{\etalchar{+}}09]{chapelle2009tighter}
Olivier Chapelle, Chuong~B Do, Choon~H Teo, Quoc~V Le, and Alex~J Smola.
\newblock Tighter bounds for structured estimation.
\newblock In {\em Advances in Neural Information Processing Systems}, pages
  281--288, 2009.

\bibitem[CFG14]{emilyfox-sghmc}
Tianqi Chen, Emily Fox, and Carlos Guestrin.
\newblock Stochastic gradient {H}amiltonian {Monte Carlo}.
\newblock In {\em International Conference on Machine Learning}, pages
  1683--1691, 2014.

\bibitem[CFM{\etalchar{+}}18]{chatterji2018theory}
Niladri~S Chatterji, Nicolas Flammarion, Yi-An Ma, Peter~L Bartlett, and
  Michael~I Jordan.
\newblock On the theory of variance reduction for stochastic gradient {M}onte
  {C}arlo.
\newblock In {\em International Conference on Machine Learning}, pages
  764--773, 2018.

\bibitem[CHJ13]{Cox2013}
Sonja Cox, Martin Hutzenthaler, and Arnulf Jentzen.
\newblock Local {L}ipschitz continuity in the initial value and strong
  completeness for nonlinear stochastic differential equations.
\newblock {\em arXiv preprint arXiv:1309.5595}, 2013.

\bibitem[CHS87]{chiang1987diffusion}
Tzuu-Shuh Chiang, Chii-Ruey Hwang, and Shuenn~Jyi Sheu.
\newblock Diffusion for global optimization in $\mathbb{R}^n$.
\newblock {\em SIAM Journal on Control and Optimization}, 25(3):737--753, 1987.

\bibitem[CSWB06]{collobert2006trading}
Ronan Collobert, Fabian Sinz, Jason Weston, and L{\'e}on Bottou.
\newblock Trading convexity for scalability.
\newblock In {\em Proceedings of the 23rd International Conference on Machine
  Learning}, pages 201--208, 2006.

\bibitem[Dal17]{Dalalyan}
Arnak~S Dalalyan.
\newblock Theoretical guarantees for approximate sampling from smooth and
  log-concave densities.
\newblock {\em Journal of the Royal Statistical Society: Series B (Statistical
  Methodology)}, 79(3):651--676, 2017.

\bibitem[DFB{\etalchar{+}}14]{ding2014bayesian}
Nan Ding, Youhan Fang, Ryan Babbush, Changyou Chen, Robert~D Skeel, and Hartmut
  Neven.
\newblock Bayesian sampling using stochastic gradient thermostats.
\newblock In {\em Advances in Neural Information Processing Systems (NIPS)},
  pages 3203--3211, 2014.

\bibitem[DK19]{DK2017}
Arnak~S. Dalalyan and Avetik~G. Karagulyan.
\newblock User-friendly guarantees for the {L}angevin {M}onte {C}arlo with
  inaccurate gradient.
\newblock {\em Stochastic Processes and their Applications},
  129(12):5278--5311, 2019.

\bibitem[DKPR87]{duane1987hybrid}
Simon Duane, Anthony~D Kennedy, Brian~J Pendleton, and Duncan Roweth.
\newblock Hybrid {M}onte {C}arlo.
\newblock {\em Physics Letters B}, 195(2):216--222, 1987.

\bibitem[DKRD19]{dalalyan2019}
Arnak~S Dalalyan, Avetik Karagulyan, and Lionel Riou-Durand.
\newblock Bounding the error of discretized {L}angevin algorithms for
  non-strongly log-concave targets.
\newblock {\em arXiv:1906.08530}, 2019.

\bibitem[DLT{\etalchar{+}}18]{du2017gradient}
Simon~S Du, Jason~D Lee, Yuandong Tian, Aarti Singh, and Barnabas Poczos.
\newblock Gradient descent learns one-hidden-layer {CNN}: Don't be afraid of
  spurious local minima.
\newblock In {\em International Conference on Machine Learning}, pages
  1339--1348, 2018.

\bibitem[DRD20]{dalalyan2018kinetic}
Arnak~S Dalalyan and Lionel Riou-Durand.
\newblock On sampling from a log-concave density using kinetic {L}angevin
  diffusions.
\newblock {\em Bernoulli}, 26(3):1956--1988, 2020.

\bibitem[DS19]{DS2019}
Susanne Ditlevsen and Adeline Samson.
\newblock Hypoelliptic diffusions: filtering and inference from complete and
  partial observations.
\newblock {\em Journal of the Royal Statistical Society Series B}, 81:361--384,
  2019.

\bibitem[EGZ19]{Eberle}
Andreas Eberle, Arnaud Guillin, and Raphael Zimmer.
\newblock Couplings and quantitative contraction rates for {L}angevin dynamics.
\newblock {\em Annals of Probability}, 47(4):1982--2010, 2019.

\bibitem[FSS18]{foster2018uniform}
Dylan~J Foster, Ayush Sekhari, and Karthik Sridharan.
\newblock Uniform convergence of gradients for non-convex learning and
  optimization.
\newblock In {\em Advances in Neural Information Processing Systems}, pages
  8745--8756, 2018.

\bibitem[GGZ20]{GGZ2}
Xuefeng Gao, Mert G\"{u}rb\"{u}zbalaban, and Lingjiong Zhu.
\newblock Breaking reversibility accelerates {L}angevin dynamics for global
  non-convex optimization.
\newblock In {\em Advances in Neural Information Processing Systems}, 2020.

\bibitem[Gid85]{gidas1985nonstationary}
Basilis Gidas.
\newblock Nonstationary {M}arkov chains and convergence of the annealing
  algorithm.
\newblock {\em Journal of Statistical Physics}, 39(1-2):73--131, 1985.

\bibitem[GL16]{ghadimi2016}
Saeed Ghadimi and Guanghui Lan.
\newblock Accelerated gradient methods for nonconvex nonlinear and stochastic
  programming.
\newblock {\em Mathematical Programming}, 156(1):59--99, 2016.

\bibitem[GLM17]{ge2017learning}
Rong Ge, Jason~D Lee, and Tengyu Ma.
\newblock Learning one-hidden-layer neural networks with landscape design.
\newblock {\em arXiv preprint arXiv:1711.00501}, 2017.

\bibitem[GM91]{gelfand1991recursive}
Saul~B Gelfand and Sanjoy~K Mitter.
\newblock Recursive stochastic algorithms for global optimization in
  $\mathbb{R}^d$.
\newblock {\em SIAM Journal on Control and Optimization}, 29(5):999--1018,
  1991.

\bibitem[Haj85]{hajek1985tutorial}
Bruce Hajek.
\newblock A tutorial survey of theory and applications of simulated annealing.
\newblock In {\em 1985 24th IEEE Conference on Decision and Control}, pages
  755--760. IEEE, 1985.

\bibitem[Hal88]{hale1988asymptotic}
JK~Hale.
\newblock {\em Asymptotic Behavior of Dissipative Systems}, volume~25.
\newblock American Mathematical Society, 1988.

\bibitem[HK70]{hoerl1970ridge}
Arthur~E Hoerl and Robert~W Kennard.
\newblock Ridge regression: Biased estimation for nonorthogonal problems.
\newblock {\em Technometrics}, 12(1):55--67, 1970.

\bibitem[HKS89]{stroock-langevin-spectrum}
Richard~A Holley, Shigeo Kusuoka, and Daniel~W Stroock.
\newblock Asymptotics of the spectral gap with applications to the theory of
  simulated annealing.
\newblock {\em Journal of Functional Analysis}, 83(2):333--347, 1989.

\bibitem[HLM{\etalchar{+}}17]{hu2017group}
Yaohua Hu, Chong Li, Kaiwen Meng, Jing Qin, and Xiaoqi Yang.
\newblock Group sparse optimization via $\ell_{p,q}$ regularization.
\newblock {\em The Journal of Machine Learning Research}, 18(1):960--1011,
  2017.

\bibitem[HLSS16]{hazan2016graduated}
Elad Hazan, Kfir~Yehuda Levy, and Shai Shalev-Shwartz.
\newblock On graduated optimization for stochastic non-convex problems.
\newblock In {\em International Conference on Machine Learning}, pages
  1833--1841, 2016.

\bibitem[HN04]{herau-nier-underdamped}
Fr{\'e}d{\'e}ric H{\'e}rau and Francis Nier.
\newblock Isotropic hypoellipticity and trend to equilibrium for the
  {F}okker-{P}lanck equation with a high-degree potential.
\newblock {\em Archive for Rational Mechanics and Analysis}, 171(2):151--218,
  2004.

\bibitem[HRS16]{hardt2016train}
Moritz Hardt, Benjamin Recht, and Yoram Singer.
\newblock Train faster, generalize better: stability of stochastic gradient
  descent.
\newblock In {\em International Conference on Machine Learning}, pages
  1225--1234, 2016.

\bibitem[JT19]{JT2019}
A.~Jofr\'{e} and P.~Thompson.
\newblock On variance reduction for stochastic smooth convex optimization with
  multiplicative noise.
\newblock {\em Mathematical Programming}, 174:253--292, 2019.

\bibitem[JWHT13]{james2013introduction}
Gareth James, Daniela Witten, Trevor Hastie, and Robert Tibshirani.
\newblock {\em An Introduction to Statistical Learning}, volume 112.
\newblock Springer, 2013.

\bibitem[KGV83]{kirkpatrick1983optimization}
Scott Kirkpatrick, C~Daniel Gelatt, and Mario~P Vecchi.
\newblock Optimization by simulated annealing.
\newblock {\em Science}, 220(4598):671--680, 1983.

\bibitem[KO01]{kilmer2001choosing}
Misha~E Kilmer and Dianne~P O'Leary.
\newblock Choosing regularization parameters in iterative methods for ill-posed
  problems.
\newblock {\em SIAM Journal on matrix analysis and applications},
  22(4):1204--1221, 2001.

\bibitem[Kra40]{kramers1940brownian}
Hendrik~Anthony Kramers.
\newblock Brownian motion in a field of force and the diffusion model of
  chemical reactions.
\newblock {\em Physica}, 7(4):284--304, 1940.

\bibitem[LCZZ18]{liu2018toward}
Tianyi Liu, Zhehui Chen, Enlu Zhou, and Tuo Zhao.
\newblock Toward deeper understanding of nonconvex stochastic optimization with
  momentum using diffusion approximations.
\newblock {\em arXiv preprint arXiv:1802.05155}, 2018.

\bibitem[LFM18]{lu2018accelerating}
Haihao Lu, Robert~M Freund, and Vahab Mirrokni.
\newblock Accelerating greedy coordinate descent methods.
\newblock In {\em International Conference on Machine Learning}, pages
  3257--3266, 2018.

\bibitem[LMS15]{leimkuhler2015computation}
Benedict Leimkuhler, Charles Matthews, and Gabriel Stoltz.
\newblock The computation of averages from equilibrium and nonequilibrium
  {L}angevin molecular dynamics.
\newblock {\em IMA Journal of Numerical Analysis}, 36(1):13--79, 2015.

\bibitem[LS13]{liptser2013statistics}
Robert~S Liptser and Albert~N Shiryaev.
\newblock {\em Statistics of Random Processes: {I}. General Theory}, volume~5.
\newblock Springer Science \& Business Media, 2013.

\bibitem[LV18]{lee2018convergence}
Yin~Tat Lee and Santosh~S Vempala.
\newblock Convergence rate of {R}iemannian {H}amiltonian {M}onte {C}arlo and
  faster polytope volume computation.
\newblock In {\em Proceedings of the 50th Annual ACM SIGACT Symposium on Theory
  of Computing}, pages 1115--1121. ACM, 2018.

\bibitem[MBM18]{mei2018landscape}
Song Mei, Yu~Bai, and Andrea Montanari.
\newblock The landscape of empirical risk for nonconvex losses.
\newblock {\em The Annals of Statistics}, 46(6A):2747--2774, 2018.

\bibitem[MCF15]{ma2015complete}
Yi-An Ma, Tianqi Chen, and Emily Fox.
\newblock A complete recipe for stochastic gradient {MCMC}.
\newblock In {\em Advances in Neural Information Processing Systems (NIPS)},
  pages 2917--2925, 2015.

\bibitem[MPS18]{MPS2018}
Oren Mangoubi, S.~Pillai, Natesh, and Aaron Smith.
\newblock Does {H}amiltonian {M}onte {C}arlo mix faster than a random walk on
  multimodal densities?
\newblock {\em arXiv preprint arXiv:1808.03230}, 2018.

\bibitem[MS17]{Mangoubi-Smith17}
O.~{Mangoubi} and A.~{Smith}.
\newblock {Rapid Mixing of {H}amiltonian {M}onte {C}arlo on Strongly
  Log-Concave Distributions}.
\newblock {\em arXiv preprint arXiv: 1708.07114}, 2017.

\bibitem[MSH02]{mattingly2002ergodicity}
Jonathan~C Mattingly, Andrew~M Stuart, and Desmond~J Higham.
\newblock Ergodicity for {SDE}s and approximations: locally {L}ipschitz vector
  fields and degenerate noise.
\newblock {\em Stochastic Processes and their Applications}, 101(2):185--232,
  2002.

\bibitem[Nea10]{neal2010mcmc}
RM~Neal.
\newblock {MCMC} using {H}amiltonian dynamics. {Handbook of {M}arkov {C}hain
  {M}onte {C}arlo} ({S. Brooks, A. Gelman, G. Jones, and X.-L. Meng, eds.}),
  2010.

\bibitem[Nes83]{nesterov1983method}
Yurii~E Nesterov.
\newblock A method for solving the convex programming problem with convergence
  rate $o(1/k^2)$.
\newblock In {\em Dokl. Akad. Nauk SSSR}, volume 269, pages 543--547, 1983.

\bibitem[NS13]{nguyen2013algorithms}
Tan Nguyen and Scott Sanner.
\newblock Algorithms for direct 0--1 loss optimization in binary
  classification.
\newblock In {\em International Conference on Machine Learning}, pages
  1085--1093, 2013.

\bibitem[{\O}ks03]{Oksendal}
B.~K. {\O}ksendal.
\newblock {\em Stochastic Differential Equations: An Introduction with
  Applications}.
\newblock Springer, 2003.

\bibitem[OW19]{wright2017escape}
Michael O'Neill and Stephen~J Wright.
\newblock Behavior of accelerated gradient methods near critical points of
  nonconvex problems.
\newblock {\em Mathematical Programming}, 176:403--427, 2019.

\bibitem[Pav14]{pavliotis2014stochastic}
Grigorios~A Pavliotis.
\newblock {\em Stochastic Processes and Applications: Diffusion Processes, the
  {F}okker-{P}lanck and {L}angevin Equations}, volume~60.
\newblock Springer, 2014.

\bibitem[Pol87]{polyak1987introduction}
Boris~T Polyak.
\newblock Introduction to {O}ptimization.
\newblock {\em {O}ptimization {S}oftware. New York}, 1987.

\bibitem[PT13]{patterson-teh}
Sam Patterson and Yee~Whye Teh.
\newblock Stochastic gradient {R}iemannian {L}angevin dynamics on the
  probability simplex.
\newblock In {\em Advances in Neural Information Processing Systems (NIPS)},
  pages 3102--3110, 2013.

\bibitem[PW16]{polyanskiy2016wasserstein}
Yury Polyanskiy and Yihong Wu.
\newblock {W}asserstein continuity of entropy and outer bounds for interference
  channels.
\newblock {\em IEEE Transactions on Information Theory}, 62(7):3992--4002,
  2016.

\bibitem[ROK12]{razavi2012robust}
S~Alireza Razavi, Esa Ollila, and Visa Koivunen.
\newblock Robust greedy algorithms for compressed sensing.
\newblock In {\em 2012 Proceedings of the 20th European Signal Processing
  Conference (EUSIPCO)}, pages 969--973. IEEE, 2012.

\bibitem[RRT17]{Raginsky}
Maxim Raginsky, Alexander Rakhlin, and Matus Telgarsky.
\newblock Non-convex learning via stochastic gradient {L}angevin dynamics: a
  nonasymptotic analysis.
\newblock In {\em Conference on Learning Theory}, pages 1674--1703, 2017.

\bibitem[SBC14]{su2014differential}
Weijie Su, Stephen Boyd, and Emmanuel Candes.
\newblock A differential equation for modeling {N}esterov's accelerated
  gradient method: Theory and insights.
\newblock In {\em Advances in Neural Information Processing Systems}, pages
  2510--2518, 2014.

\bibitem[SDJS18]{shi2018understanding}
Bin Shi, Simon~S Du, Michael~I Jordan, and Weijie~J Su.
\newblock Understanding the acceleration phenomenon via high-resolution
  differential equations.
\newblock {\em arXiv preprint arXiv:1810.08907}, 2018.

\bibitem[{\c{S}im\c{s}ekli}BCR16]{simsekli2016stochastic}
Umut {\c{S}im\c{s}ekli}, Roland Badeau, Taylan Cemgil, and Ga\"{e}l Richard.
\newblock Stochastic quasi-{N}ewton {L}angevin {M}onte {C}arlo.
\newblock In {\em International Conference on Machine Learning}, pages
  642--651, 2016.

\bibitem[{\c{S}im\c{s}ekli}YN{\etalchar{+}}18]{simsekli-async-MCMC-18}
U.~{\c{S}im\c{s}ekli}, {\c C}.~{Y{\i}ld{\i}z}, T.~H. {Nguyen}, G.~{Richard},
  and A.~{Taylan Cemgil}.
\newblock {Asynchronous Stochastic Quasi-Newton {MCMC} for Non-Convex
  Optimization}.
\newblock In {\em International Conference on Machine Learning}, pages
  4674--4683, 2018.

\bibitem[SMDH13]{sutskever2013importance}
Ilya Sutskever, James Martens, George Dahl, and Geoffrey Hinton.
\newblock On the importance of initialization and momentum in deep learning.
\newblock In {\em International Conference on Machine Learning}, pages
  1139--1147, 2013.

\bibitem[SQW16]{sun2016complete}
Ju~Sun, Qing Qu, and John Wright.
\newblock Complete dictionary recovery over the sphere {I}: Overview and the
  geometric picture.
\newblock {\em IEEE Transactions on Information Theory}, 63(2):853--884, 2016.

\bibitem[TLR18]{pmlr-v75-tzen18a}
Belinda Tzen, Tengyuan Liang, and Maxim Raginsky.
\newblock Local optimality and generalization guarantees for the {L}angevin
  algorithm via empirical metastability.
\newblock In {\em Conference on Learning Theory}, pages 857--875, 2018.

\bibitem[TTV16]{teh2016consistency}
Yee~Whye Teh, Alexandre~H Thiery, and Sebastian~J Vollmer.
\newblock Consistency and fluctuations for stochastic gradient {L}angevin
  dynamics.
\newblock {\em The Journal of Machine Learning Research}, 17(1):193--225, 2016.

\bibitem[Vil08]{villani2008optimal}
C{\'e}dric Villani.
\newblock {\em Optimal Transport: Old and New}, volume 338.
\newblock Springer Science \& Business Media, 2008.

\bibitem[WCX19]{wang2019differentially}
Di~Wang, Changyou Chen, and Jinhui Xu.
\newblock Differentially private empirical risk minimization with non-convex
  loss functions.
\newblock In {\em International Conference on Machine Learning}, pages
  6526--6535, 2019.

\bibitem[Wib18]{wibisono2018sampling}
Andre Wibisono.
\newblock Sampling as optimization in the space of measures: The {L}angevin
  dynamics as a composite optimization problem.
\newblock In {\em Proceedings of the 31st Annual Conference on Learning
  Theory}, 2018.

\bibitem[WJHZ13]{exp-loss}
Xueqin Wang, Yunlu Jiang, Mian Huang, and Heping Zhang.
\newblock Robust variable selection with exponential squared loss.
\newblock {\em Journal of the American Statistical Association},
  108(502):632--643, 2013.

\bibitem[WL07]{wu2007robust}
Yichao Wu and Yufeng Liu.
\newblock Robust truncated hinge loss support vector machines.
\newblock {\em Journal of the American Statistical Association},
  102(479):974--983, 2007.

\bibitem[WRJ16]{wilson2016lyapunov}
Ashia~C Wilson, Benjamin Recht, and Michael~I Jordan.
\newblock A lyapunov analysis of momentum methods in optimization.
\newblock {\em arXiv preprint arXiv:1611.02635}, 2016.

\bibitem[WT11]{welling2011bayesian}
Max Welling and Yee~W Teh.
\newblock Bayesian learning via stochastic gradient {L}angevin dynamics.
\newblock In {\em International Conference on Machine Learning}, pages
  681--688, 2011.

\bibitem[XCZG18]{xu2018global}
Pan Xu, Jinghui Chen, Difan Zou, and Quanquan Gu.
\newblock Global convergence of {L}angevin dynamics based algorithms for
  nonconvex optimization.
\newblock In {\em Advances in Neural Information Processing Systems}, pages
  3122--3133, 2018.

\bibitem[ZLC17]{zhang-sgld}
Yuchen {Zhang}, Percy. {Liang}, and Moses. {Charikar}.
\newblock A hitting time analysis of stochastic gradient {L}angevin dynamics.
\newblock In {\em Conference on Learning Theory}, pages 1980--2022, 2017.

\bibitem[ZZLC17]{zhang2017nonconvex}
Huishuai Zhang, Yi~Zhou, Yingbin Liang, and Yuejie Chi.
\newblock A nonconvex approach for phase retrieval: Reshaped {W}irtinger flow
  and incremental algorithms.
\newblock {\em The Journal of Machine Learning Research}, 18(1):5164--5198,
  2017.

\end{thebibliography}

\newpage
\appendix

\section{Proof of Theorem~\ref{thm-first-term} 
and Corollary~\ref{coro:main}}\label{sec:proof}

We first present several technical lemmas that will be used in our analysis and review existing results for the underdamped Langevin SDE. The proof of these lemmas will be deferred to Section~\ref{sec:lemma-proof-1}. 

Our analysis for analyzing the convergence speed of the SGHMC1 algorithm and its comparison to the underdamped Langevin SDE is based on the 2-Wasserstein distance and this requires the $L^2$ norm of the iterates to be finite. In the next lemma, we show that $L^2$ norm of the both discrete and continuous dynamics are uniformly bounded over time with explicit constants. The main idea is to make use of the properties of the Lyapunov function $\mathcal{V}$ which is designed originally for the continuous-time process and show that the discrete dynamics can also be controlled by it. 


\begin{lemma}[Uniform $L^{2}$ bounds]\label{lem:L2bound} 
\begin{enumerate}
\item [$(i)$] It holds that
\begin{align}
&\sup_{t\geq 0}\mathbb{E}_{\mathbf{z}}\Vert X(t)\Vert^{2}\leq
C_{x}^{c}
:=\frac{\int_{\mathbb{R}^{2d}}\mathcal{V}(x,v)d\mu_{0}(x,v)+\frac{d+A}{\lambda}}{\frac{1}{8} (1-2 \lambda) \beta \gamma^2}<\infty, \label{def-Cx-c}
\\
&\sup_{t\geq 0}\mathbb{E}_{\mathbf{z}}\Vert V(t)\Vert^{2}\leq C_{v}^{c}
:=\frac{\int_{\mathbb{R}^{2d}}\mathcal{V}(x,v)d\mu_{0}(x,v)
+\frac{d+A}{\lambda}}{\frac{\beta}{4}(1-2\lambda)}<\infty.
\label{def-Cxv-c}
\end{align}
\item [$(ii)$] For
$0<\eta\leq\min\left\{\frac{\gamma}{K_{2}}(d/\beta+A/\beta),\frac{\gamma\lambda}{2K_{1}},\frac{2}{\gamma\lambda}\right\}$,
where
\begin{align}
&K_{1}
:=\max\left\{\frac{32M^{2}\left(\frac{1}{2}+\gamma+\delta\right)}{ (1-2 \lambda) \beta \gamma^2 },
\frac{8\left(\frac{1}{2}M+\frac{1}{4}\gamma^{2}-\frac{1}{4}\gamma^{2}\lambda+\gamma\right)}{\beta(1-2\lambda)}\right\},
\label{def-K1}
\\
&K_{2}:=2B^{2}\left(\frac{1}{2}+\gamma+\delta\right),\label{def-K2}
\end{align}
we have
\begin{align}
&\sup_{j\geq 0}\mathbb{E}_{\mathbf{z}}\Vert X_{j}\Vert^{2}
\leq
C_{x}^{d}
:=\frac{\int_{\mathbb{R}^{2d}}\mathcal{V}(x,v)\mu_{0}(dx,dv)
+\frac{4(d+A)}{\lambda}}{\frac{1}{8} (1-2 \lambda) \beta \gamma^2 }<\infty,\label{def-Cx-d}
\\
&\sup_{j\geq 0}\mathbb{E}_{\mathbf{z}}\Vert V_{j}\Vert^{2}
\leq
C_{v}^{d}
:=\frac{\int_{\mathbb{R}^{2d}}\mathcal{V}(x,v)\mu_{0}(dx,dv)
+\frac{4(d+A)}{\lambda}}{\frac{\beta}{4}(1-2\lambda)}
<\infty.\label{def-Cxv-d}
\end{align}
\end{enumerate}
\end{lemma}

Since SGHMC1 is a discretization of the underdamped SDE (except that noise is also added to the gradients), we expect SGHMC1 to follow the underdamped SDE dynamics. It is natural to seek for bounds between the probability law $\mu_{\bz,k}$ of the SGHMC1 algorithm at step $k$ with time step $\eta$ and that of the underdamped SDE at time $t=k\eta$ which we denote by $\nu_{{\mathbf{z}},k\eta}$. In our analysis, we first control the Kullback-Leibler (KL) divergence between these two, and then convert these bounds into bounds in terms of the 2-Wasserstein metric, applying an optimal transportation inequality by \cite{BV}. {Note that Bolley and Villani theorem has also been successfully applied to analyzing the SGLD dynamics in \cite{Raginsky}. However, the analysis in \cite{Raginsky} does not directly apply to our setting as underdamped dynamics require a different Lyapunov function.} This step requires an exponential integrability property of the underdamped SDE process which we establish next, before stating our result in Lemma~\ref{diff-app} about the diffusion approximation of the SGHMC1 iterates.
\begin{lemma}[Exponential integrability]\label{lem:exp-integrable}
For every $t$,
\begin{equation*}
\mathbb{E}_{\mathbf{z}}\left[e^{\alpha_{0}\Vert(X(t),V(t))\Vert^{2}}\right]
\leq
\int_{\mathbb{R}^{2d}}e^{\frac{1}{4}\alpha\mathcal{V}(x,v)}\mu_{0}(dx,dv)
+\frac{1}{4}e^{\frac{\alpha(d+A)}{3\lambda}}\alpha\gamma(d+A)t,
\end{equation*}
where
\begin{equation}\label{eq:alpha0alpha}
\alpha_{0}:=\frac{\alpha}{\frac{64}{(1-2\lambda)\beta\gamma^{2}}+\frac{32}{\beta(1-2\lambda)}},
\qquad
\alpha:=\frac{\lambda(1-2\lambda)}{12}.
\end{equation}
\end{lemma}

We showed in the above Lemma~\ref{lem:exp-integrable} that the exponential moments grow at most linearly
in time $t$, which is a strict improvement from the exponential growth in time $t$ in \cite{Raginsky}.
As a result, in the following Lemma~\ref{diff-app} for the diffusion approximation, 
our upper bound is of the order $(k\eta)^{3/2}\sqrt{\log(k\eta)}(\delta^{1/4}+\eta^{1/4})+k\eta\sqrt{\eta}$
compared to $k\eta(\delta^{1/4}+\eta^{1/4})$ in \cite{Raginsky}. 
{The method that is used in the proof of Lemma~\ref{lem:exp-integrable} for the underdamped
dynamics can indeed be adapted to the case of the overdamped dynamics to improve the results in \cite{Raginsky}.}

\begin{lemma}[Diffusion approximation] \label{diff-app}
For any $k\in\mathbb{N}$ and any $\eta\leq 1$, so that $k\eta\geq e$
and $\eta$ satisfies the condition in Part (ii) of Lemma~\ref{lem:L2bound}.
Then, we have
\begin{equation*}
\mathcal{W}_{2}(\mu_{{\mathbf{z}},k},\nu_{{\mathbf{z}},k\eta})
\leq
(C_{0}\delta^{1/4}+C_{1}\eta^{1/4})\cdot(k\eta)^{3/2}\cdot\sqrt{\log(k\eta)}
+C_{2}(k\eta)\sqrt{\eta},
\end{equation*}
where $C_{0}$, $C_{1}$ and $C_{2}$ are given by:
\begin{align}
&C_{0}=\hat{\gamma}
\cdot\left(\left(M^{2}C_{x}^{d}+B^{2}\right)\frac{\beta}{\gamma}+\sqrt{\left(M^{2}C_{x}^{d}+B^{2}\right)\frac{\beta}{\gamma}}\right)^{1/2},\label{def-hat-Czero}
\\
&C_{1}=\hat{\gamma}
\cdot\left(\left(\frac{M^{2}\beta\eta}{\gamma}+\frac{\beta\eta\gamma}{2}\right)(C_{2})^{2}
+\sqrt{\left(\frac{M^{2}\beta\eta}{\gamma}+\frac{\beta\eta\gamma}{2}\right)(C_{2})^{2}}\right)^{1/2},
\label{def-hat-Cone}
\\
&C_{2}=\left(2\gamma^{2}C_{v}^{d}
+(4+2\delta)\left(M^{2}C_{x}^{d}+B^{2}\right)
+2\gamma\beta^{-1}\right)^{1/2},
\label{def-hat-Ctwo}
\\
&\hat{\gamma}=\frac{2\sqrt{2}}{\sqrt{\alpha_{0}}}\left(\frac{5}{2}+\log\left(\int_{\mathbb{R}^{2d}}e^{\frac{1}{4}\alpha\mathcal{V}(x,v)}\mu_{0}(dx,dv)
+\frac{1}{4}e^{\frac{\alpha(d+A)}{3\lambda}}\alpha\gamma(d+A)\right)
\right)^{1/2}.
\label{def-hat-gamma}
\end{align}
\end{lemma}


\subsection{Convergence rate to the equilibrium of the underdamped SDE}\label{subsec-Eberle}
We consider the underdamped SDE and bound the 2-Wasserstein distance $\mathcal{W}_{2}(\nu_{z, t},\pi_{\mathbf{z}})$ to the equilibrium for a fix arbitrary time $t \ge 0$. Crucial to the analysis is \cite{Eberle}, which quantifies the convergence to equilibrium for underdamped Langevin diffusions. We first review the results from \cite{Eberle}. Let us recall from \eqref{eq:lyapunov} the definition of the Lyapunov function $\mathcal{V}(x,v)$:
\begin{equation*}
\mathcal{V}(x,v)=\beta F_{\mathbf{z}}(x)
+\frac{\beta}{4}\gamma^{2}(\Vert x+\gamma^{-1}v\Vert^{2}+\Vert\gamma^{-1}v\Vert^{2}-\lambda\Vert x\Vert^{2}).
\end{equation*}
For any $(x,v),(x',v')\in\mathbb{R}^{2d}$, we set:
\begin{align*}
&r((x,v),(x',v'))=\alpha_{1}\Vert x-x'\Vert+\Vert x-x'+\gamma^{-1}(v-v')\Vert,
\\
&\rho((x,v),(x',v'))=h(r((x,v),(x',v')))\cdot(1+\varepsilon_{1}\mathcal{V}(x,v)+\varepsilon_{1}\mathcal{V}(x',v')),
\end{align*}
where $\alpha_{1},\varepsilon_{1}>0$ are appropriately chosen constants,
and $h:[0,\infty)\rightarrow[0,\infty)$ is continuous, non-decreasing concave
function such that $h(0)=0$,
$h$ is $C^{2}$ on $(0,R_{1})$ for some constant $R_{1}>0$ with right-sided derivative $h'_{+}(0)=1$
and left-sided derivative $h'_{-}(R_{1})>0$ and $h$ is constant on $[R_{1},\infty)$. 
For any two probability measures $\mu,\nu$ on $\mathbb{R}^{2d}$, we define
\begin{equation}\label{def-H-rho}
\mathcal{H}_{\rho}(\mu,\nu):=
\inf_{(X,V)\sim\mu,(X',V')\sim\nu}\mathbb{E}[\rho((X,V),(X',V'))].
\end{equation}
Note that $\mathcal{H}_{\rho}$ is a semi-metric, but not necessarily a metric. A simplified version of the main result from \cite{Eberle} which will be used in our setting is given below. 

\begin{theorem}[Theorem 2.3. and Corollary 2.6. in \cite{Eberle}] \label{exp-converg} There exist constants $\alpha_{1}, \epsilon_1 \in (0, \infty)$ and a continuous non-decreasing function $h: [0, \infty) \rightarrow [0, \infty)$ with $h(0)=0$ such that we have
\begin{equation*}
\mathcal{W}_{2}(\nu_{{\mathbf{z}},k\eta},\pi_{\mathbf{z}})
\leq C\sqrt{\mathcal{H}_{\rho}(\mu_{0},\pi_{\mathbf{z}})}e^{-\mu_{\ast}k\eta}\,,
\end{equation*}
where
\begin{align}
&\mu_{\ast}=\frac{\gamma}{768}\min\{\lambda M\gamma^{-2},\Lambda^{1/2}e^{-\Lambda}M\gamma^{-2},\Lambda^{1/2}e^{-\Lambda}\}\label{def-mu-star},
\\
&C=\sqrt{2}e^{1+\frac{\Lambda}{2}}\frac{1+\gamma}{\min\{1,\alpha_{1}\}}
\sqrt{\max\{1,4(1+2\alpha_{1}+2\alpha_{1}^{2})(d+A)\beta^{-1}\gamma^{-1}\mu_{\ast}^{-1}/\min\{1,R_{1}\}\}}\label{def-capital-C},
\\
&\Lambda=\frac{12}{5}(1+2\alpha_{1}+2\alpha_{1}^{2})(d+A)M\gamma^{-2}\lambda^{-1}(1-2\lambda)^{-1},
\qquad\alpha_{1}=(1+\Lambda^{-1})M\gamma^{-2}\label{def-capital-Lambda-alpha_1},
\\
&\varepsilon_{1}=4\gamma^{-1}\mu_{\ast}/(d+A)
\label{def-eps_1},
\\
&R_{1}=4\cdot(6/5)^{1/2}(1+2\alpha_{1}+2\alpha_{1}^{2})^{1/2}(d+A)^{1/2} \beta^{-1/2}\gamma^{-1}(\lambda-2\lambda^{2})^{-1/2}\label{def-R1}.
\end{align}
\end{theorem}

We remark that the definitions of $\Lambda, \alpha_1$
in \eqref{def-capital-Lambda-alpha_1} are coupled
and there exists $\alpha_{1}\in(0,\infty)$ so that
$\Lambda, \alpha_1$ in \eqref{def-capital-Lambda-alpha_1} are
well defined; see Theorem 2.3. in \cite{Eberle}.
In order to get explicit performance bounds, we also derive an upper bound for $\mathcal{H}_{\rho}(\mu_{0},\pi_{\mathbf{z}})$ in the next lemma. It is based on the (integrability properties) structure of the stationary distribution $\pi_\bz$ and the Lyapunov function $\mathcal{V}$ that controls the $L^2$ norm of the initial distribution $\mu_0$.

\begin{lemma}[Bounding initialization error]\label{lem:nu0-piz}
If parts $(i)$, $(ii)$, $(iii)$ and $(iv)$ of Assumption~\ref{assumptions} hold, then we have
\begin{align}
\mathcal{H}_{\rho}(\mu_{0},\pi_{\mathbf{z}})
&\leq \overline{\mathcal{H}}_{\rho}(\mu_{0}):=
R_{1}+R_{1}\varepsilon_{1}\max\left\{M+\frac{1}{2} \beta \gamma^{2},\frac{3}{4} \beta \right\}\Vert(x,v)\Vert_{L^{2}(\mu_{0})}^{2}
\nonumber
\\
&\qquad
+R_{1}\varepsilon_{1}\left(M+\frac{1}{2} \beta \gamma^{2}\right) \frac{b+d/\beta}{m}
+R_{1}\varepsilon_{1}\frac{3}{4} d
+2R_{1}\varepsilon_{1}\left(\beta A_{0}+\frac{\beta B^{2}}{2M}\right), \label{ineq-H-rho-upper}
\end{align}
where $\Vert(x,v)\Vert_{L^{2}(\mu_{0})}^{2}:=\int_{\mathbb{R}^{2d}}\Vert(x,v)\Vert^{2}\mu_{0}(dx,dv)$.
\end{lemma}

\subsection{Proof of Theorem~\ref{thm-first-term}}

As the function $F_{\bz}$ satisfies the conditions in Lemma~\ref{lem:expec-wass} in Section~\ref{sec:support} with $c_1 = M$ and $c_2= B$  (Lemma~\ref{lem:gradient-bound} in Section~\ref{sec:support}), and the probability measures
$\mu_{k,{\mathbf{z}}}, \pi_{\mathbf{z}}$ have finite second moments (Lemma~\ref{lem:L2bound}), we can apply Lemma~\ref{lem:expec-wass} and deduce that
\begin{equation}\label{F:mu:pi}
 \left| \int_{\mathbb{R}^{d}\times\mathbb{R}^{d}} F_{\bz}(x) \mu_{k,{\mathbf{z}}}(dx,dv) - \int_{\mathbb{R}^{d}\times\mathbb{R}^{d}} F_{\bz}(x) \pi_{\mathbf{z}}(dx,dv) \right| \le (M \sigma + B) \cdot \mathcal{W}_{2}(\mu_{{\mathbf{z}},k} ,\pi_{\mathbf{z}}).
\end{equation}
Here, one can obtain from Lemma~\ref{lem:L2bound} and Theorem~\ref{exp-converg} (convergence in 2-Wasserstein distance implies convergence of second moments) that
\begin{equation}\label{defn:sigma}
 \sigma^2 = \max\left\{C_{x}^{c},C_{x}^{d}\right\}=C_{x}^{d}.
\end{equation}
Now, by Lemma~\ref{diff-app} and Theorem~\ref{exp-converg}, we have
\begin{align*} \label{eq:w-d}
\mathcal{W}_{2}(\mu_{{\mathbf{z}},k} ,\pi_{\mathbf{z}})
&\leq \mathcal{W}_{2}(\mu_{{\mathbf{z}},k},\nu_{{\mathbf{z}},k\eta}) + \mathcal{W}_{2}(\nu_{{\mathbf{z}},k\eta}, \pi_{\bz}) \\
&\leq (C_{0}\delta^{1/4}+C_{1}\eta^{1/4})\cdot(k\eta)^{3/2}\cdot\sqrt{\log(k\eta)} 
+C_{2}(k\eta)\sqrt{\eta}
+ C\sqrt{\mathcal{H}_{\rho}(\mu_{0},\pi_{\mathbf{z}})}e^{-\mu_* k\eta}.
\end{align*}
It then follows from \eqref{F:mu:pi} that
\begin{align*}
&\left| \int_{\mathbb{R}^{d}\times\mathbb{R}^{d}} F_{\bz}(x) \mu_{k, \bz}(dx,dv) - \int_{\mathbb{R}^{d}\times\mathbb{R}^{d}} F_{\bz}(x) \pi_{\bz}(dx,dv) \right|
\\
&\leq (M \sigma +B) \cdot \left( C\sqrt{\mathcal{H}_{\rho}(\mu_{0},\pi_{\mathbf{z}})}e^{-\mu_* k\eta}
+(C_{0}\delta^{1/4}+C_{1}\eta^{1/4})\cdot(k\eta)^{3/2}\cdot\sqrt{\log(k\eta)}
+C_{2}(k\eta)\sqrt{\eta}\right).
\end{align*}
Let $k\eta\geq e$, and
\begin{equation*}
k\eta=\frac{1}{\mu_{\ast}}\log\left(\frac{1}{\varepsilon}\right).
\end{equation*}
Then for any $\eta$ satisfying the condition in Lemma~\ref{lem:L2bound} and
$\eta\leq\left(\frac{\varepsilon}{(\log(1/\varepsilon))^{3/2}}\right)^{4}$,
we have
\begin{align*}
&\left| \int_{\mathbb{R}^{d}\times\mathbb{R}^{d}} F_{\bz}(x) \mu_{k,\bz}(dx,dv) - \int_{\mathbb{R}^{d}\times\mathbb{R}^{d}} F_{\bz}(x) \pi_{\bz}(dx,dv) \right|
\\
&\leq
(M \sigma +B) \cdot \Bigg( C\sqrt{\mathcal{H}_{\rho}(\mu_{0},\pi_{\mathbf{z}})}\varepsilon
 + \left(\frac{C_{0}}{\mu_{\ast}^{3/2}}(\log(1/\varepsilon))^{3/2}\delta^{1/4}
 +\frac{C_{1}}{\mu_{\ast}^{3/2}}\varepsilon
 \right)\sqrt{\log(\mu_{\ast}^{-1}\log(\varepsilon^{-1}))}
 \\
 &\qquad\qquad\qquad\qquad
 +\frac{C_{2}}{\mu_{\ast}}\frac{\epsilon^{2}}{(\log(1/\varepsilon))^{2}}\Bigg).
\end{align*}
 The proof is therefore complete.

\subsection{Proof of Corollary~\ref{coro:main}}

With a slight abuse of notations, consider the random elements $(\hat X, \hat V)$ and $(\hat{X}^*, \hat{V}^*)$ with $\text{Law}((\hat X, \hat V) | \mathbf{Z}={\mathbf{z}}) = \mu_{{\mathbf{z}}, k}$ and $\text{Law}((\hat{X}^*, \hat{V}^*) | \mathbf{Z}={\mathbf{z}}) = \pi_{\mathbf{z}}$.
Then we can decompose the expected population risk of $\hat X$ (which has the same distribution as $X_k$) as follows:
\begin{equation} \label{eq:decomp}
\mathbb{E}F(\hat{X})-F^{\ast} = \left( \mathbb{E}F(\hat{X}) - \mathbb{E}F(\hat{X}^*) \right) +  \left( \mathbb{E}F(\hat{X}^*) -  \mathbb{E}F_\mathbf{Z}(\hat{X}^*) \right) + \left( \mathbb{E}F_{\mathbf{Z}}(\hat{X}^*) - F^* \right)\,.
\end{equation}

The first term in \eqref{eq:decomp} can be
written as:
\begin{equation*}
 \mathbb{E}F(\hat{X}) - \mathbb{E}F(\hat{X}^*)  = \int_{\mathcal{Z}^n} P^n(d{\mathbf{z}}) \left( \int_{\mathbb{R}^{2d}} F_{\bz}(x) \mu_{k,{\mathbf{z}}}(dx,dv) - \int_{\mathbb{R}^{2d}} F_{\bz}(x) \pi_{\mathbf{z}}(dx,dv) \right),
\end{equation*}
where $P^n$ is the product measure of independent random variables $Z_1,\ldots,Z_n$. Then it follows from Theorem~\ref{thm-first-term} and Lemma~\ref{lem:nu0-piz} that
\begin{align*}
&\mathbb{E}F(\hat{X}) - \mathbb{E}F(\hat{X}^*) \leq \overline{\mathcal{J}}_0(\varepsilon) + \mathcal{J}_1(\varepsilon).
\end{align*}


Next, we bound the second and third terms in \eqref{eq:decomp}.
Note that
\begin{align*}
\int_{\mathbb{R}^{2d}}F_{\mathbf{z}}(x)\pi_{\mathbf{z}}(dx,dv)
=\int_{\mathbb{R}^{d}}F_{\mathbf{z}}(x)\pi_{\mathbf{z}}(dx),
\end{align*}
where $\pi_{\mathbf{z}}(dx)=\Lambda_{\mathbf{z}} e^{-\beta F_{\mathbf{z}}(x)}dx $ and $\Lambda_{\mathbf{z}}=\int_{\mathbb{R}^{d}}e^{-\beta F_{\mathbf{z}}(x)}dx$. The distribution $\pi_{\bz}(dx)$, i.e., the $x-$marginal of $\pi_{\bz}(dx, dv)$, is the same as the stationary distribution of the overdamped Langevin SDE in \eqref{eq:overdamped}.
Therefore the second term and the third term in \eqref{eq:decomp}
can be bounded the same as in \cite{Raginsky} for the overdamped dynamics.

Specifically, the second term in \eqref{eq:decomp} can be bounded as
\begin{equation*}
\mathbb{E}F(\hat{X}^*) -  \mathbb{E}F_{\mathbf{Z}}(\hat{X}^*)
\leq
\frac{4\beta c_{LS}}{n}
\left(\frac{M^{2}}{m}(b+d/\beta)+B^{2}\right) = \mathcal{J}_3 (n),
\end{equation*}
by applying Lemma~\ref{lem:stability}, and the last term in \eqref{eq:decomp}
can be bounded as
\begin{eqnarray*}
\mathbb{E}F_{\mathbf{Z}}(\hat{X}^*) - F^* &=& \mathbb{E} \left[ F_{\mathbf{Z}} (\hat{X}^*) - \min_{x \in \mathbb{R}^d} F_{\mathbf{Z}} (x) \right] + \mathbb{E}\left[  \min_{x \in \mathbb{R}^d} F_{\mathbf{Z}} (x) - F_{\mathbf{Z}} (x^*) \right] \\
&\le & \mathbb{E} \left[ F_{\mathbf{Z}}(\hat{X}^*) - \min_{x \in \mathbb{R}^d} F_{\mathbf{Z}}(x) \right] \leq\mathcal{J}_2,
\end{eqnarray*}
where $x^*$ is any minimizer of $F(x)$, i.e., $F(x^*) = F^*$,
and the last step is due to Lemma~\ref{lem:thirdbound}.
The proof is complete.

\section{Proof of Theorem~\ref{thm-first-term-2} 
and Corollary~\ref{coro:main-2}}\label{sec:proof:2}

The proof of Theorem~\ref{thm-first-term-2} (Corollary~\ref{coro:main-2}) is similar to the proof of  Theorem~\ref{thm-first-term} (Corollary~\ref{coro:main}). There are two key new results that we need to establish: 
a uniform (in time) $L^{2}$ bound for the SGHMC2 iterates $(\hat{X}_{k},\hat{V}_{k})$, and the diffusion approximation that characterizes the 2-Wasserstein distance between the SGHMC2 iterates and the continuous-time underdampled Langevin diffusion. We summarize these two results in the following two lemmas and defer their proofs to Section~\ref{lemma-proof-2}. With these two lemmas, Theorem~\ref{thm-first-term-2} and Corollary~\ref{coro:main-2} readily follow and we omit the proof details.

\begin{lemma}[Uniform $L^{2}$ bounds for SGHMC2 iterates]\label{lem:L2bound:Jordan} ~\\
For 
$0<\eta\leq\min\left\{1,\frac{\gamma}{\hat{K}_{2}}(d/\beta+A/\beta),\frac{\gamma\lambda}{2\hat{K}_{1}},\frac{2}{\gamma\lambda}\right\}$,
where
\begin{align}
&\hat{K}_{1}:=K_{1}+Q_{1}\frac{4}{1-2\lambda}+Q_{2}\frac{8}{(1-2\lambda)\gamma^{2}},\label{def-K1-hat}
\\
&\hat{K}_{2}:=K_{2}+Q_{3},\label{def-K2-hat}
\end{align}
where $K_{1}$, $K_{2}$ are defined in \eqref{def-K1} and \eqref{def-K2},
and
\begin{align}
&Q_{1}:=\frac{1}{2}c_{0}
\Bigg((5M+4-2\gamma+(c_{0}+\gamma^{2}))
+(1+\gamma)\left(\frac{5}{2}+c_{0}(1+\gamma)\right)
+2\gamma^{2}\lambda\Bigg),\label{def-Q1}
\\
&Q_{2}:=\frac{1}{2}c_{0}
\Bigg[\Bigg((1+\gamma)\left(c_{0}(1+\gamma)+\frac{5}{2}\right)
+c_{0}+2
+\lambda\gamma^{2}
+2(Mc_{0}+M+1)\Bigg)\left(2(1+\delta)M^{2}\right)\nonumber
\\
&\qquad\qquad\qquad
+\left(2M^{2}+\gamma^{2}\lambda+\frac{3}{2}\gamma^{2}(1+\gamma)\right)\Bigg],\label{def-Q2}
\\
&Q_{3}:=c_{0}
\Bigg((1+\gamma)\left(c_{0}(1+\gamma)+\frac{5}{2}\right)
+c_{0}+2
+\lambda\gamma^{2}
+2(Mc_{0}+M+1)\Bigg)(1+\delta)B^{2}+c_{0}B^{2}\nonumber
\\
&\qquad\qquad\qquad
+\frac{1}{2}\gamma^{3}\beta^{-1}c_{22}
+\gamma^{2}\beta^{-1}c_{12}+M\gamma\beta^{-1}c_{22},\label{def-Q3}
\end{align}
where 
\begin{equation}\label{def-c0-c12-c22}
c_{0}:=1+\gamma^{2},\qquad c_{12}:=\frac{d}{2},\qquad c_{22}:=\frac{d}{3},
\end{equation}
we have
\begin{equation}
\sup_{j\geq 0}\mathbb{E}_{\mathbf{z}}\Vert\hat{X}_{j}\Vert^{2}
\leq
C_{x}^{d},
\qquad
\sup_{j\geq 0}\mathbb{E}_{\mathbf{z}}\Vert\hat{V}_{j}\Vert^{2}
\leq
C_{v}^{d},
\end{equation}
where $C_{x}^{d}$ and $C_{v}^{d}$ are defined in \eqref{def-Cx-d} and \eqref{def-Cxv-d}.
\end{lemma}


Next, let us provide a diffusion approximation between
the SGHMC2 algorithm $(\hat{X}_{k},\hat{V}_{k})$ and
the continuous time underdamped diffusion process $(X(k\eta),V(k\eta))$,
and we use $\hat{\mu}_{\mathbf{z},k}$ to denote
the law of $(\hat{X}_{k},\hat{V}_{k})$ and $\nu_{\mathbf{z},k}$
to denote the law of $(X(k\eta),V(k\eta))$.

\begin{lemma}[Diffusion approximation] \label{diff-app-2}
For any $k\in\mathbb{N}$ and any $\eta$, so that $k\eta\geq e$
and $\eta$ satisfies the condition in Lemma~\ref{lem:L2bound:Jordan}, we have
\begin{equation*}
\mathcal{W}_{2}(\hat{\mu}_{{\mathbf{z}},k},\nu_{{\mathbf{z}},k\eta})
\leq(C_{0}\delta^{1/4}+\hat{C}_{1}\eta^{1/2})\cdot\sqrt{k\eta}\cdot\sqrt{\log(k\eta)},
\end{equation*}
where $C_{0}$ is defined in \eqref{def-hat-Czero} and $\hat{C}_{1}$ is given by:
\begin{align}
&\hat{C}_{1}:=\hat{\gamma}
\cdot\Bigg(\frac{3\beta M^{2}}{2\gamma}
\bigg(C_{v}^{d}+\left(2(1+\delta)M^{2}C_{x}^{d}+2(1+\delta)B^{2}\right)+\frac{2d\gamma\beta^{-1}}{3}\bigg)
\nonumber
\\
&\qquad\qquad\qquad
+\sqrt{\frac{3\beta M^{2}}{2\gamma}
\bigg(C_{v}^{d}+\left(2(1+\delta)M^{2}C_{x}^{d}+2(1+\delta)B^{2}\right)+\frac{2d\gamma\beta^{-1}}{3}\bigg)}\Bigg)^{1/2},
\label{def-hat-Jone-2}
\end{align}
where $\hat{\gamma}$ is defined in \eqref{def-hat-gamma}.
\end{lemma}

\section{Proofs of Lemmas in Section~\ref{sec:proof}} \label{sec:lemma-proof-1}

\subsection{Proof of Lemma~\ref{lem:L2bound}}

(i) We first prove the continuous--time case. The main idea is to use the following Lyapunov function (see \eqref{eq:lyapunov})
introduced in \cite{Eberle} for the underdamped Langevin diffusion:
\begin{equation} \label{eq:V1}
\mathcal{V}(x,v)=\beta F_{\mathbf{z}}(x)
+\frac{\beta}{4}\gamma^{2}(\Vert x+\gamma^{-1}v\Vert^{2}+\Vert\gamma^{-1}v\Vert^{2}-\lambda\Vert x\Vert^{2}).
\end{equation}
Lemma~1.3 in \cite{Eberle} showed that if the drift condition in \eqref{eq:drift} holds, then
\begin{equation} \label{eq:V2}
\mathcal{L} \mathcal{V} \le \gamma (d+A - \lambda \mathcal{V}),
\end{equation}
where $\mathcal{L}$ is the infinitesimal generator of the underdamped Langevin diffusion $(X,V)$ defined in \eqref{eq:VL}--\eqref{eq:XL}:
\begin{equation} \label{eq:V3}
\mathcal{L}\mathcal{V}
=-(\gamma v+\nabla F_{\mathbf{z}}(x))\nabla_{v}\mathcal{V}
+\gamma\beta^{-1}\Delta_{v}\mathcal{V}
+v\nabla_{x}\mathcal{V}.
\end{equation}

To show part (i), we first note that for $\lambda \le \frac{1}{4},$
\begin{align} \label{eq:V-estimate}
\mathcal{V}(x, v) & \ge \beta F_{\mathbf{z}}(x)
+\frac{\beta}{4}(1-2 \lambda )\gamma^{2}(\Vert x+\gamma^{-1}v\Vert^{2}+\Vert\gamma^{-1}v\Vert^{2}) \nonumber\\
&\geq
\max\left\{\frac{1}{8} (1-2 \lambda) \beta \gamma^2 \Vert x\Vert^2,\frac{\beta}{4}(1-2\lambda)\Vert v\Vert^{2}\right\}.
\end{align}

Now let us set for each $t \ge 0,$
\begin{equation}\label{eq:Lt}
L(t)  := \mathbb{E}_{\mathbf{z}} [ \mathcal{V}(X(t),V(t)) ],
\end{equation}
and we will provide an upper bound for $L(t)$.

First, we can compute that
\begin{equation}\label{eq:vV}
\nabla_{v}\mathcal{V}=\beta v+\frac{\beta\gamma}{2}x,
\end{equation}
By It\^{o}'s formula and \eqref{eq:vV},
\begin{align*}
d(e^{\gamma\lambda t}\mathcal{V}(X(t),V(t)))
&=\gamma\lambda e^{\gamma\lambda t}\mathcal{V}(X(t),V(t))dt
+e^{\gamma\lambda t}\mathcal{L}\mathcal{V}(X(t),V(t))dt
\\
&\qquad
+e^{\gamma\lambda t}\left(\beta V(t)+\frac{\beta\gamma}{2}X(t)\right)\cdot\sqrt{2\gamma\beta^{-1}}dB(t),
\end{align*}
which together with \eqref{eq:V2} implies that
\begin{align}\label{eq:ito}
e^{\gamma\lambda t}\mathcal{V}(X(t),V(t))
&\leq\mathcal{V}(X(0),V(0))+\gamma(d+A)\int_{0}^{t}e^{\lambda\gamma s}ds
\nonumber\\
&\qquad
-\int_{0}^{t}e^{\gamma\lambda s}\left(\beta V(s)+\frac{\beta\gamma}{2}X(s)\right)\cdot\sqrt{2\gamma\beta^{-1}}d B(s).
\end{align}
Note that $\nabla F_{\bz}(x)$ is Lipschitz continuous by part (ii) of Assumption~\ref{assumptions}, and hence
$(X(t),V(t))$ is the unique strong solution of the SDE \eqref{eq:VL}-\eqref{eq:XL},
and thus $\mathbb{E}[\int_{0}^{T}\Vert V(t)\Vert^{2}+\Vert X(t)\Vert^{2}dt]<\infty$
for every $T>0$ (See e.g. \cite{Oksendal}).
Therefore, for every $T>0$, we have
\begin{equation*}
\int_{0}^{T}e^{2\gamma\lambda s}\left\Vert\beta V(s)+\frac{\beta\gamma}{2}X(s)\right\Vert^{2}(2\gamma\beta^{-1})ds<\infty,
\end{equation*}
and hence $\int_{0}^{t}e^{\gamma\lambda s}\left( \beta V(s)+\frac{\beta\gamma}{2}X(s)\right)\cdot\sqrt{2\gamma\beta^{-1}}B(s)$
is a martingale. Then we can infer from \eqref{eq:ito} and \eqref{eq:Lt} that for any $t \ge 0$,
\begin{equation*}
L(t)=\mathbb{E}_{\mathbf{z}} [ \mathcal{V}(X(t),V(t)) ]  \le L(0) e^{-\gamma \lambda t} + \frac{d+A}{\lambda} (1- e^{-\gamma \lambda t}).
\end{equation*}
In combination with \eqref{eq:V-estimate}, we obtain that $(X,V)$ are uniformly (in time) $L^2$ bounded.
Indeed, we have
\begin{align*}
&\frac{1}{8} (1-2 \lambda) \beta \gamma^2 \mathbb{E}_{\mathbf{z}}\Vert X(t)\Vert^2
\leq\mathbb{E}_{\mathbf{z}}[\mathcal{V}(X_{0},V_{0})]+\frac{d+A}{\lambda},
\\
&\frac{\beta}{4}(1-2\lambda)\mathbb{E}_{\mathbf{z}}\Vert V(t)\Vert^2
\leq\mathbb{E}_{\mathbf{z}}[\mathcal{V}(X_{0},V_{0})]+\frac{d+A}{\lambda}.
\end{align*}
The proof of part (i) is complete by noting that $\mathbb{E}_{\mathbf{z}}[\mathcal{V}(X_{0},V_{0})]$ is finite from part (v) of Assumption~\ref{assumptions}.

(ii) Next, we prove the uniform (in time) $L^2$ bounds for $(X_k, V_k)$.
Let us recall the dynamics:
\begin{align}
&V_{k+1}=V_{k}-\eta[\gamma V_{k}+g(X_{k},U_{\bz,k})]+\sqrt{2\gamma\beta^{-1}\eta}\xi_{k}, \label{eq:VL1}
\\
&X_{k+1}=X_{k}+\eta V_{k}, \label{eq:XL1}
\end{align}
where $\mathbb{E}g(x,U_{{\mathbf{z}}, k}) = \nabla F_{\mathbf{z}} (x)$ for any $x$. We again use the Lyapunov function $\mathcal{V}(x,v)$
in \eqref{eq:V1}, and set for each $k=0,1, \ldots,$
\begin{equation} \label{eq:L2k}
L_2(k) = \mathbb{E}_{\mathbf{z}} \mathcal{V}(X_k, V_k)/\beta= \mathbb{E}_{\mathbf{z}}  \left[F_{\mathbf{z}}(X_k)
+\frac{1}{4}\gamma^{2}\left(\Vert X_k +\gamma^{-1}V_k \Vert^{2}+\Vert\gamma^{-1}V_k \Vert^{2}-\lambda\Vert X_k\Vert^{2}\right) \right].
\end{equation}
We show below that one can find explicit constants $K_1, K_2 >0$, such that
\begin{equation*}
(L_2(k+1) - L_2(k) )/\eta \le \gamma (d/\beta+A/\beta- \lambda L_2(k)) + (K_{1}L_{2}(k)+K_{2})\cdot \eta.
\end{equation*}

We proceed in several steps in upper bounding $L_2(k+1)$.

First, by using the independence of $V_{k}-\eta[\gamma V_{k}+g_{k}(X_{k},U_{z,k})]$ and $\xi_{k}$,
we can obtain from \eqref{eq:VL1} that
\begin{align*}
&\mathbb{E}_{\mathbf{z}}\Vert V_{k+1}\Vert^{2}
\\
&=\mathbb{E}_{\mathbf{z}}\Vert V_{k}-\eta[\gamma V_{k}+g_{k}(X_{k},U_{\bz,k})]\Vert^{2}+2\gamma\beta^{-1}\eta\mathbb{E}_{\mathbf{z}}\Vert\xi_{k}\Vert^{2}
\\
&=\mathbb{E}_{\mathbf{z}}\Vert V_{k}-\eta[\gamma V_{k}+g_{k}(X_{k},U_{\bz,k})]\Vert^{2}+2\gamma\beta^{-1}\eta d
\\
&=\mathbb{E}_{\mathbf{z}}\Vert V_{k}-\eta[\gamma V_{k}+\nabla F_{\mathbf{z}}(X_{k})]\Vert^{2}+2\gamma\beta^{-1}\eta d
+\eta^{2}\mathbb{E}_{\mathbf{z}}\Vert\nabla F_{\mathbf{z}}(X_{k})-g_{k}(X_{k},U_{\bz,k})\Vert^{2}
\\
&\leq
(1-\eta\gamma)^{2}\mathbb{E}_{\mathbf{z}}\Vert V_{k}\Vert^{2}
-2\eta(1-\eta\gamma)\mathbb{E}_{\mathbf{z}}[\langle V_{k},\nabla F_{\mathbf{z}}(X_{k})\rangle]
+\eta^{2}\mathbb{E}_{\mathbf{z}}\Vert\nabla F_{\mathbf{z}}(X_{k})\Vert^{2}+2\gamma\beta^{-1}\eta d
\\
&\qquad\qquad
+2\delta\eta^{2}M^{2}\mathbb{E}_{\mathbf{z}}\Vert X_{k}\Vert^{2}
+2\delta\eta^{2}B^{2}
\\
&\leq
(1-\eta\gamma)^{2}\mathbb{E}_{\mathbf{z}}\Vert V_{k}\Vert^{2}
-2\eta (1-\eta\gamma)\mathbb{E}_{\mathbf{z}}[\langle V_{k},\nabla F_{\mathbf{z}}(X_{k})\rangle]
\\
&\qquad
+\eta^{2}(M^{2}\mathbb{E}_{\mathbf{z}}\Vert X_{k}\Vert^{2}
+B^{2}+2MB\mathbb{E}_{\mathbf{z}}\Vert X_{k}\Vert)+2\gamma\beta^{-1}\eta d
\\
&\qquad\qquad
+2\delta\eta^{2}M^{2}\mathbb{E}_{\mathbf{z}}\Vert X_{k}\Vert^{2}
+2\delta\eta^{2}B^{2},
\end{align*}
where we have used part (iv) of Assumption~\ref{assumptions} and Lemma~\ref{lem:gradient-bound} in Section~\ref{sec:support} in the Appendix.
By using $|x| \le \frac{|x|^2+1}{2}$, we immediately get
\begin{align}
\mathbb{E}_{\mathbf{z}}\Vert V_{k+1}\Vert^{2}
\nonumber
&\leq
(1-\eta\gamma)^{2}\mathbb{E}_{\mathbf{z}}\Vert V_{k}\Vert^{2}
-2\eta \mathbb{E}_{\mathbf{z}}[\langle V_{k},\nabla F_{\mathbf{z}}(X_{k})\rangle] + 2 \eta^2 \gamma \mathbb{E}_{\mathbf{z}}[\langle V_{k},\nabla F_{\mathbf{z}}(X_{k})\rangle]
\\
&\qquad
+ \left(\eta^{2}M^{2} + \eta^2 MB + \delta\eta^{2}M^{2} \right) \mathbb{E}_{\mathbf{z}}\Vert X_{k}\Vert^{2} +  (\eta^2 MB + 2\gamma\beta^{-1}\eta d + 2\delta\eta^{2}B^{2}).\label{eqn:first}
\end{align}

Second, we can compute from \eqref{eq:XL1} that
\begin{equation}\label{eqn:second}
\mathbb{E}_{\mathbf{z}}\Vert X_{k+1}\Vert^{2}
=\mathbb{E}_{\mathbf{z}}\Vert X_{k}\Vert^{2}
+2\eta\mathbb{E}_{\mathbf{z}}[\langle X_{k},V_{k}\rangle]
+\eta^{2}\mathbb{E}_{\mathbf{z}}\Vert V_{k}\Vert^{2}.
\end{equation}

Third, note that
\begin{align*}
F_{\mathbf{z}}(X_{k+1})&= F_{\mathbf{z}}(X_{k} + \eta V_k)
=F_{\mathbf{z}}(X_{k}) + \int_{0}^1 \langle \nabla F_{\mathbf{z}}(X_{k} + \tau \eta V_k), \eta V_k \rangle d \tau,
\end{align*}
which immediately suggests that
\begin{align*}
\left|F_{\mathbf{z}}(X_{k+1})-  F_{\mathbf{z}}(X_{k} ) -  \left\langle \nabla F_{\mathbf{z}}(X_{k} ), \eta V_k  \right\rangle \right|
&=\left|\int_{0}^1 \langle \nabla F_{\mathbf{z}}(X_{k} + \tau \eta V_k) - \nabla F_{\mathbf{z}}(X_{k} ), \eta V_k \rangle d \tau \right|\\
&\le \int_{0}^1 \left\Vert \nabla F_{\mathbf{z}}(X_{k} + \tau \eta V_k) - \nabla F_{\mathbf{z}}(X_{k} ) \right\Vert \cdot \left\Vert\eta V_k \right\Vert d \tau \\
& \le \frac{1}{2} M \eta^2 \Vert V_k \Vert^2,
\end{align*}
where the last inequality is due to the $M-$smoothness of $F_{\mathbf{z}}$. This implies
\begin{equation}\label{eqn:third}
 \mathbb{E}_{\mathbf{z}} F_{\mathbf{z}}(X_{k+1})-   \mathbb{E}_{\mathbf{z}} F_{\mathbf{z}}(X_{k} ) \le  \eta  \mathbb{E}_{\mathbf{z}} \langle \nabla F_{\mathbf{z}}(X_{k} ), V_k  \rangle
+  \frac{1}{2} M \eta^2 \mathbb{E}_{\mathbf{z}} \Vert V_k \Vert^2.
\end{equation}

Finally, we can compute that
\begin{align}
&\mathbb{E}_{\mathbf{z}}\left\Vert X_{k+1} +\gamma^{-1}V_{k+1} \right\Vert^{2}
\nonumber
\\
&= \mathbb{E}_{\mathbf{z}} \left\Vert X_{k} +\gamma^{-1}V_{k} - \eta \gamma^{-1} g(X_{k},U_{\bz,k})+\sqrt{2\gamma^{-1}\beta^{-1}\eta}\xi_{k}\right\Vert^{2}
\nonumber
\\
&= \mathbb{E}_{\mathbf{z}} \left\Vert X_{k} +\gamma^{-1}V_{k} - \eta \gamma^{-1} g(X_{k},U_{\bz,k})\right\Vert^{2}
+2\gamma^{-1}\beta^{-1}\eta d
\nonumber
\\
&= \mathbb{E}_{\mathbf{z}} \left\Vert X_{k} +\gamma^{-1}V_{k} - \eta \gamma^{-1} \nabla F_{\mathbf{z}}(X_{k})\right\Vert^{2}
+2\gamma^{-1}\beta^{-1}\eta d
\nonumber
\\
&\qquad\qquad\qquad\qquad
+\mathbb{E}_{\mathbf{z}}\left\Vert\eta\gamma^{-1}g(X_{k},U_{\bz,k})
-\eta\gamma^{-1}\nabla F_{\mathbf{z}}(X_{k})\right\Vert^{2}
\nonumber
\\
&\leq
\mathbb{E}_{\mathbf{z}} \left\Vert X_{k} +\gamma^{-1}V_{k} - \eta \gamma^{-1} \nabla F_{\mathbf{z}}(X_{k})\right\Vert^{2}
+2\gamma^{-1}\beta^{-1}\eta d
+2\eta^{2}\gamma^{-2}\delta(M^{2}\mathbb{E}_{\mathbf{z}}\left\Vert X_{k}\right\Vert^{2}+B^{2})
\nonumber
\\
& = \mathbb{E}_{\mathbf{z}}\left\Vert X_{k} +\gamma^{-1}V_{k} \right\Vert^{2}
- 2 \eta \gamma^{-1} \mathbb{E}_{\mathbf{z}}\langle X_{k} +\gamma^{-1} V_{k},\nabla F_{\mathbf{z}}(X_{k})\rangle
\nonumber
\\
& \qquad +\eta^{2} \gamma^{-2}   \mathbb{E}_{\mathbf{z}}\left\Vert \nabla F_{\mathbf{z}}(X_{k}) \right\Vert^2
+2\gamma^{-1}\beta^{-1}\eta d
+2\eta^{2}\gamma^{-2}\delta(M^{2}\mathbb{E}_{\mathbf{z}}\Vert X_{k}\Vert^{2}+B^{2}),\label{eqn:fourth}
\end{align}
where we have used part (iv) of Assumption~\ref{assumptions} in the inequality above.

Combining the equations \eqref{eqn:first}, \eqref{eqn:second}, \eqref{eqn:third} and \eqref{eqn:fourth},
we get
\begin{align}
&(L_2(k+1) - L_2(k) )/\eta
\nonumber\\
&=\bigg(\mathbb{E}_{\mathbf{z}}[F_{\mathbf{z}}(X_{k+1})]-\mathbb{E}_{\mathbf{z}}[F_{\mathbf{z}}(X_{k})]
+\frac{1}{4}\gamma^{2}\left(\mathbb{E}_{\mathbf{z}}\Vert X_{k+1}+\gamma^{-1}V_{k+1}\Vert^{2}
-\mathbb{E}_{\mathbf{z}}\Vert X_{k}+\gamma^{-1}V_{k}\Vert^{2}\right)
\nonumber\\
&\qquad
+\frac{1}{4}\left(\mathbb{E}_{\mathbf{z}}\Vert V_{k+1}\Vert^{2}-\mathbb{E}_{\mathbf{z}}\Vert V_{k}\Vert^{2}\right)
-\frac{1}{4}\gamma^{2}\lambda\left(\mathbb{E}_{\mathbf{z}}\Vert X_{k+1}\Vert^{2}
-\mathbb{E}_{\mathbf{z}}\Vert X_{k}\Vert^{2}\right)\bigg)\bigg/\eta
 \nonumber \\
&\leq
\mathbb{E}_{\mathbf{z}} \langle \nabla F_{\mathbf{z}}(X_{k} ), V_k  \rangle
+\frac{1}{2} M \eta\mathbb{E}_{\mathbf{z}} \Vert V_k \Vert^2
- \frac{1}{2}\gamma\mathbb{E} \left\langle X_{k} +\gamma^{-1} V_{k},\nabla F_{\mathbf{z}}(X_{k})\right\rangle
\nonumber
\\
& \qquad
+\frac{1}{4}\eta \mathbb{E} \Vert \nabla F_{\mathbf{z}}(X_{k}) \Vert^2
+\frac{1}{2}\gamma\beta^{-1}d
+\frac{1}{2}\eta\delta(M^{2}\mathbb{E}\Vert X_{k}\Vert^{2}+B^{2})
\nonumber \\
&\qquad\qquad
+\frac{1}{4}(-2\gamma+\eta\gamma^{2})\mathbb{E}_{\mathbf{z}}\Vert V_{k}\Vert^{2}
-\frac{1}{2}(1-\eta\gamma)\mathbb{E}_{\mathbf{z}}[\langle V_{k},\nabla F_{\mathbf{z}}(X_{k})\rangle]
\nonumber \\
&\qquad
+\frac{1}{4}\eta(M^{2}\mathbb{E}_{\mathbf{z}}\Vert X_{k}\Vert^{2}
+B^{2}+2MB\mathbb{E}_{\mathbf{z}}\Vert X_{k}\Vert)
+\frac{1}{2}\gamma\beta^{-1}d
\nonumber\\
&\qquad\qquad
+\frac{1}{2}\delta\eta M^{2}\mathbb{E}_{\mathbf{z}}\Vert X_{k}\Vert^{2}
+\frac{1}{2}\delta\eta B^{2}
-\frac{1}{2}\gamma^{2}\lambda\mathbb{E}_{\mathbf{z}}\langle X_{k},V_{k}\rangle
-\frac{1}{4}\gamma^{2}\lambda\eta\mathbb{E}_{\mathbf{z}}\Vert V_{k}\Vert^{2}
\nonumber \\
&=  - \frac{\gamma}{2} \mathbb{E}_{\mathbf{z}} \langle \nabla F_{\mathbf{z}}(X_{k} ), X_k  \rangle  - \frac{\gamma}{2} \mathbb{E}_{\mathbf{z}} \Vert V_k \Vert^2
- \frac{\gamma^2 \lambda}{2} \mathbb{E}_{\mathbf{z}} \langle X_{k}, V_k  \rangle + \gamma\beta^{-1}d+\mathcal{E}_{k}\eta
\nonumber \\
&\leq
-\gamma\lambda\mathbb{E}_{\mathbf{z}}[F_{\mathbf{z}}(X_{k})]
-\frac{1}{4}\lambda\gamma^{3}\mathbb{E}_{\mathbf{z}}\Vert X_{k}\Vert^{2}
+\gamma A/\beta
- \frac{\gamma}{2} \mathbb{E}_{\mathbf{z}} \Vert V_k \Vert^2
- \frac{\gamma^2 \lambda}{2} \mathbb{E}_{\mathbf{z}} \langle X_{k}, V_k  \rangle + \gamma\beta^{-1}d+\mathcal{E}_{k}\eta, \label{eq:disc-bound}
\end{align}
where we used the drift condition \eqref{eq:drift} in the last inequality, and
\begin{align*}
\mathcal{E}_{k}&:=\left(\frac{1}{2}M+\frac{1}{4}\gamma^{2}-\frac{1}{4}\gamma^{2}\lambda\right)\mathbb{E}_{\mathbf{z}}\Vert V_{k}\Vert^{2}
+\frac{1}{4}\mathbb{E}_{\mathbf{z}}\Vert\nabla F_{\mathbf{z}}(X_{k})\Vert^{2}+\delta(M^{2}\mathbb{E}\Vert X_{k}\Vert^{2}+B^{2})
\\
&\qquad
+\frac{1}{2}\gamma\mathbb{E}_{\mathbf{z}}[\langle V_{k},\nabla F_{\mathbf{z}}(X_{k})\rangle]
+\frac{1}{4}(M^{2}\mathbb{E}_{\mathbf{z}}\Vert X_{k}\Vert^{2}
+B^{2}+2MB\mathbb{E}_{\mathbf{z}}\Vert X_{k}\Vert).
\end{align*}
We can upper bound $\mathcal{E}_{k}$ as follows:
\begin{align*}
\mathcal{E}_{k}&
\leq
\left(\frac{1}{2}M+\frac{1}{4}\gamma^{2}-\frac{1}{4}\gamma^{2}\lambda\right)\mathbb{E}_{\mathbf{z}}\Vert V_{k}\Vert^{2}
+\frac{1}{4}\mathbb{E}_{\mathbf{z}}\Vert\nabla F_{\mathbf{z}}(X_{k})\Vert^{2}+\delta(M^{2}\mathbb{E}_{\mathbf{z}}\Vert X_{k}\Vert^{2}+B^{2})
\\
&\qquad
+\gamma\mathbb{E}_{\mathbf{z}}\Vert V_{k}\Vert^{2}
+\gamma\mathbb{E}_{\mathbf{z}}\Vert\nabla F_{\mathbf{z}}(X_{k})\Vert^{2}
+\frac{1}{4}\mathbb{E}_{\mathbf{z}}(M\Vert X_{k}\Vert+B)^{2}
\\
&\leq
\left(\frac{1}{2}M+\frac{1}{4}\gamma^{2}-\frac{1}{4}\gamma^{2}\lambda+\gamma\right)\mathbb{E}_{\mathbf{z}}\Vert V_{k}\Vert^{2}
+\delta(M^{2}\mathbb{E}_{\mathbf{z}}\Vert X_{k}\Vert^{2}+B^{2})
\\
&\qquad
+\left(\frac{1}{4}+\gamma\right)\mathbb{E}_{\mathbf{z}}(M\Vert X_{k}\Vert+B)^{2}
+\frac{1}{4}\mathbb{E}_{\mathbf{z}}(M\Vert X_{k}\Vert+B)^{2}
\\
&\leq
\left(\frac{1}{2}M+\frac{1}{4}\gamma^{2}-\frac{1}{4}\gamma^{2}\lambda+\gamma\right)\mathbb{E}_{\mathbf{z}}\Vert V_{k}\Vert^{2}
\\
&\qquad\qquad\qquad\qquad
+2M^{2}\left(\frac{1}{2}+\gamma+\delta\right)\mathbb{E}_{\mathbf{z}}\Vert X_{k}\Vert^{2}
+2B^{2}\left(\frac{1}{2}+\gamma+\delta\right).
\end{align*}
Since $\lambda\leq\frac{1}{4}$, we obtain from \eqref{eq:V-estimate} and \eqref{eq:L2k} that
\begin{align}
L_{2}(k)
&\geq\max\left\{\frac{1}{8} (1-2 \lambda) \gamma^2 \mathbb{E}_{\mathbf{z}}\Vert X_{k}\Vert^2,\frac{1}{4}(1-2\lambda)\mathbb{E}_{\mathbf{z}}\Vert V_{k}\Vert^{2}\right\}
\label{eq:bound-L2k} \\
&\geq
\frac{1}{16} (1-2 \lambda)\gamma^2 \mathbb{E}_{\mathbf{z}}\Vert X_{k}\Vert^2
+\frac{1}{8}(1-2\lambda)\mathbb{E}_{\mathbf{z}}\Vert V_{k}\Vert^{2}. \nonumber
\end{align}
Therefore,
\begin{equation} \label{eq:Ek-bound}
\mathcal{E}_{k}\leq K_{1}L_{2}(k)+K_{2},
\end{equation}
where we recall from \eqref{def-K1} and \eqref{def-K2} that
\begin{equation*}
K_{1}
=\max\left\{\frac{2M^{2}\left(\frac{1}{2}+\gamma+\delta\right)}{\frac{1}{16} (1-2 \lambda)\gamma^2 },
\frac{\left(\frac{1}{2}M+\frac{1}{4}\gamma^{2}-\frac{1}{4}\gamma^{2}\lambda+\gamma\right)}{\frac{1}{8}(1-2\lambda)}\right\},
\end{equation*}
and
\begin{equation*}
K_{2}=2B^{2}\left(\frac{1}{2}+\gamma+\delta\right).
\end{equation*}
Moreover, since $\lambda\leq\frac{1}{4}$, we infer from the definition of $L_2(k)$ in \eqref{eq:L2k} that
\begin{align}
L_{2}(k)
&=\mathbb{E}_{\mathbf{z}}[F_{\mathbf{z}}(X_{k})]+\frac{1}{4}\gamma^{2}(1-\lambda)
\mathbb{E}_{\mathbf{z}}\Vert X_{k}\Vert^{2}
+\frac{1}{2}\gamma\mathbb{E}_{\mathbf{z}}[\langle X_{k},V_{k}\rangle]
+\frac{1}{2}\mathbb{E}_{\mathbf{z}}\Vert V_{k}\Vert^{2}
\nonumber \\
&\leq
\mathbb{E}_{\mathbf{z}}[F_{\mathbf{z}}(X_{k})]+\frac{1}{4}\gamma^{2}\mathbb{E}_{\mathbf{z}}\Vert X_{k}\Vert^{2}
+\frac{1}{2}\gamma\mathbb{E}_{\mathbf{z}}[\langle X_{k},V_{k}\rangle]
+\frac{1}{2\lambda}\mathbb{E}_{\mathbf{z}}\Vert V_{k}\Vert^{2}. \nonumber
\end{align}
Together with \eqref{eq:disc-bound} and \eqref{eq:Ek-bound}, we deduce that
\begin{equation*}
(L_{2}(k+1)-L_{2}(k))/\eta
\leq\gamma(d/\beta+A/\beta-\lambda L_{2}(k))+(K_{1}L_{2}(k)+K_{2})\eta.
\end{equation*}
For $0 < \eta\leq\min\left\{\frac{\gamma}{K_{2}}(d/\beta+A/\beta),\frac{\gamma\lambda}{2K_{1}}\right\}$,
we get
\begin{equation*}
(L_{2}(k+1)-L_{2}(k))/\eta
\leq
2\gamma(d/\beta+A/\beta)-\frac{1}{2}\gamma\lambda L_{2}(k),
\end{equation*}
which implies
\begin{equation*}
L_2(k+1) \le \rho L_2(k) +K ,
\end{equation*}
where
\begin{equation*}
\rho:= 1- \eta \gamma \lambda/2,
\qquad
K:=2\eta\gamma(d/\beta+A/\beta),
\end{equation*}
and we have $\rho\in[0,1)$, where we used the assumption that $\eta\leq\frac{2}{\gamma\lambda}$.
It follows that
\begin{align*}
L_{2}(k)\leq L_{2}(0)+\frac{K}{1-\rho}
=\mathbb{E}_{\mathbf{z}}\left[\mathcal{V}(X_{0},V_{0})/\beta\right]
+\frac{4(d/\beta+A/\beta)}{\lambda}.
\end{align*}
The result then follows from the inequality above and \eqref{eq:bound-L2k}.

\subsection{Proof of Lemma~\ref{lem:exp-integrable}}
{From} \eqref{eq:V1}--\eqref{eq:V3}, we can directly obtain that
\begin{align}
\mathcal{L}e^{\alpha\mathcal{V}}
&=\left[-(\gamma v+\nabla F_{\mathbf{z}}(x))\alpha\nabla_{v}\mathcal{V}
+\gamma\beta^{-1}\alpha\Delta_{v}\mathcal{V}
+\gamma\beta^{-1}\alpha^{2}\Vert\nabla_{v}\mathcal{V}\Vert^{2}
+v\alpha\nabla_{x}\mathcal{V}\right]e^{\alpha\mathcal{V}}
\nonumber
\\
&=\left[\alpha\mathcal{L}\mathcal{V}
+\gamma\beta^{-1}\alpha^{2}\Vert\nabla_{v}\mathcal{V}\Vert^{2}\right]e^{\alpha\mathcal{V}}
\nonumber
\\
&\leq\left[\alpha\gamma d+\alpha\gamma A-\alpha\gamma\lambda\mathcal{V}+\alpha^{2}\gamma\beta^{-1}\Vert\nabla_{v}\mathcal{V}\Vert^{2}\right]e^{\alpha\mathcal{V}}. \label{eq:generator-ineq}
\end{align}
Moreover, we recall from \eqref{eq:vV} that
\begin{equation*}
\nabla_{v}\mathcal{V}=\beta v+\frac{\beta\gamma}{2}x,
\end{equation*}
and thus
\begin{equation*}
\Vert\nabla_{v}\mathcal{V}\Vert^{2}
\leq 2\beta^{2}\Vert v\Vert^{2}+\frac{\beta^{2}\gamma^{2}}{2}\Vert x\Vert^{2}.
\end{equation*}
We recall from \eqref{eq:V-estimate} that
\begin{equation*}
\mathcal{V}(x,v)
\geq\max\left\{\frac{1}{8} (1-2 \lambda) \beta \gamma^2 \Vert x\Vert^2,\frac{\beta}{4}(1-2\lambda)\Vert v\Vert^{2}\right\}.
\end{equation*}
Therefore, we have
\begin{equation}
\Vert\nabla_{v}\mathcal{V}\Vert^{2}
\leq\left[\frac{8\beta^{2}}{\beta(1-2\lambda)}
+\frac{4\beta^{2}\gamma^{2}}{(1-2\lambda)\beta\gamma^{2}}\right]\mathcal{V}
=\frac{12\beta}{1-2\lambda}\mathcal{V}. \label{eq:bound-gradV}
\end{equation}
By choosing:
\begin{equation} \label{eq:alpha}
\alpha=\frac{\lambda\beta}{\frac{12\beta}{1-2\lambda}}
=\frac{\lambda(1-2\lambda)}{12},
\end{equation}
we get
\begin{equation}\label{eq:keyineq}
\mathcal{L}e^{\alpha\mathcal{V}}
\leq\alpha\gamma(d+A)e^{\alpha\mathcal{V}}.
\end{equation}
Since $\mathcal{L}e^{\alpha\mathcal{V}}=\left[\mathcal{L}\alpha\mathcal{V}+\gamma\beta^{-1}\Vert\nabla_{v}\alpha\mathcal{V}\Vert^{2}\right]e^{\alpha\mathcal{V}}$, we have showed that
\begin{equation*}
\mathcal{L}\alpha\mathcal{V}+\gamma\beta^{-1}\Vert\nabla_{v}\alpha\mathcal{V}\Vert^{2}
\leq\alpha\gamma(d+A).
\end{equation*}
Applying an exponential integrability result, e.g. Corollary 2.4. in \cite{Cox2013}, we get
\begin{equation*}
\mathbb{E}\left[e^{\alpha\mathcal{V}(X(t),V(t))}\right]
\leq \mathbb{E}\left[e^{\alpha\mathcal{V}(X(0),V(0))}\right]e^{\alpha\gamma(d+A)t}.
\end{equation*}
That is,
\begin{equation}\label{eq:finiteness}
\mathbb{E}_{\mathbf{z}}\left[e^{\alpha\mathcal{V}(X(t),V(t))}\right]
\leq\int_{\mathbb{R}^{2d}}e^{\alpha\mathcal{V}(x,v)+\alpha\gamma(d+A)t}\mu_{0}(dx,dv)<\infty.
\end{equation}

Next, applying It\^{o}'s formula to $e^{\frac{1}{4}\alpha\mathcal{V}(X(t),V(t))}$,
we obtain
\begin{align}
e^{\frac{1}{4}\alpha\mathcal{V}(X(t),V(t))}
&=e^{\frac{1}{4}\alpha\mathcal{V}(X(0),V(0))}
+\int_{0}^{t}\mathcal{L}e^{\frac{1}{4}\alpha\mathcal{V}(X(s),V(s))}ds
\nonumber
\\
&\qquad\qquad\qquad
+\int_{0}^{t}\frac{1}{2}\left(\beta V(s)+\frac{\beta\gamma}{2}X(s)\right)
e^{\frac{1}{4}\alpha\mathcal{V}(X(s),V(s))}\cdot dB(s).
\label{eq:takeE}
\end{align}
For every $T>0$,
\begin{align*}
&\int_{0}^{T}\mathbb{E}\left\Vert\frac{1}{2}\left(\beta V(s)+\frac{\beta\gamma}{2}X(s)\right)
e^{\frac{1}{4}\alpha\mathcal{V}(X(s),V(s))}\right\Vert^{2}ds
\\
&\leq
\frac{\beta^{2}}{2}
\int_{0}^{T}
\mathbb{E}\left[\left(\Vert V(s)\Vert^{2}+\gamma^{2}\Vert X(s)\Vert^{2}\right)
e^{\frac{1}{2}\alpha\mathcal{V}(X(s),V(s))}\right]ds
\\
&\leq
\frac{6\beta}{1-2\lambda}
\int_{0}^{T}
\mathbb{E}\left[\mathcal{V}(X(s),V(s))
e^{\frac{1}{2}\alpha\mathcal{V}(X(s),V(s))}\right]ds
\\
&\leq
\frac{12\beta}{1-2\lambda}
\int_{0}^{T}
\mathbb{E}\left[e^{\alpha\mathcal{V}(X(s),V(s))}\right]ds<\infty,
\end{align*}
where we used \eqref{eq:V-estimate} and \eqref{eq:finiteness}.
Thus, $\int_{0}^{t}\frac{1}{2}\left(\beta V(s)+\frac{\beta\gamma}{2}X(s)\right)
e^{\frac{1}{4}\alpha\mathcal{V}(X(s),V(s))}\cdot dB(s)$
is a martingale. By taking expectations on both hand sides
of \eqref{eq:takeE}, we get
\begin{equation}\label{eq:follow}
\mathbb{E}\left[e^{\frac{1}{4}\alpha\mathcal{V}(X(s),V(s))}\right]
=\mathbb{E}\left[e^{\frac{1}{4}\alpha\mathcal{V}(X(0),V(0))}\right]
+\int_{0}^{t}\mathbb{E}\left[\mathcal{L}e^{\frac{1}{4}\alpha\mathcal{V}(X(s),V(s))}\right]ds.
\end{equation}
From \eqref{eq:generator-ineq}, \eqref{eq:bound-gradV} and \eqref{eq:alpha}, we can infer that
\begin{align*}
\mathcal{L}e^{\frac{1}{4}\alpha\mathcal{V}}
&\leq\left(\frac{1}{4}\alpha\gamma(d+A)
-\frac{1}{4}\alpha\gamma\lambda\mathcal{V}
+\gamma\beta^{-1}\frac{\alpha^{2}}{16}\Vert\nabla_{v}\mathcal{V}\Vert^{2}
\right)e^{\frac{1}{4}\alpha\mathcal{V}}
\\
&\leq\left(\frac{1}{4}\alpha\gamma(d+A)
-\frac{3}{16}\alpha\gamma\lambda\mathcal{V}\right)e^{\frac{1}{4}\alpha\mathcal{V}}
\\
&\leq\frac{1}{4}\alpha\gamma(d+A)e^{\frac{\alpha(d+A)}{3\lambda}}\,,
\end{align*}
where in the last inequality we used the facts that $\mathcal{V}\geq 0$ and $\frac{1}{4}\alpha\gamma(d+A)
-\frac{3}{16}\alpha\gamma\lambda\mathcal{V} \ge 0$ if and only if
$\mathcal{V} \le \frac{4(d+A)}{3\lambda}$. Therefore, it follows from \eqref{eq:follow} that
\begin{equation*}
\mathbb{E}\left[e^{\frac{1}{4}\alpha\mathcal{V}(X(s),V(s))}\right]
\leq\mathbb{E}\left[e^{\frac{1}{4}\alpha\mathcal{V}(X(0),V(0))}\right]
+\frac{1}{4}e^{\frac{\alpha(d+A)}{3\lambda}}\alpha\gamma(d+A)t.
\end{equation*}

Finally, by \eqref{eq:V-estimate} again,
\begin{equation*}
\Vert(x,v)\Vert^{2}
\leq
2\Vert x\Vert^{2}+2\Vert v\Vert^{2}
\leq
\left[\frac{16}{(1-2\lambda)\beta\gamma^{2}}+\frac{8}{\beta(1-2\lambda)}\right]\mathcal{V}(x,v).
\end{equation*}
Hence, the conclusion follows.

\subsection{Proof of Lemma~\ref{diff-app}}

The proof is inspired by the proof of Lemma 7 in \cite{Raginsky}
although more delicate in our setting.
Note that the main technical difficulty here is that
the underdamped Langevin diffusion is a hypoelliptic diffusion,
i.e. the diffusion matrix of the stochastic
differential equation defining the multidimensional diffusion process is not of full
rank, but its solutions admit a smooth density, see \cite{DS2019}.
In our case, there is no Brownian noise in $dX(t)$ term in \eqref{eq:XL}
and the underdamped Langevin diffusion \eqref{eq:VL}-\eqref{eq:XL} is hypoelliptic.
Consider the following continuous-time interpolation of $(X_k, V_k)$:
\begin{align}
&\overline{V}(t)=V_{0}-\int_{0}^{t}\gamma\overline{V}(\lfloor s/\eta\rfloor\eta)ds
-\int_{0}^{t}g(\overline{X}(\lfloor s/\eta\rfloor\eta),\overline{U}_{\mathbf{z}}(s))ds+\sqrt{2\gamma \beta^{-1}}\int_{0}^{t}dB(s),\label{Euler:continuous:V}
\\
&\overline{X}(t)=X_{0}+\int_{0}^{t}\overline{V}(\lfloor s/\eta\rfloor\eta)ds,\label{Euler:continuous:X}
\end{align}
where $\overline{U}_{\mathbf{z}}(t):=U_{\bz,k}$ for $k\eta\leq t<(k+1)\eta$. 
Then $(\overline{X}(k \eta), \overline{V}(k \eta))$ and $(X_k, V_k)$
have the same distribution $\mu_{\bz, k}$ for each $k \ge 0$. 
Since there is no Brownian noise in $dX(t)$ term in \eqref{eq:XL}
and $d\overline{X}(t)$ term in \eqref{Euler:continuous:X},
and their dynamics are different,
there does not exist a solution to equation (7.115) in Theorem 7.18 in \cite{liptser2013statistics},
and one can not apply Girsanov theorem to compute the relative entropy
between $(X(t),V(t))$ and $(\overline{V}(t),\overline{X}(t))$, 
which is the main technical difficulty here. 
To overcome this challenge, we define an auxiliary diffusion process $(\tilde{X}(t),\tilde{V}(t))$:
\begin{align}
&\tilde{V}(t)=V_{0}-\int_{0}^{t}\gamma\tilde{V}(\lfloor s/\eta\rfloor\eta)ds
\nonumber
\\
&\qquad\qquad\qquad
-\int_{0}^{t}g\left(X_{0}+\int_{0}^{\lfloor s/\eta\rfloor\eta}\tilde{V}(\lfloor u/\eta\rfloor\eta)du,\overline{U}_{\mathbf{z}}(s)\right)ds+\sqrt{2\gamma \beta^{-1}}\int_{0}^{t}dB(s),
\label{eq:tildeV}
\\
&\tilde{X}(t)=X_{0}+\int_{0}^{t}\tilde{V}(s)ds,
\label{eq:tildeX}
\end{align}
which serves as a bridge between the underdamped Langevin diffusion $(X(k\eta),V(k\eta))$
and the discrete time SGHMC1 iterates $(X_{k},V_{k})$.
Then, it is easy to see that $\tilde{V}(k\eta)$ has the same distribution as $V_{k}$, though $\tilde{X}(k\eta)$ is not distributed the same as $X_{k}$.
Since the drift term in \eqref{eq:tildeX} in the auxiliary diffusion process $(\tilde{X}(t),\tilde{V}(t))$ has the same dynamics as the $dX(t)$ term in \eqref{eq:XL},
Girsanov theorem is applicable according to Theorem 7.18 in \cite{liptser2013statistics}.

Let $\mathbb{P}$ be the probability measure associated with the underdamped Langevin diffusion $(X(t),V(t))$ in \eqref{eq:VL}--\eqref{eq:XL}
and $\tilde{\mathbb{P}}$ be the probability measure associated with the $(\tilde{X}(t),\tilde{V}(t))$ process in \eqref{eq:tildeV}--\eqref{eq:tildeX}.
Let $\mathcal{F}_{t}$ be the natural filtration up to time $t$.
Then, the Radon-Nikodym derivative of $\mathbb{P}$ w.r.t. $\tilde{\mathbb{P}}$
is given by the Girsanov theorem (see e.g. Section 7.6 in \cite{liptser2013statistics}):
\begin{align*}
\frac{d\mathbb{P}}{d\tilde{\mathbb{P}}}
\bigg|_{\mathcal{F}_{t}}
&=e^{-\sqrt{\frac{\beta}{2\gamma}}\int_{0}^{t}\left(\gamma\tilde{V}(s)-\gamma\tilde{V}(\lfloor s/\eta\rfloor\eta)
+\nabla F_{\mathbf{z}}(\tilde{X}(s))-g\left(X_{0}+\int_{0}^{\lfloor s/\eta\rfloor\eta}\tilde{V}(\lfloor u/\eta\rfloor\eta)du,\overline{U}_{\mathbf{z}}(s)\right)\right)\cdot dB(s)}
\\
&\qquad\qquad
\cdot e^{-\frac{\beta}{4\gamma}\int_{0}^{t}\left\Vert\gamma\tilde{V}(s)-\gamma\tilde{V}(\lfloor s/\eta\rfloor\eta)
+\nabla F_{\mathbf{z}}(\tilde{X}(s))-g\left(X_{0}+\int_{0}^{\lfloor s/\eta\rfloor\eta}\tilde{V}(\lfloor u/\eta\rfloor\eta)du,\overline{U}_{\mathbf{z}}(s)\right)\right\Vert^{2}ds}.
\end{align*}
Then by writing $\mathbb{P}_{t}$ and $\tilde{\mathbb{P}}_{t}$
as the probability measures $\mathbb{P}$ and $\tilde{\mathbb{P}}$ conditional on the filtration $\mathcal{F}_{t}$, 
\begin{align*}
&D(\tilde{\mathbb{P}}_{t}\Vert\mathbb{P}_{t})
\\
&:=-\int d\tilde{\mathbb{P}}_{t}\log\frac{d\mathbb{P}_{t}}{d\tilde{\mathbb{P}}_{t}}
\\
&=\frac{\beta}{4\gamma}\int_{0}^{t}\mathbb{E}_{\mathbf{z}}\left\Vert
\gamma\tilde{V}(s)-\gamma\tilde{V}(\lfloor s/\eta\rfloor\eta)
+\nabla F_{\mathbf{z}}(\tilde{X}(s))-g\left(X_{0}+\int_{0}^{\lfloor s/\eta\rfloor\eta}\tilde{V}(\lfloor u/\eta\rfloor\eta)du,\overline{U}_{\mathbf{z}}(s)\right)\right\Vert^{2}ds
\\
&\leq\frac{\beta}{2\gamma}\int_{0}^{t}\mathbb{E}_{\mathbf{z}}\left\Vert\nabla F_{\mathbf{z}}\left(X_{0}+\int_{0}^{\lfloor s/\eta\rfloor\eta}\tilde{V}(u)du\right)-g\left(X_{0}+\int_{0}^{\lfloor s/\eta\rfloor\eta}\tilde{V}(\lfloor u/\eta\rfloor\eta)du,\overline{U}_{\mathbf{z}}(s)\right)\right\Vert^{2}ds
\\
&\qquad\qquad
+\frac{\beta}{2\gamma}\int_{0}^{t}\mathbb{E}_{\mathbf{z}}\left\Vert
\gamma\tilde{V}(s)-\gamma\tilde{V}(\lfloor s/\eta\rfloor\eta)\right\Vert^{2}ds
\\
&\leq
\frac{\beta}{\gamma}\int_{0}^{t}\mathbb{E}_{\mathbf{z}}\left\Vert\nabla F_{\mathbf{z}}\left(X_{0}+\int_{0}^{\lfloor s/\eta\rfloor\eta}\tilde{V}(u)du\right)-\nabla F_{\mathbf{z}}\left(X_{0}+\int_{0}^{\lfloor s/\eta\rfloor\eta}\tilde{V}(\lfloor u/\eta\rfloor\eta)du\right)\right\Vert^{2}ds
\\
&\qquad
+\frac{\beta}{\gamma}\int_{0}^{t}\mathbb{E}_{\mathbf{z}}\left\Vert\nabla F_{\mathbf{z}}\left(X_{0}+\int_{0}^{\lfloor s/\eta\rfloor\eta}\tilde{V}(\lfloor u/\eta\rfloor\eta)du\right)
-g\left(X_{0}+\int_{0}^{\lfloor s/\eta\rfloor\eta}\tilde{V}(\lfloor u/\eta\rfloor\eta)du,\overline{U}_{\mathbf{z}}(s)\right)\right\Vert^{2}ds
\\
&\qquad\qquad
+\frac{\beta}{2\gamma}\int_{0}^{t}\mathbb{E}_{\mathbf{z}}\left\Vert
\gamma\tilde{V}(s)-\gamma\tilde{V}(\lfloor s/\eta\rfloor\eta)\right\Vert^{2}ds,
\end{align*}
which implies that
\begin{align}
&D(\tilde{\mathbb{P}}_{k\eta}\Vert\mathbb{P}_{k\eta})
\nonumber
\\
&\leq
\frac{\beta\eta}{\gamma}\sum_{j=0}^{k-1}
\mathbb{E}_{\mathbf{z}}\left\Vert\nabla F_{\mathbf{z}}\left(X_{0}+\int_{0}^{j\eta}\tilde{V}(u)du\right)-\nabla F_{\mathbf{z}}\left(X_{0}+\int_{0}^{j\eta}\tilde{V}(\lfloor u/\eta\rfloor\eta)du\right)\right\Vert^{2}
\nonumber 
\\
&\qquad
+\frac{\beta\eta}{\gamma}\sum_{j=0}^{k-1}
\mathbb{E}_{\mathbf{z}}\left\Vert\nabla F_{\mathbf{z}}\left(X_{0}+\int_{0}^{j\eta}\tilde{V}(\lfloor u/\eta\rfloor\eta)du\right)
-g\left(X_{0}+\int_{0}^{j\eta}\tilde{V}(\lfloor u/\eta\rfloor\eta)du,U_{\mathbf{z},j}\right)\right\Vert^{2}
\nonumber
\\
&\qquad\qquad
+\frac{\beta}{2\gamma}\sum_{j=0}^{k-1}\int_{j\eta}^{(j+1)\eta}\mathbb{E}_{\mathbf{z}}\left\Vert
\gamma\tilde{V}(s)-\gamma\tilde{V}(\lfloor s/\eta\rfloor\eta)\right\Vert^{2}ds. \label{eq:rel-entropy}
\end{align}
We first bound the first term in \eqref{eq:rel-entropy}:
\begin{align*}
&\frac{\beta\eta}{\gamma}\sum_{j=0}^{k-1}
\mathbb{E}_{\mathbf{z}}\left\Vert\nabla F_{\mathbf{z}}\left(X_{0}+\int_{0}^{j\eta}\tilde{V}(u)du\right)-\nabla F_{\mathbf{z}}\left(X_{0}+\int_{0}^{j\eta}\tilde{V}(\lfloor u/\eta\rfloor\eta)du\right)\right\Vert^{2}
\\
&
\leq
M^{2}\frac{\beta\eta}{\gamma}\sum_{j=0}^{k-1}
\mathbb{E}_{\mathbf{z}}\left\Vert\int_{0}^{j\eta}\left(\tilde{V}(u)-\tilde{V}(\lfloor u/\eta\rfloor\eta)\right)du\right\Vert^{2}
\\
&
\leq
M^{2}\frac{\beta\eta}{\gamma}\sum_{j=0}^{k-1}j\eta
\int_{0}^{j\eta}\mathbb{E}_{\mathbf{z}}\left\Vert\tilde{V}(u)-\tilde{V}(\lfloor u/\eta\rfloor\eta)\right\Vert^{2}du
\\
&=M^{2}\frac{\beta\eta}{\gamma}\sum_{j=0}^{k-1}j\eta\sum_{i=0}^{j-1}
\int_{i\eta}^{(i+1)\eta}\mathbb{E}_{\mathbf{z}}\left\Vert\tilde{V}(u)-\tilde{V}(\lfloor u/\eta\rfloor\eta)\right\Vert^{2}du\,,
\end{align*}
where we used part (ii) of Assumption~\ref{assumptions} Cauchy-Schwarz inequality.

For $i\eta<u\leq(i+1)\eta$, we have
\begin{equation}
\tilde{V}(u)-\tilde{V}(\lfloor u/\eta\rfloor\eta)
=-(u-i\eta)\gamma V_{i}
-(u-i\eta) g\left(X_{i},U_{\mathbf{z},i}\right)
+\sqrt{2\gamma\beta^{-1}}(B(u)-B(i\eta)),
\end{equation}
in distribution. Therefore,
\begin{align}
&\mathbb{E}_{\mathbf{z}}\left\Vert\tilde{V}(u)-\tilde{V}(\lfloor u/\eta\rfloor\eta)\right\Vert^{2}
\nonumber
\\
&=(u-i\eta)^{2}\mathbb{E}_{\mathbf{z}}\left\Vert\gamma V_{i}
+g\left(X_{i},U_{\mathbf{z},i}\right)\right\Vert^{2}
+2\gamma\beta^{-1}(u-i\eta)
\nonumber
\\
&=(u-i\eta)^{2}\mathbb{E}_{\mathbf{z}}\left\Vert\gamma V_{i}
+\nabla F_{\mathbf{z}}(X_{i})\right\Vert^{2}
+(u-i\eta)^{2}\mathbb{E}_{\mathbf{z}}\left\Vert \nabla F_{\mathbf{z}}(X_{i})-g\left(X_{i},U_{\mathbf{z},i}\right)\right\Vert^{2}
+2\gamma\beta^{-1}(u-i\eta)
\nonumber
\\
&\leq
2\eta^{2}\mathbb{E}_{\mathbf{z}}\left\Vert\gamma V_{i}\right\Vert^{2}
+2\eta^{2}\mathbb{E}_{\mathbf{z}}\left\Vert\nabla F_{\mathbf{z}}(X_{i})\right\Vert^{2}
+\eta^{2}2\delta(M^{2}\mathbb{E}_{\mathbf{z}}\left\Vert X_{i}\right\Vert^{2}+B^{2})
+2\gamma\beta^{-1}\eta
\nonumber
\\
&\leq
2\gamma^{2}\eta^{2}\mathbb{E}_{\mathbf{z}}\left\Vert V_{i}\right\Vert^{2}
+4\eta^{2}\left(M^{2}\mathbb{E}_{\mathbf{z}}\left\Vert X_{i}\right\Vert^{2}+B^{2}\right)
+\eta^{2}2\delta(M^{2}\mathbb{E}_{\mathbf{z}}\left\Vert X_{i}\right\Vert^{2}+B^{2})
+2\gamma\beta^{-1}\eta.\label{V:difference}
\end{align}
This implies that
\begin{align*}
&M^{2}\frac{\beta\eta}{\gamma}\sum_{j=0}^{k-1}j\eta\sum_{i=0}^{j-1}
\int_{i\eta}^{(i+1)\eta}\mathbb{E}_{\mathbf{z}}\left\Vert\tilde{V}(u)-\tilde{V}(\lfloor u/\eta\rfloor\eta)\right\Vert^{2}du
\\
&\leq
M^{2}\frac{\beta}{\gamma}(k\eta)^{3}
\left(2\gamma^{2}\eta^{2}\sup_{j\geq 0}\mathbb{E}_{\mathbf{z}}\left\Vert V_{j}\right\Vert^{2}
+(4+2\delta)\eta^{2}\left(M^{2}\sup_{j\geq 0}\mathbb{E}_{\mathbf{z}}\left\Vert X_{j}\right\Vert^{2}+B^{2}\right)
+2\gamma\beta^{-1}\eta\right).
\end{align*}

We can also bound the second term in \eqref{eq:rel-entropy}:
\begin{align*}
&\frac{\beta\eta}{\gamma}\sum_{j=0}^{k-1}
\mathbb{E}_{\mathbf{z}}\left\Vert\nabla F_{\mathbf{z}}\left(X_{0}+\int_{0}^{j\eta}\tilde{V}(\lfloor u/\eta\rfloor\eta)du\right)
-g\left(X_{0}+\int_{0}^{j\eta}\tilde{V}(\lfloor u/\eta\rfloor\eta)du,U_{\mathbf{z},j}\right)\right\Vert^{2}
\\
&\leq\frac{\beta\eta}{\gamma}
\sum_{j=0}^{k-1}2\delta\left(M^{2}\mathbb{E}_{\mathbf{z}}\left\Vert X_{0}+\int_{0}^{s}\tilde{V}(\lfloor u/\eta\rfloor\eta)du\right\Vert^{2}+B^{2}\right)
\\
&=
\frac{\beta\eta}{\gamma}
\sum_{j=0}^{k-1}2\delta\left(M^{2}\mathbb{E}_{\mathbf{z}}\left\Vert X_{j}\right\Vert^{2}+B^{2}\right)
\\
&\leq
\frac{2\beta\delta}{\gamma}k\eta
\left(M^{2}\sup_{j\geq 0}\mathbb{E}_{\mathbf{z}}\left\Vert X_{j}\right\Vert^{2}+B^{2}\right),
\end{align*}
where the first inequality follows from part (iv) of Assumption~\ref{assumptions}. 

Finally, let us bound the third term in \eqref{eq:rel-entropy} as follows:
\begin{align*}
&\frac{\beta}{2\gamma}\sum_{j=0}^{k-1}\int_{j\eta}^{(j+1)\eta}\mathbb{E}_{\mathbf{z}}\left\Vert
\gamma\tilde{V}(s)-\gamma\tilde{V}(\lfloor s/\eta\rfloor\eta)\right\Vert^{2}ds
\\
&\leq
\frac{\beta}{2\gamma}(k\eta)\gamma^{2}
\left(
2\gamma^{2}\eta^{2}C_{v}^{d}
+(4+2\delta)\eta^{2}\left(M^{2}C_{x}^{d}+B^{2}\right)
+2\gamma\beta^{-1}\eta\right),
\end{align*} 
where we used the estimate in \eqref{V:difference}.

Hence, together with Lemma~\ref{lem:L2bound}, we conclude that that
\begin{align*}
D(\tilde{\mathbb{P}}_{k\eta}\Vert\mathbb{P}_{k\eta})
&\leq
M^{2}\frac{\beta}{\gamma}(k\eta)^{3}
\left(2\gamma^{2}\eta^{2}C_{v}^{d}
+(4+2\delta)\eta^{2}\left(M^{2}C_{x}^{d}+B^{2}\right)
+2\gamma\beta^{-1}\eta\right)
\\
&\qquad\qquad\qquad
+\frac{2\beta\delta}{\gamma}k\eta
\left(M^{2}C_{x}^{d}+B^{2}\right)
\\
&\qquad
+\frac{\beta}{2\gamma}(k\eta)\gamma^{2}
\left(
2\gamma^{2}\eta^{2}C_{v}^{d}
+(4+2\delta)\eta^{2}\left(M^{2}C_{x}^{d}+B^{2}\right)
+2\gamma\beta^{-1}\eta\right).
\end{align*}
We can then apply the following result of \cite{BV}, 
that is, for any two Borel probability measures $\mu,\nu$ on $\mathbb{R}^{2d}$ with finite second moments,
\begin{equation*}
\mathcal{W}_{2}(\mu,\nu)\leq C_{\nu}
\left[\sqrt{D(\mu\Vert\nu)}+\left(\frac{D(\mu\Vert\nu)}{2}\right)^{1/4}\right],
\end{equation*}
where
\begin{equation*}
C_{\nu}=2\inf_{\lambda>0}\left(\frac{1}{\lambda}\left(\frac{3}{2}+\log\int_{\mathbb{R}^{2d}}e^{\lambda\Vert w\Vert^{2}}\nu(dw)\right)\right)^{1/2}.
\end{equation*}
From the exponential integrability of the measure $\nu_{\bz,k\eta}$ in Lemma~\ref{lem:exp-integrable}, we have
\begin{align*}
C_{\nu_{\bz,k\eta}}
&\leq
2\left(\frac{1}{\alpha_{0}}\left(\frac{3}{2}+\log\int_{\mathbb{R}^{2d}}e^{\alpha_{0}\Vert(x,v)\Vert^{2}}\nu_{\bz,k\eta}(dx,dv)\right)\right)^{1/2}
\\
&\leq
2\left(\frac{1}{\alpha_{0}}\left(\frac{3}{2}+\log\left(\int_{\mathbb{R}^{2d}}e^{\frac{1}{4}\alpha\mathcal{V}(x,v)}\mu_{0}(dx,dv)+\frac{1}{4}e^{\frac{\alpha(d+A)}{3\lambda}}\alpha\gamma(d+A)k\eta\right)\right)\right)^{1/2}.
\end{align*}
Hence
\begin{align}
\mathcal{W}_{2}^{2}(\tilde{\mathbb{P}}_{k\eta},\nu_{\bz,k\eta})
\nonumber
&\leq 
\frac{4}{\alpha_{0}}\left(\frac{3}{2}+\log\left(\int_{\mathbb{R}^{2d}}e^{\frac{1}{4}\alpha\mathcal{V}(x,v)}\mu_{0}(dx,dv)
+\frac{1}{4}e^{\frac{\alpha(d+A)}{3\lambda}}\alpha\gamma(d+A)k\eta\right)\right)
\\
&\qquad\qquad\qquad
\cdot\left[\sqrt{D(\tilde{\mathbb{P}}_{k\eta}\Vert\mathbb{P}_{k\eta})}+\left(\frac{D(\tilde{\mathbb{P}}_{k\eta}\Vert\mathbb{P}_{k\eta})}{2}\right)^{1/4}\right]^{2},
\label{W2estimate}
\end{align}
where
\begin{align}
D(\tilde{\mathbb{P}}_{k\eta}\Vert\mathbb{P}_{k\eta})
&\leq
(k\eta)^{3}\Bigg[\left(\frac{M^{2}\beta\eta}{\gamma}+\frac{\beta\eta\gamma}{2}\right)
\left(2\gamma^{2}\eta C_{v}^{d}
+(4+2\delta)\eta\left(M^{2}C_{x}^{d}+B^{2}\right)
+2\gamma\beta^{-1}\right)
\nonumber
\\
&\qquad\qquad\qquad
+\frac{\beta\delta}{\gamma}
\left(M^{2}C_{x}^{d}+B^{2}\right)\Bigg].
\nonumber
\end{align}
Note that $\eta\leq 1$ so that
$2\gamma^{2}\eta C_{v}^{d}
+(4+2\delta)\eta\left(M^{2}C_{x}^{d}+B^{2}\right)
+2\gamma\beta^{-1}\leq(C_{2})^{2}$,
where $C_{2}$ is defined in \eqref{def-hat-Ctwo}.
Then, we have
\begin{equation*}
D(\tilde{\mathbb{P}}_{k\eta}\Vert\mathbb{P}_{k\eta})
\leq
(k\eta)^{3}\Bigg[\left(\frac{M^{2}\beta\eta}{\gamma}+\frac{\beta\eta\gamma}{2}\right)(C_{2})^{2}
+\frac{\beta\delta}{\gamma}
\left(M^{2}C_{x}^{d}+B^{2}\right)\Bigg].
\end{equation*}

By using $(x+y)^{2}\leq 2(x^{2}+y^{2})$, we get
\begin{align}
\mathcal{W}_{2}^{2}(\tilde{\mathbb{P}}_{k\eta},\nu_{\bz,k\eta})
\nonumber
&\leq 
\frac{8}{\alpha_{0}}\left(\frac{3}{2}+\log\left(\int_{\mathbb{R}^{2d}}e^{\frac{1}{4}\alpha\mathcal{V}(x,v)}\mu_{0}(dx,dv)
+\frac{1}{4}e^{\frac{\alpha(d+A)}{3\lambda}}\alpha\gamma(d+A)k\eta\right)\right)
\nonumber
\\
&\qquad\qquad\qquad
\cdot\left[D(\tilde{\mathbb{P}}_{k\eta}\Vert\mathbb{P}_{k\eta})+\sqrt{D(\tilde{\mathbb{P}}_{k\eta}\Vert\mathbb{P}_{k\eta})}\right].
\label{eq:W2estimate}
\end{align}

Since $k\eta\geq e>1$, we get
\begin{align}
&\frac{8}{\alpha_{0}}\left(\frac{3}{2}+\log\left(\int_{\mathbb{R}^{2d}}e^{\frac{1}{4}\alpha\mathcal{V}(x,v)}\mu_{0}(dx,dv)
+\frac{1}{4}e^{\frac{\alpha(d+A)}{3\lambda}}\alpha\gamma(d+A)k\eta\right)\right)
\nonumber
\\
&\leq\frac{8}{\alpha_{0}}\left(\frac{3}{2}+\log\left(\int_{\mathbb{R}^{2d}}e^{\frac{1}{4}\alpha\mathcal{V}(x,v)}\mu_{0}(dx,dv)
+\frac{1}{4}e^{\frac{\alpha(d+A)}{3\lambda}}\alpha\gamma(d+A)\right)
+\log(k\eta)\right)
\nonumber
\\
&\leq\frac{8}{\alpha_{0}}\left(\frac{3}{2}+\log\left(\int_{\mathbb{R}^{2d}}e^{\frac{1}{4}\alpha\mathcal{V}(x,v)}\mu_{0}(dx,dv)
+\frac{1}{4}e^{\frac{\alpha(d+A)}{3\lambda}}\alpha\gamma(d+A)\right)
+1\right)\log(k\eta)\,,
\label{eq:log}
\end{align}
and
\begin{align*}
&D(\tilde{\mathbb{P}}_{k\eta}\Vert\mathbb{P}_{k\eta})+\sqrt{D(\tilde{\mathbb{P}}_{k\eta}\Vert\mathbb{P}_{k\eta})}
\\
&\leq
\left(\left(\frac{M^{2}\beta\eta}{\gamma}+\frac{\beta\eta\gamma}{2}\right)(C_{2})^{2}+\sqrt{\left(\frac{M^{2}\beta\eta}{\gamma}+\frac{\beta\eta\gamma}{2}\right)(C_{2})^{2}}\right)(k\eta)^{3}\eta^{1/2}
\\
&\qquad\qquad
+\left(\left(M^{2}C_{x}^{d}+B^{2}\right)\frac{\beta}{\gamma}+\sqrt{\left(M^{2}C_{x}^{d}+B^{2}\right)\frac{\beta}{\gamma}}\right)(k\eta)^{3}\sqrt{\delta},
\end{align*}
which implies that
\begin{equation*}
\mathcal{W}_{2}^{2}(\tilde{\mathbb{P}}_{k\eta},\nu_{\bz,k\eta})
\leq(C_{0}^{2}\sqrt{\delta}+C_{1}^{2}\sqrt{\eta})(k\eta)^{3}\log(k\eta),
\end{equation*}
where $C_{0}$ and $C_{1}$ are defined in \eqref{def-hat-Czero} and \eqref{def-hat-Cone}. The result then follows from the fact that $\sqrt{x+y} \le \sqrt{x} + \sqrt{y}$ for non-negative real numbers $x$ and $y$.

Finally, let us provide a bound on $\mathcal{W}_{2}(\mu_{\bz,k},\tilde{\mathbb{P}}_{k\eta})$.
Note that by the definition of $\tilde{V}$, 
we have that $\left(X_{0}+\int_{0}^{k\eta}\tilde{V}(\lfloor s/\eta\rfloor\eta)ds,\tilde{V}(k\eta)\right)$
has the same law as $\mu_{\bz,k}$, and we can compute that
\begin{align*}
&\mathbb{E}_{\mathbf{z}}\left\Vert\tilde{X}(k\eta)-X_{0}-\int_{0}^{k\eta}\tilde{V}(\lfloor s/\eta\rfloor\eta)ds\right\Vert^{2}
\\
&=\mathbb{E}_{\mathbf{z}}\left\Vert\int_{0}^{k\eta}\tilde{V}(s)-\tilde{V}(\lfloor s/\eta\rfloor\eta)ds\right\Vert^{2}
\\
&\leq
k\eta\int_{0}^{k\eta}\mathbb{E}_{\mathbf{z}}\left\Vert\tilde{V}(s)-\tilde{V}(\lfloor s/\eta\rfloor\eta)\right\Vert^{2}ds
\\
&\leq
(k\eta)^{2}\eta\left(2\gamma^{2}\eta C_{v}^{d}+(4+2\delta)\eta\left(M^{2}C_{x}^{d}+B^{2}\right)+2\gamma\beta^{-1}\right)
\leq(k\eta)^{2}\eta(C_{2})^{2}.
\end{align*}
where we used the assumption $\eta\leq 1$
so that $2\gamma^{2}\eta C_{v}^{d}
+(4+2\delta)\eta\left(M^{2}C_{x}^{d}+B^{2}\right)
+2\gamma\beta^{-1}\leq(C_{2})^{2}$ in the last inequality above, 
where $C_{2}$ is defined in \eqref{def-hat-Ctwo}.
Therefore,
\begin{equation*}
\mathcal{W}_{2}(\mu_{\bz,k},\tilde{\mathbb{P}}_{k\eta})
\leq
C_{2}k\eta\sqrt{\eta}.
\end{equation*}
The proof is complete.

\subsection{Proof of Lemma~\ref{lem:nu0-piz}}
We recall first from \eqref{eq:V-estimate} that
\begin{equation*}
\mathcal{V}(x,v)
\geq\max\left\{\frac{1}{8} (1-2 \lambda) \beta \gamma^2 \Vert x\Vert^2,\frac{\beta}{4}(1-2\lambda)\Vert v\Vert^{2}\right\}.
\end{equation*}
Since $\int_{\mathbb{R}^{2d}}e^{\alpha\mathcal{V}(x,v)}\mu_{0}(dx,dv)<\infty$
with $\alpha>0$, we have $\Vert(x,v)\Vert_{L^{2}(\mu_{0})}<\infty$.

Next, let us notice that by the concavity of the function $h$, we have (see \cite{Eberle})
\begin{equation*}
h(r) \le \min\{r, h(R_1) \} \le \min\{r,R_1 \}, \quad \text{for any $r \ge 0$}.
\end{equation*}
 It follows that
\begin{align*}
\rho((x,v),(x',v'))
&\leq\min\{r((x,v),(x',v')),R_{1}\}(1+\varepsilon_{1}\mathcal{V}(x,v)+\varepsilon_{1}\mathcal{V}(x',v'))
\\
&\leq R_{1}(1+\varepsilon_{1}\mathcal{V}(x,v)+\varepsilon_{1}\mathcal{V}(x',v')).
\end{align*}
Moreover, by the definition of $\mathcal{V}$ in \eqref{eq:lyapunov} and Lemma~\ref{lem:gradient-bound}, we deduce that
\begin{align*}
\mathcal{V}(x,v)
&\leq \beta \left(
\frac{M}{2}\Vert x\Vert^{2}+B\Vert x\Vert+A_{0} \right)
+\frac{1}{4}\beta \gamma^{2}(\Vert x+\gamma^{-1}v\Vert^{2}+\Vert\gamma^{-1}v\Vert^{2}-\lambda\Vert x\Vert^{2})
\\
&\leq \beta \left(
\frac{M}{2}\Vert x\Vert^{2}+B\Vert x\Vert+A_{0} \right)
+\frac{1}{4} \beta \gamma^{2}(2\Vert x\Vert^{2}+2\gamma^{-2}\Vert v\Vert^{2}+\Vert\gamma^{-1}v\Vert^{2}-\lambda\Vert x\Vert^{2})
\\
&\leq \beta \left(
M\Vert x\Vert^{2}+A_{0}+\frac{B^{2}}{2M} \right)
+\frac{1}{4}\beta \gamma^{2}(2\Vert x\Vert^{2}+2\gamma^{-2}\Vert v\Vert^{2}+\Vert\gamma^{-1}v\Vert^{2}-\lambda\Vert x\Vert^{2})
\\
&\leq
\left(\beta M+\frac{1}{2} \beta \gamma^{2}\right)\Vert x\Vert^{2}+\frac{3}{4}
\beta \Vert v\Vert^{2}+\beta A_{0}+\frac{\beta B^{2}}{2M}.
\end{align*}
Therefore, we obtain
\begin{align} \label{eq:H-rho}
&\mathcal{H}_{\rho}(\mu_{0},\pi_{\mathbf{z}})
\nonumber \\
&\leq
R_{1}+R_{1}\varepsilon_{1}\left(\left(M+\frac{1}{2}\beta \gamma^{2}\right)\int_{\mathbb{R}^{2d}}\Vert x\Vert^{2}\mu_{0}(dx,dv)
+\frac{3}{4}\beta \int_{\mathbb{R}^{2d}}\Vert v\Vert^{2}\mu_{0}(dx,dv)+\beta A_{0}+\frac{\beta B^{2}}{2M}\right)
\nonumber \\
&\qquad
+R_{1}\varepsilon_{1}\left(\left(M+\frac{1}{2}\beta \gamma^{2}\right)\int_{\mathbb{R}^{2d}}\Vert x\Vert^{2}\pi_{\mathbf{z}}(dx,dv)
+\frac{3}{4}\beta \int_{\mathbb{R}^{2d}}\Vert v\Vert^{2}\pi_{\mathbf{z}}(dx,dv)+\beta A_{0}+\frac{\beta B^{2}}{2M}\right).
\end{align}
It has been shown in \cite[Section 3.5]{Raginsky} that
\begin{equation*}
\int_{\mathbb{R}^{2d}}\Vert x\Vert^{2}\pi_{\mathbf{z}}(dx,dv)
\le \frac{b+d/\beta}{m}.
\end{equation*}
In addition, from the explicit expression of $\pi_{\mathbf{z}}(dx,dv)$ in \eqref{eq:inv-distr}, we have
\begin{equation*}
\int_{\mathbb{R}^{2d}}\Vert v\Vert^{2}\pi_{\mathbf{z}}(dx,dv)
=(2\pi \beta^{-1})^{-d/2}\int_{\mathbb{R}^{d}}\Vert v\Vert^{2}e^{-\Vert v\Vert^{2}/(2\beta^{-1})}dv = d/\beta.
\end{equation*}
Hence, the conclusion follows from \eqref{eq:H-rho}.

\section{Proofs of Lemmas in Section~\ref{sec:proof:2}} \label{lemma-proof-2}

\subsection{Proof of Lemma~\ref{lem:L2bound:Jordan}}

Before we proceed to the proof of Lemma~\ref{lem:L2bound:Jordan},
let us state two technical lemmas, which will be used in
the proof of Lemma~\ref{lem:L2bound:Jordan}. Recall $\psi_{0}(t)=e^{-\gamma t}$ and $\psi_{k+1}(t)=\int_{0}^{t}\psi_{k}(s)ds$, and $(\xi_{k+1},\xi'_{k+1})$ is a $2d$-dimensional centered Gaussian vector from the SGHMC2 iterates $(\hat{X}_{k}, \hat{V}_{k})$ given in \eqref{SGHMC2:V}--\eqref{SGHMC2:X}. Using the definitions, it is straightforward to establish
these two lemmas, so we omit the details of their proofs.

\begin{lemma}\label{tech:lem:2}
For any $\eta\geq 0$,
\begin{equation}
\max\left\{|\psi_{0}(\eta)-1+\gamma\eta|,|\eta-\psi_{1}(\eta)|,|\psi_{2}(\eta)|\right\}
\leq c_{0}\eta^{2},
\end{equation}
where $c_{0}:=1+\gamma^{2}$.
\end{lemma}

\begin{lemma}\label{tech:lem:1}
For any $\eta\geq 0$, 
\begin{align}
&C_{11}(\eta):=\mathbb{E}\Vert\xi_{k}\Vert^{2}\leq c_{11}\eta:=d\eta,
\\
&C_{22}(\eta):=\mathbb{E}\Vert\xi'_{k}\Vert^{2}\leq c_{22}\eta^{3}:=\frac{d}{3}\eta^{3},
\\
&C_{12}(\eta):=\mathbb{E}\langle\xi_{k},\xi'_{k}\rangle\leq c_{12}\eta^{2}:=\frac{d}{2}\eta^{2}.
\end{align}
\end{lemma}

Now, we are ready to prove Lemma~\ref{lem:L2bound:Jordan},
i.e. the uniform (in time) $L^2$ bounds for $(\hat{X}_k, \hat{V}_k)$ defined in \eqref{SGHMC2:V}--\eqref{SGHMC2:X}.  
We can rewrite the dynamics of the SGHMC2 iterates as follows:
\begin{align}
&\hat{V}_{k+1}=(1-\gamma\eta)\hat{V}_{k}-\eta g(\hat{X}_{k},U_{\mathbf{z},k})
+\hat{E}_{k}
+\sqrt{2\gamma\beta^{-1}}\xi_{k+1},
\\
&\hat{X}_{k+1}=\hat{X}_{k}+\eta\hat{V}_{k}
+\hat{E}'_{k}+\sqrt{2\gamma\beta^{-1}}\xi'_{k+1},
\end{align}
where
\begin{align}
&\hat{E}_{k}:=(\psi_{0}(\eta)-1+\gamma\eta)\hat{V}_{k}+(\eta-\psi_{1}(\eta))g(\hat{X}_{k},U_{\mathbf{z},k}),
\\
&\hat{E}'_{k}:=(\psi_{1}(\eta)-\eta)\hat{V}_{k}-\psi_{2}(\eta)g(\hat{X}_{k},U_{\mathbf{z},k}),
\end{align}
where $\mathbb{E}g(x,U_{{\mathbf{z}}, k}) = \nabla F_{\mathbf{z}} (x)$ for any $x$. We again use the Lyapunov function $\mathcal{V}(x,v)$ 
defined as before, and set for each $k=0,1, \ldots,$
\begin{equation} \label{eq:L2k:hat}
\hat{L}_2(k) = \mathbb{E}_{\mathbf{z}} \mathcal{V}(\hat{X}_k, \hat{V}_k)/\beta= \mathbb{E}_{\mathbf{z}}  \left[F_{\mathbf{z}}(\hat{X}_k)
+\frac{1}{4}\gamma^{2}\left(\Vert \hat{X}_k +\gamma^{-1}\hat{V}_k \Vert^{2}+\Vert\gamma^{-1}\hat{V}_k \Vert^{2}-\lambda\Vert \hat{X}_k\Vert^{2}\right) \right].
\end{equation}

We can compute that
\begin{align*}
\mathbb{E}F_{\mathbf{z}}(\hat{X}_{k+1})
&=\mathbb{E}F_{\mathbf{z}}\left(\hat{X}_{k}+\eta\hat{V}_{k}
+\hat{E}'_{k}+\sqrt{2\gamma\beta^{-1}}\xi'_{k+1}\right)
\\
&\leq\mathbb{E}F_{\mathbf{z}}\left(\hat{X}_{k}\right)
+\mathbb{E}\left\langle\nabla F_{\mathbf{z}}(\hat{X}_{k}),\eta\hat{V}_{k}
+\hat{E}'_{k}\right\rangle
+\frac{M}{2}\mathbb{E}\left\Vert\eta\hat{V}_{k}
+\hat{E}'_{k}+\sqrt{2\gamma\beta^{-1}}\xi'_{k+1}\right\Vert^{2}
\\
&=
\mathbb{E}F_{\mathbf{z}}\left(\hat{X}_{k}\right)
+\mathbb{E}\left\langle\nabla F_{\mathbf{z}}(\hat{X}_{k}),\eta\hat{V}_{k}\right\rangle
+\frac{M}{2}\eta^{2}\mathbb{E}\left\Vert\hat{V}_{k}\right\Vert^{2}
+\delta_{1}(k),
\end{align*}
where
\begin{align}
\delta_{1}(k)&:=\mathbb{E}\left\langle\nabla F_{\mathbf{z}}(\hat{X}_{k}),\hat{E}'_{k}\right\rangle
+\frac{M}{2}\mathbb{E}\left\Vert\hat{E}'_{k}+\sqrt{2\gamma\beta^{-1}}\xi'_{k+1}\right\Vert^{2}
+M\mathbb{E}\left\langle\eta\hat{V}_{k},\hat{E}'_{k}+\sqrt{2\gamma\beta^{-1}}\xi'_{k+1}\right\rangle
\\
&=\mathbb{E}\left\langle\nabla F_{\mathbf{z}}(\hat{X}_{k}),\hat{E}'_{k}\right\rangle
+\frac{M}{2}\mathbb{E}\left\Vert\hat{E}'_{k}\right\Vert^{2}
+M\gamma\beta^{-1}C_{22}(\eta)
+M\mathbb{E}\left\langle\eta\hat{V}_{k},\hat{E}'_{k}\right\rangle.
\end{align}
We can also compute that
\begin{align*}
&\frac{1}{4}\gamma^{2}\mathbb{E}\left\Vert\hat{X}_{k+1}+\gamma^{-1}\hat{V}_{k+1}\right\Vert^{2}
\\
&=\frac{1}{4}\mathbb{E}\left\Vert\gamma\hat{X}_{k+1}+\hat{V}_{k+1}\right\Vert^{2}
\\
&=\frac{1}{4}\mathbb{E}\left\Vert\left(\gamma\hat{X}_{k}+\hat{V}_{k}-\eta g(\hat{X}_{k},U_{\mathbf{z},k})
+\sqrt{2\gamma\beta^{-1}}\xi_{k+1}\right)
+\gamma\hat{E}'_{k}+\gamma\sqrt{2\gamma\beta^{-1}}\xi'_{k+1}
+\hat{E}_{k}\right\Vert^{2}
\\
&=\frac{1}{4}\mathbb{E}\left\Vert\gamma\hat{X}_{k}+\hat{V}_{k}-\eta g(\hat{X}_{k},U_{\mathbf{z},k})\right\Vert^{2}
+\delta_{2}(k),
\end{align*}
where
\begin{align*}
\delta_{2}(k)&:=\frac{1}{2}\gamma\beta^{-1}C_{11}(\eta)
+\frac{1}{2}\gamma^{3}\beta^{-1}C_{22}(\eta)
+\gamma^{2}\beta^{-1}C_{12}(\eta)
\\
&\qquad\qquad\qquad
+\frac{1}{4}\mathbb{E}\left\Vert\gamma\hat{E}'_{k}+\hat{E}_{k}\right\Vert^{2}
+\frac{1}{2}\mathbb{E}\left\langle\gamma\hat{X}_{k}+\hat{V}_{k}-\eta g(\hat{X}_{k},U_{\mathbf{z},k}),\gamma\hat{E}'_{k}+\hat{E}_{k}\right\rangle.
\end{align*}
We can also compute that
\begin{align*}
\frac{1}{4}\gamma^{2}\mathbb{E}\left\Vert\gamma^{-1}\hat{V}_{k+1}\right\Vert^{2}
&=\frac{1}{4}\mathbb{E}\left\Vert\hat{V}_{k+1}\right\Vert^{2}
\\
&=\frac{1}{4}\mathbb{E}\left\Vert(1-\gamma\eta)\hat{V}_{k}-\eta g(\hat{X}_{k},U_{\mathbf{z},k})
+\hat{E}_{k}
+\sqrt{2\gamma\beta^{-1}}\xi_{k+1}\right\Vert^{2}
\\
&=\frac{1}{4}\mathbb{E}\left\Vert(1-\gamma\eta)\hat{V}_{k}-\eta g(\hat{X}_{k},U_{\mathbf{z},k})
+\hat{E}_{k}\right\Vert^{2}
+\frac{1}{2}\gamma\beta^{-1}C_{11}(\eta)
\\
&=\frac{1}{4}\mathbb{E}\left\Vert(1-\gamma\eta)\hat{V}_{k}-\eta g(\hat{X}_{k},U_{\mathbf{z},k})\right\Vert^{2}
+\delta_{3}(k),
\end{align*}
where
\begin{equation}
\delta_{3}(k):=\frac{1}{4}\mathbb{E}\left\Vert\hat{E}_{k}\right\Vert^{2}
+\frac{1}{2}\mathbb{E}\left\langle(1-\gamma\eta)\hat{V}_{k}-\eta g(\hat{X}_{k},U_{\mathbf{z},k}),\hat{E}_{k}\right\rangle
+\frac{1}{2}\gamma\beta^{-1}C_{11}(\eta).
\end{equation}

Finally, we can compute that
\begin{align*}
-\frac{1}{4}\gamma^{2}\lambda\mathbb{E}\left\Vert\hat{X}_{k+1}\right\Vert^{2}
&=-\frac{1}{4}\gamma^{2}\lambda\mathbb{E}\left\Vert\hat{X}_{k}+\eta\hat{V}_{k}
+\hat{E}'_{k}+\sqrt{2\gamma\beta^{-1}}\xi'_{k+1}\right\Vert^{2}
\\
&=-\frac{1}{4}\gamma^{2}\lambda\mathbb{E}\left\Vert\hat{X}_{k}+\eta\hat{V}_{k}\right\Vert^{2}
-\frac{1}{4}\gamma^{2}\lambda\mathbb{E}\left\Vert\hat{E}'_{k}+\sqrt{2\gamma\beta^{-1}}\xi'_{k+1}\right\Vert^{2}
\\
&\qquad\qquad
-\frac{1}{2}\gamma^{2}\lambda\mathbb{E}\left\langle\hat{X}_{k}+\eta\hat{V}_{k},\hat{E}'_{k}+\sqrt{2\gamma\beta^{-1}}\xi'_{k+1}\right\rangle
\\
&\leq
-\frac{1}{4}\gamma^{2}\lambda\mathbb{E}\left\Vert\hat{X}_{k}+\eta\hat{V}_{k}\right\Vert^{2}
+\delta_{4}(k),
\end{align*}
where
\begin{equation}
\delta_{4}(k):=-\frac{1}{2}\gamma^{2}\lambda\mathbb{E}\left\langle\hat{X}_{k}+\eta\hat{V}_{k},\hat{E}'_{k}\right\rangle.
\end{equation}
By following the proofs of the $L_{2}$ uniform bound for SGHMC1 iterates, 
we get
\begin{equation*}
\frac{\hat{L}_2(k+1) - \hat{L}_2(k)}{\eta} 
\leq \gamma (A/\beta- \lambda \hat{L}_2(k)) + (K_{1}\hat{L}_{2}(k)+K_{2})\cdot\eta
+\frac{\delta_{1}(k)+\delta_{2}(k)+\delta_{3}(k)+\delta_{4}(k)}{\eta},
\end{equation*}
where $K_1$ and $K_2$ are given in \eqref{def-K1} and \eqref{def-K2}.

Next, we can estimate that
\begin{align*}
\delta_{1}(k)&=\mathbb{E}\left\langle\nabla F_{\mathbf{z}}(\hat{X}_{k}),\hat{E}'_{k}\right\rangle
+\frac{M}{2}\mathbb{E}\left\Vert\hat{E}'_{k}\right\Vert^{2}
+M\mathbb{E}\left\langle\eta\hat{V}_{k},\hat{E}'_{k}\right\rangle
+M\gamma\beta^{-1}C_{22}(\eta)
\\
&\leq
c_{0}\eta^{2}\mathbb{E}\left[\Vert\nabla F_{\mathbf{z}}(\hat{X}_{k})\Vert\cdot\left(\Vert\hat{V}_{k}\Vert+\Vert g(\hat{X}_{k},U_{\mathbf{z},k})\Vert\right)\right]
+\frac{M}{2}c_{0}^{2}\eta^{4}\mathbb{E}\left(\Vert\hat{V}_{k}\Vert+\Vert g(\hat{X}_{k},U_{\mathbf{z},k})\Vert\right)^{2}
\\
&\qquad\qquad
+Mc_{0}\eta^{2}\mathbb{E}\left[\Vert\hat{V}_{k}\Vert\cdot\left(\Vert\hat{V}_{k}\Vert+\Vert g(\hat{X}_{k},U_{\mathbf{z},k})\Vert\right)\right]
+M\gamma\beta^{-1}c_{22}\eta^{3}
\\
&\leq
\frac{1}{2}c_{0}\eta^{2}\mathbb{E}\Vert\nabla F_{\mathbf{z}}(\hat{X}_{k})\Vert^{2}
+Mc_{0}^{2}\eta^{4}\mathbb{E}\Vert g(\hat{X}_{k},U_{\mathbf{z},k})\Vert^{2}
+M\gamma\beta^{-1}c_{22}\eta^{3}
\\
&\qquad\qquad
+\frac{1}{2}Mc_{0}\eta^{2}(1+2\eta^{2})\mathbb{E}\Vert\hat{V}_{k}\Vert^{2}
+\frac{1}{2}(M+1)c_{0}\eta^{2}\mathbb{E}\left(\Vert\hat{V}_{k}\Vert+\Vert g(\hat{X}_{k},U_{\mathbf{z},k})\Vert\right)^{2}
\\
&\leq
\frac{1}{2}c_{0}\eta^{2}\mathbb{E}\Vert\nabla F_{\mathbf{z}}(\hat{X}_{k})\Vert^{2}
+Mc_{0}^{2}\eta^{4}\mathbb{E}\Vert g(\hat{X}_{k},U_{\mathbf{z},k})\Vert^{2}
+M\gamma\beta^{-1}c_{22}\eta^{3}
\\
&\qquad
+\frac{1}{2}Mc_{0}\eta^{2}(1+2\eta^{2})\mathbb{E}\Vert\hat{V}_{k}\Vert^{2}
+(M+1)c_{0}\eta^{2}\mathbb{E}\Vert\hat{V}_{k}\Vert^{2}
+(M+1)c_{0}\eta^{2}\mathbb{E}\Vert g(\hat{X}_{k},U_{\mathbf{z},k})\Vert^{2}
\\
&=\frac{1}{2}c_{0}\eta^{2}\mathbb{E}\Vert\nabla F_{\mathbf{z}}(\hat{X}_{k})\Vert^{2}
+c_{0}\eta^{2}(Mc_{0}\eta^{2}+M+1)\mathbb{E}\Vert g(\hat{X}_{k},U_{\mathbf{z},k})\Vert^{2}
+M\gamma\beta^{-1}c_{22}\eta^{3}
\\
&\qquad
+\frac{1}{2}c_{0}\eta^{2}(M(1+2\eta^{2})+2M+2)\mathbb{E}\Vert\hat{V}_{k}\Vert^{2},
\end{align*}
and
\begin{align*}
\delta_{2}(k)&=\frac{1}{2}\gamma\beta^{-1}C_{11}(\eta)+\frac{1}{2}\gamma^{3}\beta^{-1}C_{22}(\eta)+\gamma^{2}\beta^{-1}C_{12}(\eta)
\\
&\qquad\qquad
+\frac{1}{4}\mathbb{E}\Vert\gamma\hat{E}'_{k}+\hat{E}_{k}\Vert^{2}
+\frac{1}{2}\mathbb{E}\left\langle\gamma\hat{X}_{k}+\hat{V}_{k}-\eta g(\hat{X}_{k},U_{\mathbf{z},k}),\gamma\hat{E}'_{k}+\hat{E}_{k}\right\rangle
\\
&\leq
\frac{1}{2}\gamma\beta^{-1}c_{11}\eta+\frac{1}{2}\gamma^{3}\beta^{-1}c_{22}\eta^{3}
+\gamma^{2}\beta^{-1}c_{12}\eta^{2}
\\
&\qquad\qquad
+\frac{1}{4}c_{0}^{2}\eta^{4}(1+\gamma)^{2}\mathbb{E}\left(\Vert\hat{V}_{k}\Vert+\Vert g(\hat{X}_{k},U_{\mathbf{z},k})\Vert\right)^{2}
\\
&\qquad\qquad\qquad
+\frac{1}{2}c_{0}\eta^{2}(1+\gamma)\mathbb{E}\left[\left\Vert\gamma\hat{X}_{k}+\hat{V}_{k}-\eta g(\hat{X}_{k},U_{\mathbf{z},k})\right\Vert
\cdot\left(\Vert\hat{V}_{k}\Vert+\Vert g(\hat{X}_{k},U_{\mathbf{z},k})\Vert\right)\right]
\\
&\leq
\frac{1}{2}\gamma\beta^{-1}c_{11}\eta+\frac{1}{2}\gamma^{3}\beta^{-1}c_{22}\eta^{3}
+\gamma^{2}\beta^{-1}c_{12}\eta^{2}
\\
&\qquad\qquad
+\frac{1}{4}c_{0}\eta^{2}(1+\gamma)(1+c_{0}\eta^{2}(1+\gamma))\mathbb{E}\left(\Vert\hat{V}_{k}\Vert+\Vert g(\hat{X}_{k},U_{\mathbf{z},k})\Vert\right)^{2}
\\
&\qquad\qquad\qquad
+\frac{1}{4}c_{0}\eta^{2}(1+\gamma)\mathbb{E}\left\Vert\gamma\hat{X}_{k}+\hat{V}_{k}-\eta g(\hat{X}_{k},U_{\mathbf{z},k})\right\Vert^{2}
\\
&\leq
\frac{1}{2}\gamma\beta^{-1}c_{11}\eta+\frac{1}{2}\gamma^{3}\beta^{-1}c_{22}\eta^{3}
+\gamma^{2}\beta^{-1}c_{12}\eta^{2}
\\
&\qquad\qquad
+\frac{1}{2}c_{0}\eta^{2}(1+\gamma)(1+c_{0}\eta^{2}(1+\gamma))\mathbb{E}\Vert\hat{V}_{k}\Vert^{2}
\\
&\qquad\qquad\qquad
+\frac{1}{2}c_{0}\eta^{2}(1+\gamma)(1+c_{0}\eta^{2}(1+\gamma))\mathbb{E}\Vert g(\hat{X}_{k},U_{\mathbf{z},k})\Vert^{2}
\\
&\qquad\qquad\qquad
+\frac{3}{4}c_{0}\eta^{2}(1+\gamma)\gamma^{2}\mathbb{E}\Vert\hat{X}_{k}\Vert^{2}
+\frac{3}{4}c_{0}\eta^{2}(1+\gamma)\mathbb{E}\Vert\hat{V}_{k}\Vert^{2}
\\
&\qquad\qquad
+\frac{3}{4}c_{0}\eta^{2}(1+\gamma)\eta^{2}\mathbb{E}\Vert g(\hat{X}_{k},U_{\mathbf{z},k})\Vert^{2}
\\
&=\frac{1}{2}\gamma\beta^{-1}c_{11}\eta+\frac{1}{2}\gamma^{3}\beta^{-1}c_{22}\eta^{3}
+\gamma^{2}\beta^{-1}c_{12}\eta^{2}
\\
&\qquad\qquad
+\frac{1}{2}c_{0}\eta^{2}(1+\gamma)\left(\frac{5}{2}+c_{0}\eta^{2}(1+\gamma)\right)\mathbb{E}\Vert\hat{V}_{k}\Vert^{2}
\\
&\qquad\qquad\qquad
+\frac{1}{2}c_{0}\eta^{2}(1+\gamma)\left(1+c_{0}\eta^{2}(1+\gamma)+\frac{3}{2}\eta^{4}\right)\mathbb{E}\Vert g(\hat{X}_{k},U_{\mathbf{z},k})\Vert^{2}
\\
&\qquad\qquad\qquad
+\frac{3}{4}c_{0}\eta^{2}(1+\gamma)\gamma^{2}\mathbb{E}\Vert\hat{X}_{k}\Vert^{2},
\end{align*}
and we can compute that
\begin{align*}
\delta_{3}(k)&=\frac{1}{4}\mathbb{E}\Vert\hat{E}_{k}\Vert^{2}
+\frac{1}{2}\mathbb{E}\left\langle(1-\gamma\eta)\hat{V}_{k}-\eta g(\hat{X}_{k},U_{\mathbf{z},k}),\hat{E}_{k}\right\rangle
+\frac{1}{2}\gamma\beta^{-1}C_{11}(\eta)
\\
&\leq\frac{1}{4}c_{0}^{2}\eta^{4}\mathbb{E}\left(\Vert\hat{V}_{k}\Vert+\Vert g(\hat{X}_{k},U_{\mathbf{z},k})\Vert\right)^{2}
\\
&\qquad
+\frac{1}{2}c_{0}\eta^{2}\mathbb{E}\left[\left\Vert(1-\gamma\eta)\hat{V}_{k}-\eta g(\hat{X}_{k},U_{\mathbf{z},k})\right\Vert
\cdot\left(\Vert\hat{V}_{k}\Vert+\Vert g(\hat{X}_{k},U_{\mathbf{z},k})\Vert\right)\right]
+\frac{1}{2}\gamma\beta^{-1}c_{11}\eta
\\
&\leq\frac{1}{4}c_{0}\eta^{2}(1+c_{0}\eta^{2})\mathbb{E}\left(\Vert\hat{V}_{k}\Vert+\Vert g(\hat{X}_{k},U_{\mathbf{z},k})\Vert\right)^{2}
\\
&\qquad
+\frac{1}{4}c_{0}\eta^{2}\mathbb{E}\left\Vert(1-\gamma\eta)\hat{V}_{k}-\eta g(\hat{X}_{k},U_{\mathbf{z},k})\right\Vert^{2}
+\frac{1}{2}\gamma\beta^{-1}c_{11}\eta
\\
&\leq\frac{1}{2}c_{0}\eta^{2}(1+c_{0}\eta^{2})\mathbb{E}\Vert\hat{V}_{k}\Vert^{2}
+\frac{1}{2}c_{0}\eta^{2}(1+c_{0}\eta^{2})\mathbb{E}\Vert g(\hat{X}_{k},U_{\mathbf{z},k})\Vert^{2}
\\
&\qquad
+\frac{1}{2}c_{0}\eta^{2}(1-\gamma\eta)^{2}\mathbb{E}\Vert\hat{V}_{k}\Vert^{2}
+\frac{1}{2}c_{0}\eta^{4}\mathbb{E}\Vert g(\hat{X}_{k},U_{\mathbf{z},k})\Vert^{2}
+\frac{1}{2}\gamma\beta^{-1}c_{11}\eta
\\
&=\frac{1}{2}c_{0}\eta^{2}(2-2\gamma\eta+(c_{0}+\gamma^{2})\eta^{2})\mathbb{E}\Vert\hat{V}_{k}\Vert^{2}
\\
&\qquad\qquad
+\frac{1}{2}c_{0}\eta^{2}(1+(c_{0}+1)\eta^{2})\mathbb{E}\Vert g(\hat{X}_{k},U_{\mathbf{z},k})\Vert^{2}
+\frac{1}{2}\gamma\beta^{-1}c_{11}\eta,
\end{align*}
and finally we can compute that
\begin{align*}
\delta_{4}(k)&=-\frac{1}{2}\gamma^{2}\lambda\mathbb{E}\left\langle\hat{X}_{k}+\eta\hat{V}_{k},\hat{E}'_{k}\right\rangle
\\
&\leq\frac{1}{2}\gamma^{2}\lambda c_{0}\eta^{2}\mathbb{E}
\left[\left\Vert\hat{X}_{k}+\eta\hat{V}_{k}\right\Vert\cdot
\left(\Vert\hat{V}_{k}\Vert+\Vert g(\hat{X}_{k},U_{\mathbf{z},k})\Vert\right)\right]
\\
&\leq
\frac{1}{4}\gamma^{2}\lambda c_{0}\eta^{2}\mathbb{E}
\Vert\hat{X}_{k}+\eta\hat{V}_{k}\Vert^{2}
+\frac{1}{4}\gamma^{2}\lambda c_{0}\eta^{2}\mathbb{E}
\left(\Vert\hat{V}_{k}\Vert+\Vert g(\hat{X}_{k},U_{\mathbf{z},k})\Vert\right)^{2}
\\
&\leq
\frac{1}{2}\gamma^{2}\lambda c_{0}\eta^{2}\mathbb{E}
\Vert\hat{X}_{k}\Vert^{2}
+\frac{1}{2}\gamma^{2}\lambda c_{0}\eta^{2}(1+\eta^{2})\mathbb{E}\Vert\hat{V}_{k}\Vert^{2}
+\frac{1}{2}\gamma^{2}\lambda c_{0}\eta^{2}\mathbb{E}\Vert g(\hat{X}_{k},U_{\mathbf{z},k})\Vert^{2}.
\end{align*}

Putting everything together, we have
\begin{align*}
&\frac{\hat{L}_2(k+1) - \hat{L}_2(k)}{\eta} 
\\
&\leq \gamma (A/\beta- \lambda \hat{L}_2(k)) + (K_{1}\hat{L}_{2}(k)+K_{2})\cdot\eta
+\frac{\delta_{1}(k)+\delta_{2}(k)+\delta_{3}(k)+\delta_{4}(k)}{\eta}
\\
&\leq
\gamma ((d+A)/\beta- \lambda \hat{L}_2(k)) + (K_{1}\hat{L}_{2}(k)+K_{2})\cdot\eta
\\
&\qquad
+\frac{1}{2}c_{0}\eta
\Bigg((M(1+2\eta^{2})+2M+4-2\gamma\eta+(c_{0}+\gamma^{2})\eta^{2})
\\
&\qquad\qquad
+(1+\gamma)\left(\frac{5}{2}+c_{0}\eta^{2}(1+\gamma)\right)
+\gamma^{2}\lambda(1+\eta^{2})\Bigg)\mathbb{E}\Vert\hat{V}_{k}\Vert^{2}
\\
&\qquad
+\frac{1}{2}c_{0}\eta
\Bigg((1+\gamma)\left(1+c_{0}\eta^{2}(1+\gamma)+\frac{3}{2}\eta^{4}\right)
+1+(c_{0}+1)\eta^{2}
\\
&\qquad\qquad\qquad
+\lambda\gamma^{2}
+2(Mc_{0}\eta^{2}+M+1)\Bigg)\mathbb{E}\Vert g(\hat{X}_{k},U_{\mathbf{z},k})\Vert^{2}
\\
&\qquad\qquad\qquad
+\frac{1}{2}c_{0}\eta\mathbb{E}\Vert\nabla F_{\mathbf{z}}(\hat{X}_{k})\Vert^{2}
+\frac{1}{2}\gamma^{2}c_{0}\eta\left(\lambda+\frac{3}{2}(1+\gamma)\right)\mathbb{E}\Vert\hat{X}_{k}\Vert^{2}
\\
&\qquad
+\frac{1}{2}\gamma^{3}\beta^{-1}c_{22}\eta^{2}
+\gamma^{2}\beta^{-1}c_{12}\eta
+M\gamma\beta^{-1}c_{22}\eta^{2},
\end{align*}
where we used the fact that $c_{11}=d$.
Moreover,
\begin{equation*}
\mathbb{E}\Vert\nabla F_{\mathbf{z}}(\hat{X}_{k})\Vert^{2}
\leq
\mathbb{E}(M\Vert\hat{X}_{k}\Vert+B)^{2}
\leq
2M^{2}\mathbb{E}\Vert\hat{X}_{k}\Vert^{2}+2B^{2},
\end{equation*}
and
\begin{align}
\mathbb{E}\Vert g(\hat{X}_{k},U_{\mathbf{z},k})\Vert^{2}
&=\mathbb{E}\Vert\nabla F_{\mathbf{z}}(\hat{X}_{k})\Vert^{2}
+\mathbb{E}\Vert g(\hat{X}_{k},U_{\mathbf{z},k})-\nabla F_{\mathbf{z}}(\hat{X}_{k})\Vert^{2}
\nonumber
\\
&\leq
\mathbb{E}\Vert\nabla F_{\mathbf{z}}(\hat{X}_{k})\Vert^{2}
+2\delta M^{2}\mathbb{E}\Vert\hat{X}_{k}\Vert^{2}
+2\delta B^{2}
\nonumber
\\
&\leq
2(1+\delta)M^{2}\mathbb{E}\Vert\hat{X}_{k}\Vert^{2}
+2(1+\delta)B^{2}.\label{g:estimate}
\end{align}
Therefore, we have
\begin{align*}
\frac{\hat{L}_2(k+1) - \hat{L}_2(k)}{\eta} 
&\leq
\gamma ((d+A)/\beta- \lambda \hat{L}_2(k)) + (K_{1}\hat{L}_{2}(k)+K_{2})\cdot\eta
\\
&\qquad
+\frac{1}{2}c_{0}\eta
\Bigg((M(1+2\eta^{2})+2M+4-2\gamma\eta+(c_{0}+\gamma^{2})\eta^{2})
\\
&\qquad\qquad
+(1+\gamma)\left(\frac{5}{2}+c_{0}\eta^{2}(1+\gamma)\right)
+\gamma^{2}\lambda(1+\eta^{2})\Bigg)\mathbb{E}\Vert\hat{V}_{k}\Vert^{2}
\\
&\qquad
+\frac{1}{2}c_{0}\eta
\Bigg[\Bigg((1+\gamma)\left(1+c_{0}\eta^{2}(1+\gamma)+\frac{3}{2}\eta^{4}\right)
+1+(c_{0}+1)\eta^{2}
\\
&\qquad\qquad\qquad
+\lambda\gamma^{2}
+2(Mc_{0}\eta^{2}+M+1)\Bigg)\left(2(1+\delta)M^{2}\right)
\\
&\qquad\qquad\qquad
+\left(2M^{2}+\gamma^{2}\lambda+\frac{3}{2}\gamma^{2}(1+\gamma)\right)\Bigg]\mathbb{E}\Vert\hat{X}_{k}\Vert^{2}
\\
&\qquad
+c_{0}\eta
\Bigg((1+\gamma)\left(1+c_{0}\eta^{2}(1+\gamma)+\frac{3}{2}\eta^{4}\right)
+1+(c_{0}+1)\eta^{2}
\\
&\qquad\qquad\qquad
+\lambda\gamma^{2}
+2(Mc_{0}\eta^{2}+M+1)\Bigg)(1+\delta)B^{2}
+c_{0}B^{2}\eta
\\
&\qquad
+\frac{1}{2}\gamma^{3}\beta^{-1}c_{22}\eta^{2}
+\gamma^{2}\beta^{-1}c_{12}\eta
+M\gamma\beta^{-1}c_{22}\eta^{2},
\end{align*}

By applying the assumption $\eta\leq 1$, we have
\begin{align*}
\frac{\hat{L}_2(k+1) - \hat{L}_2(k)}{\eta} 
&\leq
\gamma ((d+A)/\beta- \lambda \hat{L}_2(k)) + (K_{1}\hat{L}_{2}(k)+K_{2})\cdot\eta
\\
&\qquad
+\eta Q_{1}\mathbb{E}\Vert\hat{V}_{k}\Vert^{2}
+\eta Q_{2}\mathbb{E}\Vert\hat{X}_{k}\Vert^{2}
+\eta Q_{3},
\end{align*}
where the constants $Q_1, Q_2, Q_3$ are given in \eqref{def-Q1}--\eqref{def-Q3}.
Let us recall that for $\lambda\leq\frac{1}{4}$,
\begin{equation*}
\mathcal{V}(x,v)
\geq\max\left\{\frac{1}{8}(1-2\lambda)\beta\gamma^{2}\Vert x\Vert^{2},
\frac{\beta}{4}(1-2\lambda)\Vert v\Vert^{2}\right\}.
\end{equation*}
Thus, we have
\begin{align*}
\frac{\hat{L}_2(k+1) - \hat{L}_2(k)}{\eta} 
&\leq
\gamma ((d+A)/\beta- \lambda \hat{L}_2(k)) + (K_{1}\hat{L}_{2}(k)+K_{2})\cdot\eta
\\
&\qquad
+\eta\left(Q_{1}\frac{4}{1-2\lambda}+Q_{2}\frac{8}{(1-2\lambda)\gamma^{2}}\right)\hat{L}_{2}(k)
+\eta Q_{3},
\end{align*}
Therefore, for 
\begin{equation}
0<\eta\leq\min\left\{\frac{\gamma}{\hat{K}_{2}}(d/\beta+A/\beta),\frac{\gamma\lambda}{2\hat{K}_{1}}\right\},
\end{equation}
where $\hat{K}_{1}:=K_{1}+\frac{4Q_{1}}{1-2\lambda}+\frac{8Q_{2}}{(1-2\lambda)\gamma^{2}},$ and $\hat{K}_{2}:=K_{2}+Q_{3},$
we get
\begin{equation*}
(\hat{L}_{2}(k+1)-\hat{L}_{2}(k))/\eta
\leq
2\gamma(d/\beta+A/\beta)-\frac{1}{2}\gamma\lambda \hat{L}_{2}(k).
\end{equation*}
This implies $\hat{L}_2(k+1) \le \rho \hat{L}_2(k) +K$, where
$\rho:= 1- \eta \gamma \lambda/2\in[0,1)$, where we used
the assumption $\eta\leq\frac{2}{\gamma\lambda}$, 
and $K:=2\eta\gamma(d/\beta+A/\beta).$
It follows that
\begin{align*}
\hat{L}_{2}(k)\leq \hat{L}_{2}(0)+\frac{K}{1-\rho}
=\mathbb{E}_{\mathbf{z}}\left[\mathcal{V}(\hat{X}_{0},\hat{V}_{0})/\beta\right]
+\frac{4(d/\beta+A/\beta)}{\lambda}.
\end{align*}
The uniform $L^2$ bounds then readily follow.

\subsection{Proof of Lemma~\ref{diff-app-2}}
We follow similar steps as in the proof of Lemma 7 in \cite{Raginsky}.
Recall that with the same initialization, the SGHMC2 iterates $(\hat{X}_{k}, \hat{V}_{k})$ has the same distribution as $(\hat{X}(k \eta), \hat{V}(k \eta))$ where $(\hat{X}(\cdot), \hat{V}(\cdot))$ is a continuous-time process satisfying
\begin{align} 
&d\hat{V}(t)=-\gamma\hat{V}(t)dt- g(\hat{X}(\lfloor t/\eta\rfloor\eta),U_{\mathbf{z}}(t))dt
+\sqrt{2\gamma \beta^{-1}}dB(t), 
\\
&d\hat{X}(t)=\hat{V}(t)dt, 
\end{align}
Let $\mathbb{P}$ be the probability measure associated with the underdamped Langevin diffusion $(X(t),V(t))$ in \eqref{eq:VL}--\eqref{eq:XL}
and $\hat{\mathbb{P}}$ be the probability measure associated with the $(\hat{X}(t),\hat{V}(t))$ process.
Let $\mathcal{F}_{t}$ be the natural filtration up to time $t$.
Then, the Radon-Nikodym derivative of $\mathbb{P}$ w.r.t. $\hat{\mathbb{P}}$
is given by the Girsanov theorem (see e.g. Section 7.6 in \cite{liptser2013statistics}):
\begin{equation*}
\frac{d\mathbb{P}}{d\hat{\mathbb{P}}}
\bigg|_{\mathcal{F}_{t}}
=e^{-\sqrt{\frac{\beta}{2\gamma}}\int_{0}^{t}(\nabla F_{\mathbf{z}}(\hat{X}(s))-g(\hat{X}(\lfloor s/\eta\rfloor\eta),U_{\mathbf{z}}(s)))\cdot dB(s)
-\frac{\beta}{4\gamma}\int_{0}^{t}\Vert\nabla F_{\mathbf{z}}(\hat{X}(s))-g(\hat{X}(\lfloor s/\eta\rfloor\eta),U_{\mathbf{z}}(s))\Vert^{2}ds}.
\end{equation*}
Then by writing $\mathbb{P}_{t}$ and $\hat{\mathbb{P}}_{t}$
as the probability measures $\mathbb{P}$ and $\hat{\mathbb{P}}$ conditional on the filtration $\mathcal{F}_{t}$, 
\begin{align*}
D(\hat{\mathbb{P}}_{t}\Vert\mathbb{P}_{t})
&:=-\int d\hat{\mathbb{P}}_{t}\log\frac{d\mathbb{P}_{t}}{d\hat{\mathbb{P}}_{t}}
\\
&=\frac{\beta}{4\gamma}\int_{0}^{t}\mathbb{E}_{\mathbf{z}}\left\Vert\nabla F_{\mathbf{z}}(\hat{X}(s))-g(\hat{X}(\lfloor s/\eta\rfloor\eta),U_{\mathbf{z}}(s))\right\Vert^{2}ds.
\end{align*}
Then, we get
\begin{align}
D(\hat{\mathbb{P}}_{k\eta}\Vert\mathbb{P}_{k\eta})
&=\frac{\beta}{4\gamma}\sum_{j=0}^{k-1}
\int_{j\eta}^{(j+1)\eta}\mathbb{E}_{\mathbf{z}}\left\Vert\nabla F_{\mathbf{z}}(\hat{X}(s))-g(\hat{X}(\lfloor s/\eta\rfloor\eta),U_{\mathbf{z}}(s))\right\Vert^{2}ds
\nonumber \\
&\leq
\frac{\beta}{2\gamma}\sum_{j=0}^{k-1}
\int_{j\eta}^{(j+1)\eta}\mathbb{E}_{\mathbf{z}}\left\Vert\nabla F_{\mathbf{z}}(\hat{X}(s))-\nabla F_{\mathbf{z}}(\hat{X}(\lfloor s/\eta\rfloor\eta))\right\Vert^{2}ds
\nonumber \\
&\qquad\qquad
+\frac{\beta}{2\gamma}\sum_{j=0}^{k-1}
\int_{j\eta}^{(j+1)\eta}\mathbb{E}_{\mathbf{z}}\left\Vert\nabla F_{\mathbf{z}}(\hat{X}(\lfloor s/\eta\rfloor\eta))
-g(\hat{X}(\lfloor s/\eta\rfloor\eta),U_{\mathbf{z}}(s))\right\Vert^{2}ds. \label{eq:rel-entropy-2}
\end{align}
We first bound the first term in \eqref{eq:rel-entropy-2}. 
Before we proceed, let us notice that for any $k\eta\leq s<(k+1)\eta$,
\begin{equation}
\hat{X}(s)=\hat{X}_{k}+\psi_{1}(s-k\eta)\hat{V}_{k}-\psi_{2}(s-k\eta)g(\hat{X}_{k},U_{\mathbf{z},k})
+\sqrt{2\gamma\beta^{-1}}\xi'_{k+1,s-k\eta},
\end{equation}
in distribution, where $\xi'_{k+1,s-k\eta}$
is centered Gaussian independent of $\mathcal{F}_{k}$ and 
$\mathbb{E}\Vert\xi'_{k+1,s-k\eta}\Vert^{2}=C_{22}(s-k\eta)\leq\frac{d}{3}(s-k\eta)^{3}\leq\frac{d}{3}\eta^{3}$.
Moreover, $\psi_{1}(s-k\eta)=\int_{0}^{s-k\eta}e^{-\gamma t}dt\leq(s-k\eta)\leq\eta$,
and $\psi_{2}(s-k\eta)=\int_{0}^{s-k\eta}\psi_{1}(t)dt\leq\int_{0}^{s-k\eta}tdt\leq\eta^{2}$.
Therefore, we can compute that
\begin{align*}
&\frac{\beta}{2\gamma}\sum_{j=0}^{k-1}
\int_{j\eta}^{(j+1)\eta}\mathbb{E}_{\mathbf{z}}\left\Vert\nabla F_{\mathbf{z}}(\hat{X}(s))-\nabla F_{\mathbf{z}}(\hat{X}(\lfloor s/\eta\rfloor\eta))\right\Vert^{2}ds
\\
&
\leq
\frac{\beta M^{2}}{2\gamma}\sum_{j=0}^{k-1}
\int_{j\eta}^{(j+1)\eta}\mathbb{E}_{\mathbf{z}}\left\Vert\hat{X}(s)-\hat{X}(\lfloor s/\eta\rfloor\eta)\right\Vert^{2}ds
\\
&
=\frac{\beta M^{2}}{2\gamma}\sum_{j=0}^{k-1}
\int_{j\eta}^{(j+1)\eta}\mathbb{E}_{\mathbf{z}}\left\Vert\psi_{1}(s-j\eta)\hat{V}_{j}-\psi_{2}(s-j\eta)g(\hat{X}_{j},U_{\mathbf{z},j})
+\sqrt{2\gamma\beta^{-1}}\xi'_{j+1,s-j\eta}\right\Vert^{2}ds
\\
&\leq
\frac{3\beta M^{2}}{2\gamma}\sum_{j=0}^{k-1}
\int_{j\eta}^{(j+1)\eta}\bigg(\mathbb{E}_{\mathbf{z}}\left\Vert\psi_{1}(s-j\eta)\hat{V}_{j}\right\Vert^{2}
+\mathbb{E}_{\mathbf{z}}\left\Vert\psi_{2}(s-j\eta)g(\hat{X}_{j},U_{\mathbf{z},j})\right\Vert^{2}
\\
&\qquad\qquad\qquad
+\mathbb{E}_{\mathbf{z}}\left\Vert\sqrt{2\gamma\beta^{-1}}\xi'_{j+1,s-j\eta}\right\Vert^{2}\bigg)ds
\\
&\leq
\frac{3\beta M^{2}}{2\gamma}(k\eta)
\bigg(\eta^{2}\sup_{j\geq 0}\mathbb{E}_{\mathbf{z}}\Vert\hat{V}_{j}\Vert^{2}
+\eta^{4}\left(2(1+\delta)M^{2}\sup_{j\geq 0}\mathbb{E}\Vert\hat{X}_{j}\Vert^{2}
+2(1+\delta)B^{2}\right)
+\frac{d\eta^{3}}{3}2\gamma\beta^{-1}\bigg)
\\
&\leq
\frac{3\beta M^{2}}{2\gamma}(k\eta)\eta^{2}
\bigg(C_{v}^{d}+\left(2(1+\delta)M^{2}C_{x}^{d}+2(1+\delta)B^{2}\right)+\frac{2d\gamma\beta^{-1}}{3}\bigg),
\end{align*}
where we used \eqref{g:estimate}, the assumption $\eta\leq 1$
and Lemma~\ref{lem:L2bound:Jordan}.

We can also bound the second term in \eqref{eq:rel-entropy-2}:
\begin{align*}
&\frac{\beta}{2\gamma}\sum_{j=0}^{k-1}
\int_{j\eta}^{(j+1)\eta}\mathbb{E}_{\mathbf{z}}\left\Vert\nabla F_{\mathbf{z}}(\hat{X}(\lfloor s/\eta\rfloor\eta))-g(\hat{X}(\lfloor s/\eta\rfloor\eta),U_{\mathbf{z}}(s))\right\Vert^{2}ds
\\
&=\frac{\beta}{2\gamma}\eta\sum_{j=0}^{k-1}\mathbb{E}_{\mathbf{z}}\left\Vert\nabla F_{\mathbf{z}}(\hat{X}_{j})-g(\hat{X}_{j},U_{\mathbf{z},j})\right\Vert^{2}
\\
&\leq\frac{\beta}{2\gamma}\eta\delta
\sum_{j=0}^{k-1}2(M^{2}\mathbb{E}_{\mathbf{z}}\Vert\hat{X}_{j}\Vert^{2}+B^{2})
\\
&\leq\left(M^{2}C_{x}^{d}+B^{2}\right)\frac{\beta}{\gamma}k\eta\delta,
\end{align*}
where the first inequality follows from part (iv) of Assumption~\ref{assumptions},
and we also used Lemma~\ref{lem:L2bound:Jordan}.
Hence, we conclude that
\begin{align}
D(\hat{\mu}_{\bz,k}\Vert\nu_{\bz,k\eta})
&\leq
\frac{3\beta M^{2}}{2\gamma}(k\eta)\eta^{2}
\bigg(C_{v}^{d}+\left(2(1+\delta)M^{2}C_{x}^{d}+2(1+\delta)B^{2}\right)+\frac{2d\gamma\beta^{-1}}{3}\bigg)
\nonumber \\
&\qquad\qquad\qquad
+\left(M^{2}C_{x}^{d}+B^{2}\right)\frac{\beta}{\gamma}k\eta\delta.\label{Destimate}
\end{align} 
To complete the proof, we can follow similar steps as in the proof of Lemma~\ref{diff-app}. By using the estimate in \eqref{Destimate}, the result from \cite{BV}, and the exponential integrability of the measure $\nu_{\bz,k\eta}$ in Lemma~\ref{lem:exp-integrable}, we can infer that 
\begin{align*}
&D(\hat{\mu}_{\bz,k}\Vert\nu_{\bz,k\eta})+\sqrt{D(\hat{\mu}_{\bz,k}\Vert\nu_{\bz,k\eta})}
\\
&\leq
\Bigg(\frac{3\beta M^{2}}{2\gamma}
\bigg(C_{v}^{d}+\left(2(1+\delta)M^{2}C_{x}^{d}+2(1+\delta)B^{2}\right)+\frac{2d\gamma\beta^{-1}}{3}\bigg)
\\
&\qquad\qquad
+\sqrt{\frac{3\beta M^{2}}{2\gamma}
\bigg(C_{v}^{d}+\left(2(1+\delta)M^{2}C_{x}^{d}+2(1+\delta)B^{2}\right)+\frac{2d\gamma\beta^{-1}}{3}\bigg)}\Bigg)k\eta^{2}
\\
&\qquad\qquad\qquad
+\left(\left(M^{2}C_{x}^{d}+B^{2}\right)\frac{\beta}{\gamma}+\sqrt{\left(M^{2}C_{x}^{d}+B^{2}\right)\frac{\beta}{\gamma}}\right)k\eta\sqrt{\delta},
\end{align*}
and 
\begin{align*}
\mathcal{W}_{2}^{2}(\hat{\mu}_{\bz,k},\nu_{\bz,k\eta})
\nonumber
&\leq 
\frac{8}{\alpha_{0}}\left(\frac{3}{2}+\log\left(\int_{\mathbb{R}^{2d}}e^{\frac{1}{4}\alpha\mathcal{V}(x,v)}\mu_{0}(dx,dv)
+\frac{1}{4}e^{\frac{\alpha(d+A)}{3\lambda}}\alpha\gamma(d+A)k\eta\right)\right)
\nonumber
\\
&\qquad\qquad\qquad
\cdot\left[D(\hat{\mu}_{\bz,k}\Vert\nu_{\bz,k\eta})+\sqrt{D(\hat{\mu}_{\bz,k}\Vert\nu_{\bz,k\eta})}\right],
\end{align*}
which together implies that
\begin{equation*}
\mathcal{W}_{2}^{2}(\hat{\mu}_{\bz,k},\nu_{\bz,k\eta})
\leq(C_{0}^{2}\sqrt{\delta}+\hat{C}_{1}^{2}\eta)(k\eta)\log(k\eta),
\end{equation*}
where $C_{0}$ and $\hat{C}_{1}$ are defined in \eqref{def-hat-Czero} and \eqref{def-hat-Jone-2}. The result then follows from the fact that $\sqrt{x+y} \le \sqrt{x} + \sqrt{y}$ for non-negative real numbers $x$ and $y$.

\section{Supporting Lemmas}\label{sec:support}
In this section, we present several supporting lemmas from the existing literature. These lemmas are used in our proofs, so we include them here for the sake of completeness. 
The first lemma shows that $f$ admits lower and upper bounds that are quadratic functions.
\begin{lemma}[See {\cite[Lemma 2]{Raginsky}}]\label{lem:gradient-bound}
If parts $(i)$ and $(ii)$ of Assumption~\ref{assumptions}  hold, then for all $x\in\mathbb{R}^{d}$ and $z$,
\begin{equation*}
\Vert\nabla f(x,z)\Vert
\leq M\Vert x\Vert+B,
\end{equation*}
and
\begin{equation*}
\frac{m}{3}\Vert x\Vert^{2}
-\frac{b}{2}\log 3
\leq f(x,z)
\leq\frac{M}{2}\Vert x\Vert^{2}+B\Vert x\Vert+ A_{0}.
\end{equation*}
\end{lemma}

The next lemma shows a 2-Wasserstein continuity result for functions of quadratic growth. This lemma was also used in \cite{Raginsky} to study the SGLD dynamics.
\begin{lemma}[See \cite{polyanskiy2016wasserstein}] \label{lem:expec-wass}
Let $\mu,\nu$ be two probability measures on $\mathbb{R}^{2d}$
with finite second moments, and let $G:\mathbb{R}^{2d}\rightarrow\mathbb{R}$
be a $C^{1}$ function obeying
\begin{equation*}
\Vert\nabla G(w)\Vert
\leq c_{1}\Vert w\Vert+c_{2},
\end{equation*}
for some constants $c_{1}>0$ and $c_{2}\geq 0$.
Then,
\begin{equation*}
\left|
\int_{\mathbb{R}^{2d}}Gd\mu-\int_{\mathbb{R}^{2d}}Gd\nu\right|
\leq(c_{1}\sigma+c_{2})\mathcal{W}_{2}(\mu,\nu),
\end{equation*}
where
\begin{equation*}
\sigma^{2}=\max\left\{\int_{\mathbb{R}^{2d}}\Vert w\Vert^{2}\mu(dw),
\int_{\mathbb{R}^{2d}}\Vert w\Vert^{2}\nu(dw)\right\}.
\end{equation*}
\end{lemma}

The next lemma shows a uniform stability of $\pi_{\mathbf{z}}$. Note that the $x-$marginal of $\pi_{\mathbf{z}}(dx,dv)$ for the underdamped diffusion is the same as the stationary distribution for the overdamped diffusion studied in \cite{Raginsky}.
For two $n-$tuples $\bz =(z_1, \ldots, z_n), \overline{\bz}= (\overline{z}_1, \ldots, \overline{z}_n) \in\mathcal{Z}^n$, we say $\bz$ and $\overline{\bz}$ differ only in a single coordinate if card$|\{i: z_i \neq \overline{z}_i\}| =1$.

\begin{lemma}[Proposition 12, \cite{Raginsky}] \label{lem:stability}
For any two ${\mathbf{z}},\overline{\mathbf{z}} \in\mathcal{Z}^n$ that differ only in a single coordinate,
\begin{equation*}
\sup_{z\in\mathcal{Z}}
\left|\int_{\mathbb{R}^{2d}}f(x,z)\pi_{\mathbf{z}}(dx,dv)
-\int_{\mathbb{R}^{2d}}f(x,z)\pi_{\overline{\mathbf{z}}}(dx,dv)\right|
\leq\frac{4\beta c_{LS}}{n}
\left(\frac{M^{2}}{m}(b+d/\beta)+B^{2}\right),
\end{equation*}
where
\begin{equation*}
c_{LS}
\leq\frac{2m^{2}+8M^{2}}{m^{2}M\beta}
+\frac{1}{\lambda_{\ast}}\left(\frac{6M(d+\beta)}{m}+2\right),
\end{equation*}
where $\lambda_{\ast}$ is the uniform spectral gap for overdamped Langevin dynamics:
\begin{equation*}
\lambda_{\ast}=
\inf_{{\mathbf{z}}\in\mathcal{Z}^n}
\inf
\left\{
\frac{\beta^{-1}\int_{\mathbb{R}^{d}}\Vert\nabla g\Vert^{2}d\pi_{\mathbf{z}}}{\int_{\mathbb{R}^{d}}g^{2}d\pi_{\mathbf{z}}}:
g\in C^{1}(\mathbb{R}^{d})\cap L^{2}(\pi_{\mathbf{z}}),
g\neq 0,
\int_{\mathbb{R}^{d}}gd\pi_{\mathbf{z}}=0\right\}.
\end{equation*}
\end{lemma}

The next lemma show that for large values of $\beta$, the $x-$marginal of the stationary distribution $\pi_{\mathbf{z}}(dx,dv)$
is concentrated at the minimizer of $F_{\bz}$. Note in Proposition 11 of \cite{Raginsky}, they have
the assumption $\beta\geq 2/m$, which seems to be only used to derive their Lemma 4,
but not used in deriving their Proposition 11.

\begin{lemma}[Proposition 11, \cite{Raginsky}]\label{lem:thirdbound} It holds that
\begin{equation*}
\int_{\mathbb{R}^{2d}}F_{\mathbf{z}}(x)\pi_{\mathbf{z}}(dx,dv)-\min_{x\in\mathbb{R}^{d}}F_{\mathbf{z}}(x)
\leq\frac{d}{2 \beta}\log\left(\frac{eM}{m}\left(\frac{b\beta}{d}+1\right)\right).
\end{equation*}
\end{lemma}

\section{Proof of Proposition~\ref{prop:comparison}}\label{sec:comparison}

Let us first prove that $\lambda_{\ast}=\mathcal{O}(a^{-2})$.
We first recall that
$\lambda_{\ast}$ is the uniform spectral gap for overdamped Langevin dynamics:
\begin{equation*}
\lambda_{\ast}:=
\inf_{\mathbf{z}\in\mathcal{Z}^n}
\inf
\left\{
\frac{\beta^{-1}\int_{\mathbb{R}^{d}}\Vert\nabla g\Vert^{2}d\pi_{\mathbf{z}}}{\int_{\mathbb{R}^{d}}g^{2}d\pi_{\mathbf{z}}}:
g\in C^{1}(\mathbb{R}^{d})\cap L^{2}(\pi_{\mathbf{z}}),
g\neq 0,
\int_{\mathbb{R}^{d}}gd\pi_{\mathbf{z}}=0\right\}.
\end{equation*}

In particular, fix any $\mathbf{z}\in\mathcal{Z}^{n}$
so that for every $g\in C^{1}(\mathbb{R}^{d})\cap L^{2}(\pi_{\mathbf{z}})$,
such that $g\neq 0$,
and $\int_{\mathbb{R}^{d}}gd\pi_{\mathbf{z}}=0$, we have 
\begin{equation*}
\lambda_{\ast}
\leq\frac{\beta^{-1}\int_{\mathbb{R}^{d}}\Vert\nabla g\Vert^{2}e^{-\beta F_{\mathbf{z}}(x)}dx}
{\int_{\mathbb{R}^{d}}g^{2}e^{-\beta F_{\mathbf{z}}(x)}dx}.
\end{equation*}

It follows from Lemma~\ref{lem:gradient-bound} that
\begin{equation}\label{eqn:Fz:bounds}
\frac{m}{3}\Vert x\Vert^{2}-\frac{b}{2}\log 3
\leq F_{\mathbf{z}}(x)
\leq\frac{M^{2}}{2}\Vert x\Vert^{2}+B\Vert x\Vert+A_{0},
\end{equation}
with $m=m_{1}a^{-2}$, $M=M_{1}a^{-2}$, and $B=B_{1}a^{-1}$.

Next, let us take the test function
$g_{1}(x):=\Vert x\Vert^{2}$.
And we further define
\begin{equation}\label{def:c1:appendix}
c_{1}:=\int_{\mathbb{R}^{d}}g_{1}d\pi_{\mathbf{z}}
=\frac{\int_{\mathbb{R}^{d}}g_{1}(x)e^{-\beta F_{\mathbf{z}}(x)}dx}{\int_{\mathbb{R}^{d}}e^{-\beta F_{\mathbf{z}}(x)}dx},
\end{equation}
and we also define
\begin{equation*}
g(x):=g_{1}(x)-c_{1},
\end{equation*}
so that $g\in C^{1}(\mathbb{R}^{d})\cap L^{2}(\pi_{\mathbf{z}})$,
$g\neq 0$, and $\int_{\mathbb{R}^{d}}gd\pi_{\mathbf{z}}=0$.
Moreover, we have
\begin{equation*}
\Vert\nabla g(x)\Vert=\Vert\nabla g_{1}(x)\Vert=2\Vert x\Vert,
\qquad
\text{and}
\qquad
g_{1}(ax)=a^{2}g_{1}(x)=a^{2}\Vert x\Vert^{2}.
\end{equation*}

Next, by the definition of $c_{1}$ in \eqref{def:c1:appendix} and the bounds in \eqref{eqn:Fz:bounds}, we get
\begin{align*}
&c_{1}
\geq\frac{\int_{\mathbb{R}^{d}}\Vert x\Vert^{2}e^{-\beta(\frac{M^{2}}{2}\Vert x\Vert^{2}+B\Vert x\Vert+A_{0})}dx}{\int_{\mathbb{R}^{d}}e^{-\beta(\frac{m}{3}\Vert x\Vert^{2}-\frac{b}{2}\log 3)}dx}
=
\frac{\int_{\mathbb{R}^{d}}\Vert ax\Vert^{2}e^{-\beta(\frac{M^{2}}{2}\Vert ax\Vert^{2}+B\Vert ax\Vert+A_{0})}dx}{\int_{\mathbb{R}^{d}}e^{-\beta(\frac{m}{3}\Vert ax\Vert^{2}-\frac{b}{2}\log 3)}dx}
=a^{2}\underline{c}_{1},
\\
&c_{1}
\leq\frac{\int_{\mathbb{R}^{d}}\Vert x\Vert^{2}e^{-\beta(\frac{m}{3}\Vert x\Vert^{2}-\frac{b}{2}\log 3)}dx}{\int_{\mathbb{R}^{d}}e^{-\beta(\frac{M^{2}}{2}\Vert x\Vert^{2}+B\Vert x\Vert+A_{0})}dx}
=
\frac{\int_{\mathbb{R}^{d}}\Vert ax\Vert^{2}e^{-\beta(\frac{m}{3}\Vert ax\Vert^{2}-\frac{b}{2}\log 3)}dx}{\int_{\mathbb{R}^{d}}e^{-\beta(\frac{M^{2}}{2}\Vert ax\Vert^{2}+B\Vert ax\Vert+A_{0})}dx}
=a^{2}\overline{c}_{1},
\end{align*}
where
\begin{align*}
&\underline{c}_{1}:=\frac{\int_{\mathbb{R}^{d}}\Vert x\Vert^{2}e^{-\beta(\frac{M_{1}^{2}}{2}\Vert x\Vert^{2}+B_{1}\Vert x\Vert+A_{0})}dx}{\int_{\mathbb{R}^{d}}e^{-\beta(\frac{m_{1}}{3}\Vert x\Vert^{2}-\frac{b}{2}\log 3)}dx},
\\
&\overline{c}_{1}:=\frac{\int_{\mathbb{R}^{d}}\Vert x\Vert^{2}e^{-\beta(\frac{m_{1}}{3}\Vert x\Vert^{2}-\frac{b}{2}\log 3)}dx}{\int_{\mathbb{R}^{d}}e^{-\beta(\frac{M_{1}^{2}}{2}\Vert x\Vert^{2}+B_{1}\Vert x\Vert+A_{0})}dx}.
\end{align*}

Hence, we have
\begin{align*}
\lambda_{\ast}
&\leq\frac{\beta^{-1}\int_{\mathbb{R}^{d}}\Vert\nabla g(x)\Vert^{2}e^{-\beta(\frac{m}{3}\Vert x\Vert^{2}-\frac{b}{2}\log 3)}dx}
{\int_{\mathbb{R}^{d}}g(x)^{2}e^{-\beta(\frac{M^{2}}{2}\Vert x\Vert^{2}+B\Vert x\Vert+A_{0})}dx}
\\
&=\frac{\beta^{-1}\int_{\mathbb{R}^{d}}4\Vert x\Vert^{2}e^{-\beta(\frac{m}{3}\Vert x\Vert^{2}-\frac{b}{2}\log 3)}dx}
{\int_{\mathbb{R}^{d}}(g_{1}(x)-c_{1})^{2}e^{-\beta(\frac{M^{2}}{2}\Vert x\Vert^{2}+B\Vert x\Vert+A_{0})}dx}
\\
&=\frac{\beta^{-1}\int_{\mathbb{R}^{d}}4\Vert ax\Vert^{2}e^{-\beta(\frac{m}{3}\Vert ax\Vert^{2}-\frac{b}{2}\log 3)}dx}
{\int_{\mathbb{R}^{d}}(g_{1}(ax)-c_{1})^{2}e^{-\beta(\frac{M^{2}}{2}\Vert ax\Vert^{2}+B\Vert ax\Vert+A_{0})}dx}
\\
&\leq
\frac{\beta^{-1}\int_{\mathbb{R}^{d}}4\Vert ax\Vert^{2}e^{-\beta(\frac{m}{3}\Vert ax\Vert^{2}-\frac{b}{2}\log 3)}dx}
{\min_{a^{2}\underline{c}_{1}\leq\tilde{c}\leq a^{2}\overline{c}_{1}}
\int_{\mathbb{R}^{d}}(a^{2}\Vert x\Vert^{2}-\tilde{c})^{2}e^{-\beta(\frac{M^{2}}{2}\Vert ax\Vert^{2}+B\Vert ax\Vert+A_{0})}dx}
\\
&=
a^{-2}\frac{\beta^{-1}\int_{\mathbb{R}^{d}}4\Vert x\Vert^{2}e^{-\beta(\frac{m_{1}}{3}\Vert x\Vert^{2}-\frac{b}{2}\log 3)}dx}
{\min_{\underline{c}_{1}\leq c\leq\overline{c}_{1}}\int_{\mathbb{R}^{d}}(\Vert x\Vert^{2}-c)^{2}e^{-\beta(\frac{M_{1}^{2}}{2}\Vert x\Vert^{2}+B_{1}\Vert x\Vert+A_{0})}dx},
\end{align*}
where we used $m=m_{1}a^{-2}$, $M=M_{1}a^{-2}$, $B=B_{1}a^{-1}$
and $g_{1}(ax)=a^{2}g_{1}(x)=a^{2}\Vert x\Vert^{2}$.
Hence, we conclude that $\lambda_{\ast}=\mathcal{O}(a^{-2})$.

Next, let us prove that
$\mu_{\ast}=\Theta(a^{-1})$.
We recall that $\mu_*$ the convergence rate for underdamped Langevin dynamics is given by:
\begin{equation*}
\mu_* =\frac{\gamma}{768}\min\left\{\lambda M \gamma^{-2},\Lambda^{1/2}e^{-\Lambda}M\gamma^{-2},\Lambda^{1/2}e^{-\Lambda}\right\},
\end{equation*}
where
\begin{equation*}
\Lambda=\frac{12}{5}(1+2\alpha_{1}+2\alpha_{1}^{2})(d+A)M\gamma^{-2}\lambda^{-1}(1-2\lambda)^{-1},
\qquad
\alpha_{1}=(1+\Lambda^{-1})M\gamma^{-2},
\end{equation*}
where $\lambda,A$ come from the drift condition \eqref{eq:drift},
and from \cite{GGZ2}, we can take
\begin{equation}\label{take:A}
\lambda=\frac{1}{2}\min\left(\frac{1}{4},\frac{m}{M+\gamma^{2}/2}\right),
\qquad
A=\frac{\beta}{2}\frac{m}{M+\frac{1}{2}\gamma^{2}}\left(\frac{B^{2}}{2M+\gamma^{2}}+\frac{b}{m}\left(M+\frac{1}{2}\gamma^{2}\right)+A_{0}\right).
\end{equation}
Note that $\mu_{\ast}$ depends on the objective function $F_{\mathbf{z}}$ only via
the parameters from its properties, which is independent of $\mathbf{z}$.
Recall that $m=m_{1}a^{-2}$, $M=M_{1}a^{-2}$, $B=B_{1}a^{-1}$.
We define $\gamma=:\gamma_{1}a^{-1}$ so that $\gamma_{1}$ is independent of $a$ and
\begin{equation}\label{eqn:mu:a}
\mu_* =a^{-1}\frac{\gamma_{1}}{768}\min\left\{\lambda M_{1}\gamma_{1}^{-2},\Lambda^{1/2}e^{-\Lambda}M_{1}\gamma_{1}^{-2},\Lambda^{1/2}e^{-\Lambda}\right\},
\end{equation}
where we can check that $\lambda$, $\Lambda$ are independent of $a$. 
Then, we can see from \eqref{eqn:mu:a} that $\mu_{\ast}$ is linear in $a^{-1}$
so that we have $\mu_{\ast}=\Theta(a^{-1})$. The proof is complete.

\section{Explicit dependence of constants on key parameters}\label{app-parameter-dependence}

In this section we provide explicit dependence of constants on parameters $\beta, d, \mu_{\ast}, \lambda_{\ast}$ and $n$, which is used in Section~\ref{sec-compare}. To simplify the presentation, we use the notation $\tilde{\bigO}$, $\tilde{\Theta}$ to hide factors that depend on other parameters.

We recall the constants from Table~\ref{table_constants}.
It is easy to see that
\begin{equation*}
A=\tilde{\Theta} (\beta),
\qquad
\alpha_{1}=\tilde{\Theta}(1),
\qquad
\Lambda=\tilde{\Theta}(d+\beta),
\end{equation*}
where we take $A$ as in \eqref{take:A} and
\begin{equation}\label{eq:mu-ast}
\mu_{\ast}=\tilde{\Theta} \left(\sqrt{d+\beta}e^{-\Lambda}\right)=\tilde{\Theta}\left(\sqrt{d+\beta}e^{-\tilde{\Theta}(d+\beta)}\right).
\end{equation}
Since $\varepsilon_{1}=\tilde{\bigO}(\mu_{\ast}/(d+\beta))$,
and $\mu_{\ast}$ is exponentially small in $\beta+d$,
we get that
\begin{equation*}
\overline{\mathcal{H}}_{\rho}(\mu_{0})
=\tilde{\bigO}(R_{1})=(1+d/\beta)^{1/2}.
\end{equation*}
In addition, in view of \eqref{eq:mu-ast}, it follows that
\begin{equation*}
C=\tilde{\bigO} \left(e^{\Lambda/2}(d+\beta)^{1/2}\beta^{-1/2}\mu_{\ast}^{-1/2}\right)
=\tilde{\bigO} \left(\frac{(d+\beta)^{3/4}\beta^{-1/2}}{\mu_{\ast}}\right).
\end{equation*}
The structure of the initial distribution $\mu_{0}(dx,dv)$ would affect
the overall dependence on $\beta,d$.
Since we assumed in Section~\ref{sec-compare} that $\mu_{0}(dx,dv)$ is supported on a Euclidean ball with radius being a universal constant, then the Lyapunov function
$\mathcal{V}$ in \eqref{eq:lyapunov} is linear in $\beta$. We can then obtain
\begin{equation*}
\int_{\mathbb{R}^{2d}}\mathcal{V}(x,v)\mu_{0}(dx,dv)
=\tilde{\bigO}(\beta),
\qquad
\int_{\mathbb{R}^{2d}}e^{\alpha\mathcal{V}(x,v)}\mu_{0}(dx,dv)=e^{\tilde{\bigO}(\beta)},
\end{equation*}
It follows that
\begin{equation*}
C_{x}^{d}=\tilde{\bigO}\left((\beta+d)/\beta\right),
\qquad
C_{v}^{d}=\tilde{\bigO}\left((\beta+d)/\beta\right),
\qquad
\sigma=\tilde{\bigO}\left(\sqrt{(\beta+d)/\beta}\right).
\end{equation*}


Next, we have $\alpha_{0}=\tilde{\bigO} (\beta)$ and $\alpha=\tilde{\bigO}(1)$, and
\begin{align*}
&\hat{\gamma}=\tilde{\bigO} (\sqrt{(\beta+d)/\beta}),
\\
&C_{0}= \tilde{\bigO} \left((d+\beta)/\sqrt{\beta}\right),
\qquad
C_{1}=\tilde{\bigO} \left((d+\beta)/\sqrt{\beta}\right),
\qquad
C_{2}=\tilde{\bigO} \left(\sqrt{(d+\beta)/\beta}\right).
\end{align*}
Moreover, by the definition of $\hat{C}_{1}$ in \eqref{def-hat-Jone-2}, we get
\begin{equation*}
\hat{C}_{1}=\tilde{\bigO} \left((d+\beta)/\sqrt{\beta}\right).
\end{equation*}
Hence, from \eqref{def-J0-bar}, we obtain
\begin{equation*}
\overline{\mathcal{J}}_{0}(\varepsilon)
=\tilde{\bigO}\left(\frac{(d+\beta)^{3/2}}{\mu_{\ast} \beta^{5/4}} \varepsilon\right),
\end{equation*}
and from \eqref{eq:J1}, 
we get
\begin{equation*}
\mathcal{J}_{1}(\varepsilon)
=\tilde{\bigO}\left(
\frac{(d+\beta)^{3/2}}{\beta(\mu_{\ast})^{3/2}}
\left((\log(1/\varepsilon))^{3/2}\delta^{1/4}
+\varepsilon\right)
\sqrt{\log(\mu_{\ast}^{-1}\log(\varepsilon^{-1}))}
+\frac{d+\beta}{\beta}\frac{\varepsilon^{2}}{\mu_{\ast}(\log(1/\varepsilon))^{2}}\right)\,.
\end{equation*}
Moreover, from \eqref{eq:J1:Jordan}, we get
\begin{equation*}
\hat{\mathcal{J}}_1(\varepsilon)=\tilde{\bigO}\left(
\frac{(d+\beta)^{3/2}}{\beta\sqrt{\mu_{\ast}}}
\left(\sqrt{\log(1/\varepsilon)}\delta^{1/4}
+\varepsilon\right)
\sqrt{\log(\mu_{\ast}^{-1}\log(\varepsilon^{-1}))}\right)\,. 
\end{equation*}
Finally, from \eqref{def-J2} and \eqref{def-J3}, we have
\begin{equation*}
\mathcal{J}_{2}= \tilde{\bigO} \left(\frac{d}{\beta}\log(\beta+1)\right),
\qquad\text{and}
\qquad
\mathcal{J}_{3}(n)=\tilde{\bigO}\left(\frac{1}{n}\frac{(\beta+d)^{2}}{\lambda_{\ast}}\right).
\end{equation*}

\section{Proof of Lemma~\ref{lem-example-classification} and Lemma~\ref{lem:2}}

\subsection{Proof of Lemma~\ref{lem-example-classification}}

Since the distribution of $A_{in}$ has compact support, we have $\|a_i\| \leq D$ for some $D>0$. Let $s_i := \langle a_i, x \rangle$. 
By taking the gradient of $f(x,z_i)$ with respect to $x$, we obtain 
\begin{equation} \nabla f(x,z_i) = - 2 \left(y_i - \sigma(s_i) \right) \sigma'(s_i) a_i + \lambda_r x.\label{eqn-grad-logistic}
\end{equation}
This implies
\begin{eqnarray}\langle \nabla f(x,z_i), x \rangle &=& - 2 \left(y_i - \sigma(s_i) \right) \sigma'(s_i) s_i + \lambda_r \|x\|^2 \\
&\geq& \lambda_r \|x\|^2 - 2 (1+ \|\sigma\|_\infty) \|\sigma'\|_\infty |s_i| \\
&\geq& \lambda_r \|x\|^2 - 2 (1+ \|\sigma\|_\infty) \|\sigma'\|_\infty D \|x\|\,,
\end{eqnarray}
where we used the triangle inequality and the Cauchy-Schwartz inequality. Then, it is straightforward to check that we obtain $\langle \nabla f(x,z_i), x \rangle \geq  m\|x\|^2 - b$
for
 $$ m = \lambda_r/2, \quad b = 8 (1+\|\sigma\|_\infty)^2 \|\sigma'\|_\infty^2 D^2 / \lambda_r\,,$$
and therefore part (iii) of Assumption~\ref{assumptions} holds. Also for any $z=(a,y)$, 
$|f(0,z)| = |(y - \sigma(0))^2| \leq A_0$ for $A_0 =   (1 + \|\sigma(0)\|)^2$. 
Similarly, 
$\|\nabla f(0,z)\| = \| -2 (y - \sigma(0)) \sigma'(0) a \| \leq B_1$ for $$B_1 := 2(1+|\sigma(0)|) |\sigma'(0)| D.$$ 
Therefore, part (i) of Assumption~\ref{assumptions} holds for any  $B \geq B_1$. We also have the Hessian matrix
$$ \nabla^2 f(x,z_i) =  2 [\sigma'(s_i)]^2 a_i a_i^T - 2 (y_i - \sigma(s_i))\sigma''(s_i)a_ia_i^T +\lambda_r I_d\,,$$
where $I_d$ is the $d\times d$ identity matrix. Hence,
$ \|\nabla^2 f(x,z_i)\| \leq M_1$
where $$M_1=2\|\sigma'\|_\infty^2 D^2 + 2\left(1+\|\sigma\|_\infty\right) \|\sigma''\|_\infty D^2 + \lambda_r.$$ 
Therefore, part (ii) of Assumption~\ref{assumptions} also holds for $M=M_1$.  In particular, $\nabla f(x,z_{j})$ is also i.i.d. and we have $ \mathbb{E}[ \nabla f(x,z_{j})]= \nabla F_{\mathbf{z}}(x)$ for any $x\in\mathbb{R}^d$. Furthermore, it follows from \eqref{eqn-grad-logistic} that
$$\| \nabla f(x,z_j)\| \leq B_2 + \lambda_r \|x\|, \quad \mbox{where} \quad B_2 := 2 \left(1 + \|\sigma\|_\infty \right) \|\sigma'\|_\infty D, $$
for any $z_j$. Therefore, if we let $u_j := \nabla f(x,z_j) - \nabla F_{\mathbf{z}}(x)$, then $u_j$ are i.i.d. with $\mathbb{E}[u_j] = 0$ and
\begin{align} 
\mathbb{E}\|u_j\|^2 &\leq 2 \mathbb{E}\left[\left\Vert \nabla f(x,z_j) \right\Vert^2\right] + 2 \mathbb{E}\left[\left\|\nabla F_{\mathbf{z}}(x)\right\|^2\right]\nonumber  \\
&\leq 4 (B_2 + \lambda_r \|x\|)^2 \nonumber
\\
&\leq 8 B_2^2 + 8\lambda_r^2\|x\|^2\,,
\label{eq-variance-bound}
\end{align}
where we used Cauchy-Schwarz inequality. This implies
 \beq \mathbb{E}\left[\left\Vert g(x,U_{\mathbf{z}})-\nabla F_{\mathbf{z}}(x)\right\Vert^{2}\right]
  =  \mathbb{E} \left\|\frac{1}{n_b}\sum_{j=1}^{n_b} u_j\right\|^2 =  \frac{1}{n_b^2}\sum_{j=1}^{n_b} \mathbb{E}\|u_j\|^2
  \leq 2\delta (M^2 \|x\|^2 + B^2) 
 \eeq
for any $\delta \in [\frac{1}{4n_b},1)$, $M \geq M_2 := 4\lambda_r$, $B \geq B_2 := 4M_2$ where we used \eqref{eq-variance-bound} and the fact that $u_j$ are i.i.d. with mean zero. If we choose, for instance, $M = M_1 + M_2$, $B=\max(B_1, B_2)=B_2 $; we observe that part (i) and (iv) of Assumptions~\ref{assumptions} hold.

\subsection{Proof of Lemma~\ref{lem:2}}

We set $r_i = y_i - \langle a_i, x \rangle$ and follow a similar approach to the proof of Lemma~\ref{lem-example-classification}.

We can compute that 
$$ \nabla f(x,z_i) = -\rho'(r_i) a_i + \lambda_r x.$$
This leads to
\beq 
\langle \nabla f(x,z_i), x \rangle = -\rho'(r_i) \langle a_i,x \rangle + \lambda_r \|x\|^2 
\geq \lambda_r \|x\|^2 - \|\rho'\| D \|x\|
\geq m\|x\|^2 - b\,,
\eeq
with $$ m=\lambda_r/2, \quad b = \frac{2 \|\rho'\|_\infty^2 D^2}{\lambda_r},$$
Therefore, part (iii) of Assumption~\ref{assumptions} holds.
We have also
$$|f(0,z_i)| \leq | \rho(y_i)|, \quad \nabla f(0,z_i) = -\|\rho'(y_i) a_i\| \leq \|\rho'\|_\infty D$$ 
for any $z_i$. Therefore, part (i) of Assumption~\ref{assumptions} holds with $A_0 = \|\rho\|_\infty$ and $B =  \|\rho'\|_\infty D$. Since
    $$ \nabla^2 f(0,z_i) = \rho''(r_i) a_i a_i^T + \lambda_r I_d\,,$$
where $I_d$ is the $d\times d$ identity matrix, we also have 
    $$ \left\|\nabla^2 f(0,z_i) \right\| \leq \|\rho''\|_\infty D^2 + \lambda_r. $$
Therefore, part (ii) of Assumption~\ref{assumptions} holds for any $M \geq   \|\rho''\|_\infty D^2 + \lambda_r$. We have also
$$ \|\nabla f(x,z_i)\|\leq \|\rho'\|_\infty D + \lambda_r \|x\|.$$
Therefore, if we let $v_j = \nabla f(x,z_j) - \nabla F_{\mathbf{z}}(x)$, then $v_j$ are i.i.d. with $\mathbb{E}[v_j] = 0$ and
\begin{align} 
\mathbb{E}\|v_j\|^2 &\leq 2 \mathbb{E}\left[\left\Vert \nabla f(x,z_j) \right\Vert^2\right] + 2 \mathbb{E}\left[\left\|\nabla F_{\mathbf{z}}(x)\right\|^2\right]\nonumber  
\\
&\leq 4 \left(\|\rho'\|_\infty D + \lambda_r \|x\|\right)^2 \nonumber
\\
&\leq 8 \left(\|\rho'\|_\infty^2 D^2 + 8\lambda_r^2\|x\|^2\right), 
 \label{eq-variance-bound-2}
\end{align}
where we used Cauchy-Schwarz inequality. This implies
 \beq \mathbb{E}\left[\left\Vert g(x,U_{\mathbf{z}})-\nabla F_{\mathbf{z}}(x)\right\Vert^{2}\right]
  =  \mathbb{E} \left\|\frac{1}{n_b}\sum_{j=1}^{n_b} v_j\right\|^2 =  \frac{1}{n_b^2}\sum_{j=1}^{n_b} \mathbb{E}\|v_j\|^2
  \leq 2\delta (M^2 \|x\|^2 + B^2) 
 \eeq
for any $\delta \in [\frac{1}{4n_b},1)$ and $M_2 \geq 4\lambda_r$ and $B \geq 4\|\rho'\|_\infty D$ where we used \eqref{eq-variance-bound-2} and the fact that $v_j$ are i.i.d. with mean zero. We conclude that Assumption~\ref{assumptions} work for $M = \|\rho''\|_\infty D^2 + 5\lambda_r$ and $B =  4\|\rho'\|_\infty D$.

\end{document}